\documentclass{elsarticle}

\usepackage{hyperref}
\usepackage{float}
\usepackage{graphicx}
\usepackage{subfigure}
\usepackage{amsmath,bm}
\usepackage{epstopdf} 
\usepackage{caption}
\usepackage{mathptmx}
\usepackage{verbatim}
\usepackage{indentfirst}
\usepackage[title]{appendix}
\usepackage{amssymb}
\usepackage{cases}
\usepackage{color}

\begin{document}
\newtheorem{theorem}{Theorem}[section]
\newtheorem{example}{Example}[section]
\newtheorem{remark}{Remark}[section]
\numberwithin{equation}{section}
\numberwithin{figure}{section}
\numberwithin{table}{section}

\captionsetup[figure]{font=small,labelfont={bf},labelformat={default},labelsep=period,name={Fig.}}
\captionsetup[table]{labelfont={bf},labelformat={default},labelsep=period,name={Table}}

\begin{frontmatter}

\title{A model-data asymptotic-preserving neural network method based on micro-macro decomposition for gray radiative transfer equations}

\author[10614]{Hongyan Li}
\ead{lihongyan142869@163.com}

\author[909]{Song Jiang}
\ead{jiang@iapcm.ac.cn}

\author[909]{Wenjun Sun}
\ead{sun$\_$wenjun@iapcm.ac.cn}

\author[10614]{Liwei Xu\corref{mycorrespondingauthor}}
\cortext[mycorrespondingauthor]{Corresponding author}
\ead{xul@uestc.edu.cn}

\author[106141]{Guanyu Zhou}
\ead{zhoug@uestc.edu.cn}

\address[10614]{School of Mathematical Sciences, University of Electronic Science and Technology of China, Chengdu Sichuan 611731, China}
\address[909]{Institute of Applied Physics and Computational Mathematics, Beijing 100094, China}
\address[106141]{Institute of Fundamental and Frontier Sciences, University of Electronic Science and Technology of China, Chengdu Sichuan 610054, China}

\begin{abstract}
We propose a model-data asymptotic-preserving neural network(MD-APNN) method to solve the nonlinear gray radiative transfer equations(GRTEs).
The system is challenging to be simulated with both the traditional numerical schemes and the vanilla physics-informed neural networks(PINNs)
due to the multiscale characteristics. Under the framework of PINNs, we employ a micro-macro decomposition technique to construct a new asymptotic-preserving(AP) loss function,
which includes the residual of the governing equations in the micro-macro coupled form, the initial and boundary conditions with additional diffusion limit information, the conservation laws, and a few labeled data. A convergence analysis is performed
for the proposed method, and a number of numerical examples are presented to illustrate the efficiency of MD-APNNs,
and particularly, the importance of the AP property in the neural networks for the diffusion dominating problems.
The numerical results indicate that MD-APNNs lead to a better performance than APNNs or pure data-driven networks in
the simulation of the nonlinear non-stationary GRTEs.
\end{abstract}


\end{frontmatter}


\section{Introduction}
The gray radiative transfer equations(GRTEs) describe the radiation transport of
photons and energy exchange with background materials, which find a wide range of applications in the fields of astrophysics, inertial/magnetic confinement fusion, high-temperature flow systems, and so on \cite{chandrasekhar2013radiative,howell2020thermal,pomraning2005equations}. 
As a type of kinetic model coupled with a nonlinear material thermal energy equation,
an accurate simulation of GRTEs is nontrivial to achieve because of the high dimensionality, strongly coupled nonlinearity,
and multi-scale features caused by different opacities of the background materials.

The stochastic and deterministic methods are two widely applied conventional numerical approaches for the GRTEs.
The former, for instance the Implicit Monte Carlo(IMC) method \cite{fleck1971implicit,gentile2001implicit,mcclarren2009modified,shi2020continuous,shi2018functional,densmore2011asymptotic}, has the advantage on dealing with the high dimensionality,
and however the convergence rate is low and suffers statistical noises.
While the deterministic schemes (e.g. finite difference/element/volume methods, discontinuous Galerkin method) with the asymptotic preserving(AP) technique are designed to capture the multi-scale features and treat the coupled nonlinearity.
However, for high-dimensional cases, the discretization of the deterministic methods yields a very large algebraic system to be solved, which demands huge computational costs.
The AP schemes are first proposed for the steady neutron transport problems \cite{jin1991discrete,jin1993fully,larsen1987asymptotic,morel1989asymptotic}. For the unsteady cases, the AP schemes are constructed based on a decomposition of the distribution into the equilibrium and disturbance parts \cite{jin2000uniformly,klar1998asymptotic}.
The method has been further developed and extended to the multi-scale kinetic equations\cite{jin1999efficient,jin2010asymptotic,lemou2008new}, and to the the GRTEs, for example the AP unified gas kinetic scheme(AP-UGKS)\cite{sun2015asymptotic,xu2020asymptotic}, the high order AP discontinuous Galerkin(AP-DG) method based on micro-macro decomposition \cite{xiong2022high} and the AP IMEX(AP-IMEX) method based on $P_N$ decomposition \cite{fu2022asymptotic}.

In recent years, the deep learning method has become a competitive method for solving partial differential equations(PDEs).
The idea is to approximate the unknown using the deep neural network by constructing suitable loss functions for optimization.
Various neural network methods have been proposed, such as the physics-informed neural network(PINN) method \cite{raissi2019physics} and the deep Galerkin method(DGM) \cite{sirignano2018dgm} with the loss established on the $L^2$-residual of the PDEs,
and the deep Ritz method(DRM) \cite{yu2018deep} with the loss based on the Ritz formulation, and so on. For other methods, we refer the reader to \cite{han2018solving,zhang2021mod,weinan2021algorithms,hwang2020trend}.

The PINNs and Model-Operator-Data Network(MOD-Net) have been applied to solve steady and unsteady linear radiative transfer models\cite{zhang2021mod,mishra2021physics,lu2022solving,chen2021solving,liu2021deep,jin2021asymptotic}.
But the deep learning method for the nonlinear time-dependent GRTEs has not been fully investigated.
We find that for GRTEs, the PINN loss function deteriorates in diffusion regimes with relatively small scale parameters $\epsilon$,
since PINNs tend to learn simplified models during the training process and lose the accuracy of asymptotic limit states.
To tackle this problem, designing an asymptotic-preserving loss function is necessary \cite{jin2021asymptotic}.

In this work, we reformulate the original GRTEs into a coupled system in AP form based on the idea of micro-macro decomposition and construct the corresponding AP loss function for the multiscale problems.
Beyond the AP property, we also take into account the conservation constraints of the microscopic quantities and the boundary constraints with additional diffusion limit information in particular cases in the AP loss function.
We prove that the approximation error of APNNs is bounded by the loss function.
The numerical examples in linear or stationary cases confirm the well approximation of the proposed APNN method.
We notice that in some highly-nonlinear and non-stationary cases, using only PDE model constraints, APNNs may be poorly trained.
Since a few labeled data can be acquired via the conventional schemes on coarse grids with low computational cost,
we propose model-data asymptotic-preserving networks(MD-APNNs) by adding a few labeled data to APNNs.
The numerical experiments show better performance of MD-APNNs on the non-stationary GRTEs in both kinetic regime($\epsilon=O(1)$) and diffusion regime($\epsilon \ll 1$) than those of APNNs and pure data-driven networks.

The rest of this paper is organized as follows. In Section \ref{sec2}, we introduce the nonlinear gray radiative transfer model and its corresponding micro-macro decomposition scheme. In Section \ref{sec3}, we present the details of the MD-APNN method based on micro-macro decomposition and define the AP loss function. Section \ref{sec4} is devoted to the convergence theory. We carry out several numerical experiments in Section \ref{sec5}. The results of simulating the time-dependent linear radiative transfer model and the stationary nonlinear GRTEs confirm the AP property of our method. Then numerical examples of nonlinear time-dependent GRTEs are presented to validate the efficiency of MD-APNNs. Finally, a conclusion is given in Section \ref{sec6}.

\section{The GRTEs and micro-macro decomposition scheme}\label{sec2}
\subsection{The gray radiative transfer equations}
The gray radiative transfer equations describe the radiative transfer and the energy exchange between radiation and material. Omitting the scattering terms and the internal source, we consider the scaled form of the GRTEs with initial and boundary conditions in a bounded domain $\tau \times D \times S^2$

\begin{subnumcases}{\label{eq:GRTEs}}
\frac{\epsilon^2}{c} \frac{\partial I}{\partial t}+\epsilon \Omega \cdot \nabla_r I =\sigma \left(\frac{1}{4 \pi}acT^4-I\right), & $(t,r,\Omega) \in\tau\times D \times S^2$, \ \ \ \ \label{eq:GRTEs-a} \\
\epsilon^2 C_v\frac{\partial T}{\partial t} \equiv \epsilon^2 \frac{\partial U}{\partial t} = \sigma \left(\int_{S^2} I \text{d} \Omega -acT^4\right), & $(t,r)\in\tau\times D$, \ \ \ \ \label{eq:GRTEs-b} \\
B_{ou}I = I_\Gamma(t,r,\Omega), & $(t,r,\Omega)\in \Gamma$, \ \ \ \ \label{eq:GRTEs-c} \\
I_{in}I = I_0(r,\Omega),  & $(t,r,\Omega)\in \{0\}\times D \times S^2$, \ \ \ \ \label{eq:GRTEs-d} \\
I_{in}T=T_0(r), & $(t,r)\in \{0\} \times D$, \ \ \ \ \label{eq:GRTEs-e}
\end{subnumcases}
where $r=(x,y,z) \in D$ is the spatial position variable with boundary $\partial D $, $t \in \tau$ the time variable, and $\Omega = (\xi,\eta,\mu)$ the angular direction on the unit sphere $S^2$ along which the photons propagate.
Eqs.\eqref{eq:GRTEs-a}\eqref{eq:GRTEs-b} are the governing equations,
where $I(t,r,\Omega)$ denotes the radiation intensity, $T(t,r)$ the material temperature, $U(t,r)$ the material energy density, $\sigma(r, T)$ the opacity, $\epsilon >0$ the Knudsen number defined as the ratio of the photon mean free path over the characteristic length of space, $a$ the radiation constant and $c$ the scaled speed of light.
The material temperature $T(t,r)$ and the material energy density $U(t,r)$ satisfy the relationship
\begin{equation*}
\frac{\partial U}{\partial T}=C_v>0, 
\end{equation*}
where $C_v(r,t)$ is the scaled heat capacity.
$I_0(r,\Omega)$ and $T_0(r)$ of Eqs.\eqref{eq:GRTEs-d}-\eqref{eq:GRTEs-e} are the initial values for $I$ and $T$, respectively.
The boundary condition Eq.\eqref{eq:GRTEs-c} will be specified later.

\subsubsection{The boundary conditions}
Separate the boundary $\Gamma= \tau\times\partial D \times S^2$ into the inflow and outflow parts as follows
\[
\Gamma_{-}= \left\{(t,r,\Omega)\in \Gamma : \Omega \cdot n_{r} <0\right\}, \qquad \Gamma_{+}= \left\{(t,r,\Omega)\in \Gamma : \Omega \cdot n_{r} >0\right\},
\]
where $n_r$ is the unit outer normal of $\partial D$ at position $r$.
In this paper, we consider the following three types of boundary conditions for GRTEs.
\begin{enumerate}[(1)]
	\item The Dirichlet boundary condition
	\begin{equation}\label{eq:BC-Dir}
	I(t,r,\Omega) = I_{B}(t,r,\Omega), \quad (t,r,\Omega)\in \Gamma_{-}, 
	\end{equation}
	where $I_{B}(t,r,\Omega)$ is a given function. 
	\item The reflecting boundary condition
	\begin{equation}\label{eq:BC-ref}
	I(t,r,\Omega) = I(t,r,\Omega^{'}), \quad (t,r,\Omega)\in \Gamma_{-}.
	\end{equation}
	where $\Omega^{'} = \Omega -2n_r(n_r \cdot \Omega)$ is the reflection of $\Omega$ in the tangent plane to $\partial D$.
	\item If $D$ is symmetric, we can enforce the periodic boundary condition \cite{lemou2008new}
	\begin{equation}\label{eq:BC-per}
	I(t,r,\Omega)=I(t,S(r),\Omega), \quad (t,r,\Omega) \in \Gamma,
	\end{equation}
	where $S:\partial D \rightarrow \partial D$ is a one-to-one mapping.
\end{enumerate}

\subsubsection{The 1D case and linear model}
In one-dimensional case ($D \subset \mathbb{R}$), Eqs.\eqref{eq:GRTEs-a}-\eqref{eq:GRTEs-b} reduce to
\begin{subnumcases}{\label{eq:GRTEs-1d}}
\frac{\epsilon^2}{c} \frac{\partial I}{\partial t}+\epsilon \mu \frac{\partial I}{\partial x} =\sigma \left(\frac{1}{2}acT^4 -I\right) & $(t,x,\mu) \in \tau \times D \times [-1,1]$, \label{eq:GRTEs-1d-a} \\
\epsilon^2 C_v\frac{\partial T}{\partial t} =\sigma \left(\int_{-1}^{1}I \text{d} \mu -acT^4 \right) & $(t,x) \in \tau \times D$ . \label{eq:GRTEs-1d-b}
\end{subnumcases}
When the material temperature is the same as the radiation temperature, the gray radiative transfer equation becomes the scaled linear transport model
\begin{equation}\label{eq:GRTEs-lin}
\frac{\epsilon^2}{c}\partial_t I + \epsilon \Omega \cdot \nabla_r I = \sigma\left(\frac{1}{4\pi}\int_{S^2} I \text{d} \Omega -I\right).
\end{equation}
And in 1D case, the above equation can be simplified as
\begin{equation}\label{eq:GRTEs-1d-lin}
\frac{\epsilon^2}{c} \partial_t I+\epsilon \mu \partial_x I =\sigma \left(\frac{1}{2}\int_{-1}^{1}I \text{d} \mu -I\right).
\end{equation}

\subsection{The diffusion limit}
Eqs.\eqref{eq:GRTEs-a}-\eqref{eq:GRTEs-b} are regarded as a relaxation model for the radiation intensity in local thermodynamic equilibrium, with the emission source of the background medium given by a Planck function at the local material temperature \cite{sun2015asymptotic}, namely $\frac{1}{4 \pi} \sigma acT^4$.
As the parameter $\epsilon \rightarrow 0$, \cite[Larsen et al.]{larsen1983asymptotic} shows that, away from boundaries and initial moments, the radiation intensity $I(t,r,\Omega)$ approaches to a Planck function at the local temperature, i.e.,
\begin{equation*}
I^{(0)}=\frac{1}{4 \pi} ac\left(T^{(0)}\right)^4,
\end{equation*}
and the local temperature $T^{(0)}$ satisfies the nonlinear diffusion equation
\begin{equation}\label{eq:Diff-lim}
\frac{\partial}{\partial t}\left(C_v T^{(0)}\right) + a \frac{\partial}{\partial t}\left(T^{(0)}\right)^4 = \nabla \cdot \frac{ac}{3\sigma}
\nabla \left(T^{(0)}\right)^4.
\end{equation}
Therefore, the asymptotic preserving scheme for Eq.\eqref{eq:GRTEs} should lead to a correct discretization of the diffusion limit equation \eqref{eq:Diff-lim} when the Knudsen number $\epsilon$ is extremely small.
Moreover, the scheme should be uniformly stable with respect to $\epsilon$.

\subsection {The micro-macro decomposition method}
We denote by
\[
\left \langle I \right \rangle = \frac{1}{4 \pi}\int_{S^2} I(t,r,\Omega) \text{d} \Omega
\]
the integral average of $I$ over the angular variable $\Omega$,
and make decomposition of the radiation intensity as in \cite{xiong2022high} by
\begin{equation*}
I(t,r,\Omega) = \rho(t,r)+\frac{\epsilon}{\sqrt{\sigma_0}} g(t,r,\Omega),
\end{equation*}
where $\rho=\left \langle I \right \rangle, \left \langle g \right \rangle=0$, $\sigma _0 >0$ is a constant defined as a referred opacity, and $\sqrt{\sigma_0}$ is added when $\sigma$ is relatively large.
Integrating Eq.\eqref{eq:GRTEs-a} over the angular direction $\Omega$ of photons and subtracting it from Eq.\eqref{eq:GRTEs-a}, together with Eq.\eqref{eq:GRTEs-b},
we get the micro-macro decomposition coupled system
\begin{subnumcases}{\label{eq:GRTEs-mmd}}
\frac{1}{c}\partial_t \rho + \frac{1}{\sqrt{\sigma_0}} \left \langle \Omega \cdot \nabla_r g \right \rangle = -\frac{1}{4\pi}C_v \partial_tT, \label{eq:GRTEs-mmd-a} \\
\frac{\epsilon^2}{c}\partial_t g +\epsilon \Omega\cdot \nabla_r g-\epsilon \left \langle \Omega \cdot \nabla_r g \right \rangle + \sqrt{\sigma_0} \Omega\cdot \nabla_r \rho + \sigma g=0, \label{eq:GRTEs-mmd-b}\\
\epsilon^2 C_v \partial_tT = \sigma\left(4\pi\rho-acT^4\right). \label{eq:GRTEs-mmd-c}
\end{subnumcases}
As for the boundary conditions, first we have the Dirichlet boundary condition for Eq.\eqref{eq:GRTEs-mmd}
\begin{equation}\label{eq:BC-Dir-mmd}
\rho(t,r) + \frac{\epsilon}{\sqrt{\sigma_0}} g(t,r,\Omega)=I_{B}(t,r,\Omega), \quad (t,r,\Omega)\in \Gamma_{-}.
\end{equation}
Then the reflecting boundary conditions together with the implicit constraint on macroscopic quantity under diffusion regime are stated as follows
\begin{subequations}\label{eq:BC-ref-mmd}
	\begin{align}
	& \nabla_r \rho(t,r)=0, \qquad  \qquad \ \ \ (t,r)\in \tau\times \partial D, \label{eq:BC-ref-mmd-a} \\
	& g(t,r,\Omega) = g(t,r,\Omega^{'}), \quad (t,r,\Omega)\in \Gamma_{-}. \label{eq:BC-ref-mmd-b}
	\end{align}
\end{subequations}
And the periodic boundary conditions are given by
\begin{subequations}\label{eq:BC-per-mmd}
	\begin{align}
	& \rho(t,r)=\rho(t,S(r)),\qquad \quad (t,r)\in \tau\times \partial D, \label{eq:BC-per-mmd-a} \\
	& g(t,r,\Omega) = g(t,S(r),\Omega), \quad (t,r,\Omega) \in \Gamma. \label{eq:BC-per-mmd-b}
	\end{align}
\end{subequations}
For the initial conditions, we impose
\begin{subequations}\label{eq:IC-mmd}
	\begin{align}
	& \rho(0,r) + \frac{\epsilon}{\sqrt{\sigma_0}} g(0,r,\Omega)=I_0(r,\Omega), \quad (r,\Omega)\in D \times S^2,  \label{eq:IC-mmd-a} \\
	& T(0,r)=T_0(r), \qquad \qquad \qquad \qquad \quad  r \in D. \label{eq:IC-mmd-b}
	\end{align}
\end{subequations}
Due to the fact that $\frac{\epsilon^2}{c}\partial_t g +\epsilon \Omega\cdot \nabla_r g-\epsilon \left \langle \Omega \cdot \nabla_r g \right \rangle $ and $\epsilon^2 \partial_t T$ vanish as $\epsilon \rightarrow 0$, we obtain from Eqs.\eqref{eq:GRTEs-mmd-b} and \eqref{eq:GRTEs-mmd-c} that the limit case($\epsilon \rightarrow 0$) satisfies
\[
g=-\frac{\sqrt{\sigma_0}}{\sigma}\Omega \cdot \nabla_r \rho, \qquad \rho = \frac{1}{4\pi}acT^4,
\]
which, together with Eq.\eqref{eq:GRTEs-mmd-a}, result the nonlinear diffusion limit equation \eqref{eq:Diff-lim}. This
indicates the advantage of the micro-macro decomposition method.

\section{Asymptotic-preserving neural networks (APNNs)}\label{sec3}
Given an input $X=(t,r,\Omega) \in K \equiv \tau \times D\times S^2$, the output of the  $L$-layer deep feedforward neural network $F_{\theta}$ is a combination of the affine transformations $\{C_i\}_{i=1}^L$ and nonlinear activation functions $\{\sigma_i\}_{i=1}^L$, stated as follows
\[
F_{\theta}(X)= \sigma_{L} \circ C_L \circ \sigma_{L-1} \circ C_{L-1} ......\circ \sigma_2 \circ C_2 \circ \sigma_1 \circ C_1(X).
\]
Popular choices of $\sigma_i$ include the sigmoid function, the hyperbolic tangent function, the ReLU function, and so on.
To be more specific, the relation between the $l$-layer and $l+1$-layer($l=0,1,2,...L-1$) in $F_{\theta}$ is given by
\[
u_j^{l+1}=\sigma_{l+1}C_{l+1}(u_l) =\sigma_{l+1}\left(\sum_{i}^{m_l}W_{ij}^{l+1} u_i^l+b_j^{l+1}\right),
\]
where $u_i^l$ denotes the value of the $i$-th neuron in the $l$-th layer, $W_{ij}^{l+1}\in R^{m_{l+1}\times m_l}$ the weight between the $i$-th neuron in the $l$-th layer and the $j$-th neuron in the $l+1$-th layer, $b_j^{l+1}\in R^{m_{l+1}}$ the bias of the $j$-th neuron in the $l+1$-th layer, $m_l$ the number of neurons in the $l$-th layer and $\sigma_{l+1}$  the activation function in the $l+1$-th layer.
For simplicity, we denote the network structure by $[m_0, m_1,...,m_L]$, with the input dimension $m_0=d_{in}$ and the output dimension $m_L=d_o$.

Solving the PDE model with deep neural network methods usually consists of three parts: (1) neural network parameterization of PDE solutions; (2) construction of population or empirical loss functions; (3) selection of appropriate optimization algorithms. One of the most typical neural network methods is the PINN method, the loss function of which is constructed on the $L^2$-residuals of the PDE system. For the APNN method \cite{jin2021asymptotic}, the key idea is to design the loss function with the AP property, which means as the physical scaling parameter tends to zero, the loss function of the multiscale model is consistent with the loss function of the diffusion limit equation. In particular, APNNs reformulate the original PDE into a new system in AP form and build the $L^2$-residual loss of the new system. The boundary, initial, and conservation laws are treated as regularization terms in the new loss function.

\subsection{PINNs fail to resolve GRTEs with small scales}
In the process of applying PINNs, we approximate radiation intensity $I(t,r,\Omega)$ and material temperature $T(t,r)$ by the neural networks $I_{\theta_{1}}^{nn}(t,r,\Omega)$ and $T_{\theta_{2}}^{nn}(t,r)$, respectively.
Consider the fact that $I$ and $T$ are nonnegative variables, we add proper activation function to the output layer of $I_{\theta_{1}}^{nn}$ and $T_{\theta_{2}}^{nn}$ to ensure the nonnegativity.
The PINN loss function of the original GRTEs (Eq.\eqref{eq:GRTEs}) is constructed as the least square of the residuals of the governing equations, combining with the initial and boundary conditions as the regularization terms
\begin{equation}\label{eq:PINNs-loss}
\begin{aligned}
& L_{\text{PINNs}}^{\epsilon} = \left\| \frac{\epsilon^2}{c}\partial_t I_{\theta_1}^{nn} + \epsilon \Omega\cdot \nabla_r I_{\theta_1}^{nn}-\sigma\left(\frac{1}{4 \pi}ac (T^{nn}_{\theta_{2}})^4-I_{\theta_{1}}^{nn} \right) \right\|^2_{L^2(K)} \\
& \ + \left\| \epsilon^2 C_v \partial_tT_{\theta_{2}}^{nn} - \sigma\left(\int_{S^2} I_{\theta_{1}}^{nn}\text{d}\Omega-ac(T^{nn}_{\theta_{2}})^4\right) \right\|^2_{L^2(\tau \times D)}  +\lambda_{1}\Vert B_{ou}I_{\theta_{1}}^{nn}-I_{\Gamma} \Vert^2_{L^2(\Gamma)} \\
& \ +\lambda_{2} \left(\Vert I_{in}I_{\theta_{1}}^{nn}-I_0 \Vert^2_{L^2(D \times S^2)} +\Vert I_{in}T_{\theta_{2}}^{nn}-T_0 \Vert^2_{L^2(D)}\right) \equiv \sum_{\text{i}=1}^4 L_{\text{PINNs,i}}^\epsilon,
\end{aligned}
\end{equation}
where $\lambda_1$ and $\lambda_2$ are the weight parameters.
Let us investigate the AP property of the PINN method.
Taking the limit $\epsilon \rightarrow 0$, we see that the loss of the governing equations satisfies
\begin{equation}
\begin{aligned}
L_{\text{PINNs},g}^{\epsilon} & \equiv  L_{\text{PINNs},1} + L_{\text{PINNs},2} \rightarrow \left\| \sigma\left(\frac{1}{4 \pi}ac({T^{nn}_{\theta_{2}}})^4-I_{\theta_{1}}^{nn} \right) \right\|^2_{L^2(K)} \\
& \quad + \left\| \sigma\left(\int_{S^2} I_{\theta_{1}}^{nn}\text{d}\Omega-ac({T^{nn}_{\theta_{2}}})^4\right) \right\|^2_{L^2(\tau \times D)} \equiv L_{\text{PINNs},g}.
\end{aligned}
\end{equation}
The right-hand side can be viewed as the PINNs' loss of the following system
%
\[
\sigma\left(\frac{1}{4 \pi}acT^4-I\right)=0, \quad 
\sigma\left(\int_{S^2} I\text{d}\Omega-acT^4\right)=0,
\]
%
which further implies that $I=\rho=\frac{1}{4\pi}acT^4$.
However, it is not the desired diffusion limit equation \eqref{eq:Diff-lim}.
As a result, PINNs will fail when $\epsilon$ is sufficiently small.

\subsection{The micro-macro decomposition based APNN method}
The APNN method is to apply PINNs to solve the micro-macro decomposition system Eq.\eqref{eq:GRTEs-mmd}.
We use two $L$-layer deep neural networks $g^{nn}_{\theta_1}(t,r,\Omega) \approx g(t,r,\Omega)$, and $(\rho,T)^{nn}_{\theta_{2}}(t,r) = (\rho ^{nn}_{\theta_{21}}(t,r)$, $T^{nn}_{\theta_{22}}(t,r)) \approx (\rho(t,r), T(t,r))$.
Since $\rho$ and $T$ are nonnegative, we choose the appropriate activation function $\sigma^o(X)$ at the output layer of $(\rho, T)^{nn}_{\theta_{2}}$ to guarantee the nonnegativity.
For instance, we take $\sigma^o(X)=e^{-X}$, $\sigma^o(X)=X$, $\sigma^o(X)=\ln(1+e^{X})$.

Now we design the APNN loss $L_{\text{APNNs}}^{\epsilon}$, containing the residuals of the micro-macro governing equation, boundary condition, initial value, and the conservation law for the microscopic quantity $g$.
\begin{equation}\label{eq:APNNs-loss}
L_{\text{APNNs}}^{\epsilon} = L_{\text{APNNs},g}^{\epsilon} + L_{\text{APNNs},i}^{\epsilon} + L_{\text{APNNs},b}^{\epsilon} + L_{\text{APNNs},c}^{\epsilon},
\end{equation}
where
\[
\begin{aligned}
L_{\text{APNNs},g}^{\epsilon} = & \left\| \frac{\epsilon^2}{c}\partial_t g_{\theta_1}^{nn} + \epsilon \Omega\cdot \nabla_r g_{\theta_1}^{nn}-\epsilon \left \langle \Omega \cdot \nabla_r g_{\theta_1}^{nn} \right \rangle + \sqrt{\sigma_0} \Omega\cdot \nabla_r \rho_{\theta_{21}}^{nn} + \sigma g_{\theta_1}^{nn} \right\|^2_{L^2(K)}  \\
& + \left\| \frac{1}{c}\partial_t \rho_{\theta_{21}}^{nn} +\frac{1}{\sqrt{\sigma_0}} \left \langle \Omega \cdot \nabla_r g_{\theta_1}^{nn} \right \rangle + \frac{1}{4\pi}C_v \partial_tT_{\theta_{22}}^{nn} \right\|^2_{L^2(\tau \times D)}  \\
& + \left\| \epsilon^2 C_v \partial_tT_{\theta_{22}}^{nn} - \sigma\left(4\pi\rho_{\theta_{21}}^{nn}-ac({T^{nn}_{\theta_{22}}})^4\right) \right\|^2_{L^2(\tau \times D)},
\end{aligned}
\]
\[
L_{\text{APNNs},b}^{\epsilon} = \lambda_{1} \left\| B_{ou}\left(\rho^{nn}_{\theta_{21}} + \frac{\epsilon}{\sqrt{\sigma_0}} g^{nn}_{\theta_1}\right) - I_{\Gamma} \right\|^2_{L^2(\Gamma)},
\]
\[
L_{\text{APNNs},i}^{\epsilon} = \lambda_{2} \left( \left\| I_{in}\left(\rho^{nn}_{\theta_{21}} + \frac{\epsilon}{\sqrt{\sigma_0}} g^{nn}_{\theta_1}\right)-I_0 \right\|^2_{L^2(D \times S^2)} + \left\| I_{in}T^{nn}_{\theta_{22}}-T_0 \right\|^2_{L^2(D)} \right),
\]
\[
L_{\text{APNNs},c}^{\epsilon} = \lambda_{3} \Vert \left \langle g_{\theta_1}^{nn} \right \rangle \Vert^2_{L^2(\tau \times D)}.
\]
To demonstrate the AP property, we observe that, as $\epsilon \rightarrow 0$,
\[
\begin{aligned}
L_{\text{APNNs},g}^\epsilon \rightarrow & \left\| \sqrt{\sigma_0} \Omega\cdot \nabla_r \rho_{\theta_{21}}^{nn} + \sigma g_{\theta_1}^{nn} \right\|^2_{L^2(K)}  \\
&\quad + \left\| \frac{1}{c}\partial_t \rho_{\theta_{21}}^{nn} +\frac{1}{\sqrt{\sigma_0}} \left \langle \Omega \cdot \nabla_r g_{\theta_1}^{nn} \right \rangle + \frac{1}{4\pi}C_v \partial_tT_{\theta_{22}}^{nn} \right\|^2_{L^2(\tau \times D)} \\
&\quad + \left\| \sigma\left(4\pi\rho_{\theta_{21}}^{nn}-ac({T^{nn}_{\theta_{22}}})^4\right) \right\|^2_{L^2(\tau \times D)} \equiv L_{\text{APNNs},g},
\end{aligned}
\]
which can be regarded as the PINNs' loss of the following system
\begin{subnumcases}{\label{eq:APNN-grov-lim}}
\sqrt{\sigma_0} \Omega\cdot \nabla_r \rho + \sigma g=0, \label{eq:APNN-grov-lim-a}\\
\frac{1}{c}\partial_t \rho +\frac{1}{\sqrt{\sigma_0}} \left \langle \Omega \cdot \nabla_r g \right \rangle + \frac{1}{4\pi}C_v \partial_tT=0, \label{eq:APNN-grov-lim-b} \\
\sigma\left(4\pi\rho-acT^4\right)=0. \label{eq:APNN-grov-lim-c}
\end{subnumcases}

Substituting Eqs.\eqref{eq:APNN-grov-lim-a} and \eqref{eq:APNN-grov-lim-c} into Eq.\eqref{eq:APNN-grov-lim-b}, we get the diffusion limit equation \eqref{eq:Diff-lim}.
Therefore, the loss function $L_{\text{APNNs}}^\epsilon$ possesses the AP property and the APNN method is well applicable to the case where $\epsilon$ is extremely small.

\subsection{Discretization of the loss function}
\subsubsection{Angular integration approximation}
In practice, the angular integration terms $\langle g \rangle$, $\langle \Omega\cdot \nabla_r g \rangle$ of Eq.\eqref{eq:GRTEs-mmd} need to be approximated by proper quadrature rules
\begin{subequations}{\label{Angular-integration}}
	\begin{align}
	&\left \langle g \right \rangle=\frac{1}{4\pi}\int_{S^2} g(t,r,\Omega) \text{d}\Omega \approx \frac{1}{4\pi} \sum_{m=1}^{N_s} g(\Omega_m) \omega_m, \label{Angular-integration-a} \\
	&\left \langle \Omega\cdot \nabla_r g \right \rangle=\frac{1}{4\pi}\int_{S^2} \Omega\cdot \nabla_r g(t,r,\Omega) \text{d}\Omega \approx \frac{1}{4\pi} \sum_{m=1}^{N_s} \Omega_m \cdot \nabla_r g(\Omega_m) \omega_m, \label{Angular-integration-b}
	\end{align}
\end{subequations}
where the quadrature points $\{\Omega_m\}_{m=1}^{N_s} \subset S^2$ are also used as the training points, and $\{\omega_m\}_{m=1}^{N_s}$ denote the associated weights.
In 1D case, $\Omega = \mu \in[-1,1]$, and we can apply the Gauss-Legendre quadrature rule for approximation with the Legendre points $\mu_m \in[-1,1]$ and weights $\omega_m \in(0,1)$ ($m=1,2...,N_s$).
\subsubsection{Empirical loss functions, sample points, and labeled data}
In order to increase the accuracy in the training process associated with the highly nonlinear, nonstationary GRTEs,
we can add some labeled data,
obtained by experimental observation or by conventional schemes on a coarse grid, to the loss function as a regularization term.
As a result, we get the model-data asymptotic-preserving neural networks(MD-APNNs).

In practice, we approximate the loss function by suitable quadrature rules.
Here we sample the training points $\wp_{f}=\left\{X^f_j=(t^f_j, r^f_j, \Omega^f_j), j=1,2,...N_f\right\}$ as a low-discrepancy Sobol sequence in the computational domain, and choose the Quasi-Monte Carlo method for numerical integration \cite{sobol1967distribution}.
The empirical loss function of the governing equation is stated as follows
\[
\begin{aligned}
& \quad L_{\text{MD-APNNs},g}^{\epsilon,nn} = \frac{1}{N_{\text{int}}} \sum_{j=1}^{N_{\text{int}}} \Bigg( \Big| \frac{\epsilon^2}{c}\partial_t g_{\theta_1}^{nn}(t_j^{\text{int}},r_j^{\text{int}},\Omega_j^{\text{int}}) + \epsilon \Omega\cdot \nabla_r g_{\theta_1}^{nn}(t_j^{\text{int}},r_j^{\text{int}},\Omega_j^{\text{int}})  \\
& -\epsilon \frac{1}{4\pi}\sum_{m=1}^{N_s}\Omega_m \cdot \nabla_r g_{\theta_1}^{nn}(t_j^{\text{int}},r_j^{\text{int}},\Omega_m^{\text{int}})\omega_m + \sqrt{\sigma_0} \Omega\cdot \nabla_r \rho_{\theta_{21}}^{nn}(t_j^{\text{int}},r_j^{\text{int}}) + \sigma g_{\theta_1}^{nn}(t_j^{\text{int}},r_j^{\text{int}},\Omega_j^{\text{int}}) \Big|^2  \\
& + \Big| \frac{1}{c}\partial_t \rho_{\theta_{21}}^{nn}(t_j^{\text{int}},r_j^{\text{int}}) +\frac{1}{\sqrt{\sigma_0}} \frac{1}{4\pi}\sum_{m=1}^{N_s}\Omega_m \cdot \nabla_r g_{\theta_1}^{nn}(t_j^{\text{int}},r_j^{\text{int}},\Omega_m^{\text{int}})\omega_m + \frac{1}{4\pi}C_v \partial_tT_{\theta_{22}}^{nn}(t_j^{\text{int}},r_j^{\text{int}}) \Big|^2 \\
& + \Big| \epsilon^2 C_v \partial_tT_{\theta_{22}}^{nn}(t_j^{\text{int}},r_j^{\text{int}}) - \sigma\left(4\pi\rho_{\theta_{21}}^{nn}(t_j^{\text{int}},r_j^{\text{int}})-ac({T^{nn}_{\theta_{22}}})^4(t_j^{\text{int}},r_j^{\text{int}})\right) \Big|^2 \Bigg),
\end{aligned}
\]
where $\{(t_j^{\text{int}},r_j^{\text{int}},\Omega_j^{\text{int}})\}_{j=1}^{N_{\text{int}}}$ are interior Sobol sequence points, $\{\Omega_m^{\text{int}} \}_{m=1}^{N_s}$ are the Gauss-quadrature points and $\{\omega_m\}_{m=1}^{N_s}$ are the corresponding quadrature weights.

We also have the empirical loss functions for the boundary condition, the initial value, and conservation law of $\langle g \rangle$
\[
L_{\text{MD-APNNs},b}^{\epsilon,nn} = \lambda_{1}\frac{1}{N_{\text{sb}}}\sum_{j=1}^{N_{\text{sb}}} \left| B_{ou}\left(\rho_{\theta_{21}}^{nn}(t_j^{\text{sb}},r_j^{\text{sb}})+\frac{\epsilon}{\sqrt{\sigma_{0}}}g_{\theta_{1}}^{nn}(t_j^{\text{sb}},r_j^{\text{sb}},\Omega_j^{\text{sb}})\right)-I_{\Gamma}(t_j^{\text{sb}},r_j^{\text{sb}},\Omega_j^{\text{sb}}) \right|^2,
\]
\[
\begin{aligned}
L_{\text{MD-APNNs},i}^{\epsilon,nn} & = \lambda_{2}\frac{1}{N_{\text{tb}}}\sum_{j=1}^{N_{\text{tb}}} \Bigg( \Big| I_{in}\left(\rho_{\theta_{21}}^{nn}(0,r_j^{\text{tb}})+\frac{\epsilon}{\sqrt{\sigma_{0}}}g_{\theta_{1}}^{nn}(0,r_j^{\text{tb}},\Omega_j^{\text{tb}})\right)-I_0(r_j^{\text{tb}},\Omega_j^{\text{tb}}) \Big|^2 \\
& \quad + \Big| I_{in}T_{\theta_{22}}^{nn}(0,r_j^{\text{tb}})-T_0(r_j^{\text{tb}}) \Big|^2 \Bigg),
\end{aligned}
\]
\[
L_{\text{MD-APNNs},c}^{\epsilon,nn} = \lambda_{3} \frac{1}{N_{\text{int}}}\sum_{j=1}^{N_{\text{int}}} \left( \left|\frac{1}{4\pi}\sum_{m=1}^{N_s} g_{\theta_1}^{nn}(t_j^{\text{int}},r_j^{\text{int}},\Omega_m^{\text{int}})\omega_m \right|^2 \right),
\]
where $\{ (t_j^{\text{sb}},r_j^{\text{sb}},\Omega_j^{\text{sb}}) \}_{j=1}^{N_{\text{tb}}}$ are the boundary Sobol sequence points, $\{(0,r_j^{\text{tb}},\Omega_j^{\text{tb}})\}_{j=1}^{N_{\text{tb}}}$ the initial Sobol sequence points.

Suppose we have already obtained some label data $\{T^*(t_i,r_i)\}_{i=1}^{N_0}$ of the material temperature $T$ via conventional schemes on a coarse grid.
We introduce the loss function of the data regularization
\[
L_{\text{MD-APNNs},l}^{\epsilon,nn} = \lambda_0\frac{1}{N_{0}}\sum_{i=1}^{N_0} \left( \left| T_{\theta_{22}}^{nn}(t_i,r_i)-T^*(t_i,r_i) \right|^2 \right).
\]
Here, $\lambda_0,\lambda_{1},\lambda_{2},\lambda_{3}$ are adjustable hyperparameters.
Summing up the above empirical loss functions, we get the total empirical loss function for MD-APNNs
\begin{equation}
L_{\text{MD-APNNs}}^{\epsilon,nn} = \sum_{\mathrm{i}=\{g,b,i,c,l\}} L_{\text{MD-APNNs,i}}^{\epsilon,nn} .
\end{equation}
Now, we are in the position to find the solution of the minimization problem
\begin{equation}
\theta_1^{*},\theta_2^{*}= \mathop{\arg\min}\limits_{\theta_1,\theta_2}\left( L_{\text{MD-APNNs}}^{\epsilon,\text{nn}}\right).
\end{equation}
In simulation, we adopt appropriate optimization algorithms to minimize the non-convex loss function $L_{\text{MD-APNNs}}^{\epsilon,\text{nn}}$.

\section{The convergence anlaysis}\label{sec4}
\noindent \textbf{Assumption 4.1} Assume there exist positive constants $\sigma_{\min}$, $\sigma_{\max}$ and $T_{\min}$, $T_{\max}$ such that
\begin{subequations}{\label{sigmaT-range}}
	\begin{align}
	0 \textless \sigma_{\min} \le \sigma(r,T) \le \sigma_{\max} \textless \infty \quad \forall(r,T)\in (D\times R^{+}), \label{sigmaT-range-a} \\
	0 \textless  T_{\min} \le T(t,r), T^*(t,r) \le T_{\max} \textless \infty \quad \forall(t,r)\in (\tau\times D). \label{sigmaT-range-b}
	\end{align}
\end{subequations}

\begin{remark}
	{\rm We aim to demonstrate the convergence of the proposed method. Let $(g^{*}, \rho^{*}, T^{*})$ be the solution of Eq.\eqref{eq:GRTEs-mmd}. To simplify the notation, here and hereafter, we use $A \lesssim B$ to express $A \le CB$ with a constant $C$ independent of the size of the neural network. Suppose $(g, \rho, T)$ is obtained by minimizing the loss $L_{\text{APNNs}}^{\epsilon, nn}$ (no labeled data are used) such that $L_{\text{APNNs}}^{\epsilon, nn} \le \zeta \le 1$. Our goal is to claim that $\Vert (g, \rho, T)- (g^{*}, \rho^{*}, T^{*}) \Vert_{L^{\infty}(\tau; L^2(D\times S^2))}^2 \lesssim \zeta$ holds. Note that if the label data are involved, the same assertion still holds because $L_{\text{MD-APNNs}}^{\epsilon, nn} \ge L_{\text{APNNs}}^{\epsilon, nn}$.}
\end{remark}

\begin{remark}
	{\rm If $(g^{*}, \rho^{*}, T^{*})$ is sufficiently smooth, according to the approximation theory of NN \cite{hornik1989multilayer}, for arbitray $\zeta \textgreater 0$, exist a parameter set $(\theta_{1},\theta_{21},\theta_{22})$ such that $\Vert (g_{\theta_{1}}^{nn}, \rho_{\theta_{21}}^{nn}, T_{\theta_{22}}^{nn})- (g^{*}, \rho^{*}, T^{*}) \Vert_{C^1(\tau \times D \times S^2)} \lesssim \zeta$, which implies $L_{\text{APNNs}}^{\epsilon}(g_{\theta_{1}}^{nn}, \rho_{\theta_{21}}^{nn}, T_{\theta_{22}}^{nn}) \lesssim \zeta$. Moreover, if the errors of numerical integration (the angular integration formula \eqref{Angular-integration} and the quasi-Monte-Carlo method to approximate $L_{\text{APNNs}}^{\epsilon}$ by $L_{\text{APNNs}}^{\epsilon,nn}$ ) are small enough, then we achieve $L_{\text{APNNs}}^{\epsilon,nn}(g_{\theta_{1}}^{nn}, \rho_{\theta_{21}}^{nn}, T_{\theta_{22}}^{nn}) \lesssim \zeta$.
		In summary, for smooth $(g^{*}, \rho^{*}, T^{*})$ and arbitary $\zeta \textgreater 0$, there exists a NN $(g_{\theta_{1}}^{nn}, \rho_{\theta_{21}}^{nn}, T_{\theta_{22}}^{nn})$ such that $L_{\text{APNNs}}^{\epsilon,nn} \le \zeta$. Conversely, for $(g_{\theta_{1}}^{nn}, \rho_{\theta_{21}}^{nn}, T_{\theta_{22}}^{nn})$ satisfying $L_{\text{APNNs}}^{\epsilon,nn} \le \zeta$. Suppose $(g_{\theta_{1}}^{nn}, \rho_{\theta_{21}}^{nn}, T_{\theta_{22}}^{nn})$ is smooth, and the errors of quadrature rules and sufficiently small, then $L_{\text{APNNs}}^{\epsilon,nn} \le \zeta$ also implies $L_{\text{APNNs}}^{\epsilon} \le \zeta$. Now the key question is: Does $L_{\text{APNNs}}^{\epsilon} \le \zeta$ ensures $\Vert (g, \rho, T)- (g^{*}, \rho^{*}, T^{*}) \Vert_{L^{\infty}(\tau; L^2(D\times S^2))}^2 \lesssim \zeta$? The answer is stated in the following theorem.}
\end{remark}

\begin{theorem}\label{th:error}
	Let the solutions to Eq.\eqref{eq:GRTEs} satisfy $ I\in C^1(\tau \times D; L^2(D\times S^2)) $, $T \in C^1(\tau; L^2(D))$. 
	Then we have
	\begin{equation}{\label{theorem}}
	\Vert I-I^*\Vert_{L^{\infty}\left(\tau;{L^2(D\times S^2)}\right)}^2 + \Vert T-T^* \Vert_{L^{\infty}\left(\tau;{L^2(D)}\right)}^2 \le C(\epsilon)L_{\text{APNNs}}^{\epsilon},
	\end{equation}
	where $C(\epsilon)$ is a constant depending on $(a,c,\sigma,C_v,\alpha_1,\alpha_2,T_{\min},T_{\max},\sigma_{\min},\sigma_{\max})$.
\end{theorem}

The detailed proof of Theorem~\ref{th:error} is provided in Appendix A.

\section{Numerical Experiment}
\label{sec5}
We present several examples to make a comparison on the performance of PINNs, APNNs, and MD-APNNs.
Firstly, we apply both PINNs and APNNs to the linear radiative transfer equations and the nonlinear stationary GRTEs,
and the results confirm that the APNN method has advantage in solving diffusive scaling problems.
Secondly, we carry out numerical simulation of non-stationary nonlinear GRTEs in transport and diffusion regimes by Data-driven networks, APNNs, and MD-APNNs with a few labeled data. It shows that the MD-APNN method has a better performance.

For all examples, we adopt the Adam algorithm with an initial learning rate $0.001$ in the optimization process of the loss function, take $\tanh x$ as the activation function of all hidden layers, choose $\sigma_0=1$ and set the number of Gauss-quadrature points $N_s=10$.
We calculate $L^2$-norm relative errors of $I$ and $T$ in the form
\[
L^2_{\text{error}}(I) = \frac{\Vert I^{nn}-I^{\text{Ref}} \Vert_{L^{2}(K)}}{\Vert I^{\text{Ref}} \Vert_{L^{2}(K)}}, \quad L^2_{\text{error}}(T) = \frac{\Vert T^{nn}-T^{\text{Ref}}\Vert_{L^{2}(\tau \times D)}}{\Vert T^{\text{Ref}} \Vert_{L^{2}(\tau \times D)}}.
\]
\subsection{Time-dependent linear transport equations}
\subsubsection{Diffusion regime}\label{ex:1}
\begin{figure}[!htbp]
	\setlength{\abovecaptionskip}{0.cm}
	\setlength{\belowcaptionskip}{-0.cm}
	\centering
	{
		\begin{minipage}{2.2in}
			\centering
			\includegraphics[width=2.2in]{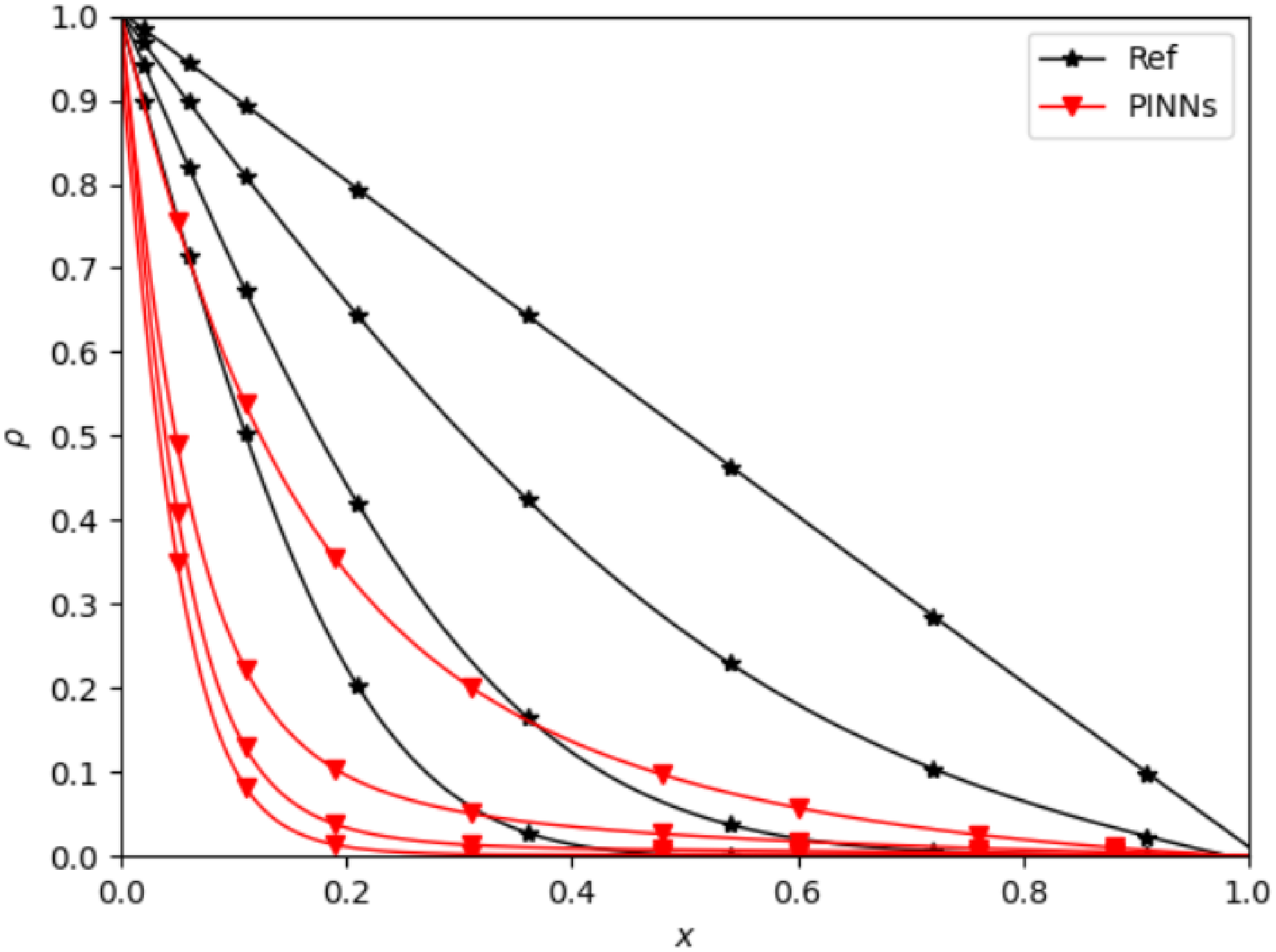}
		\end{minipage}
	}
	{
		\begin{minipage}{2.2in}
			\centering
			\includegraphics[width=2.2in]{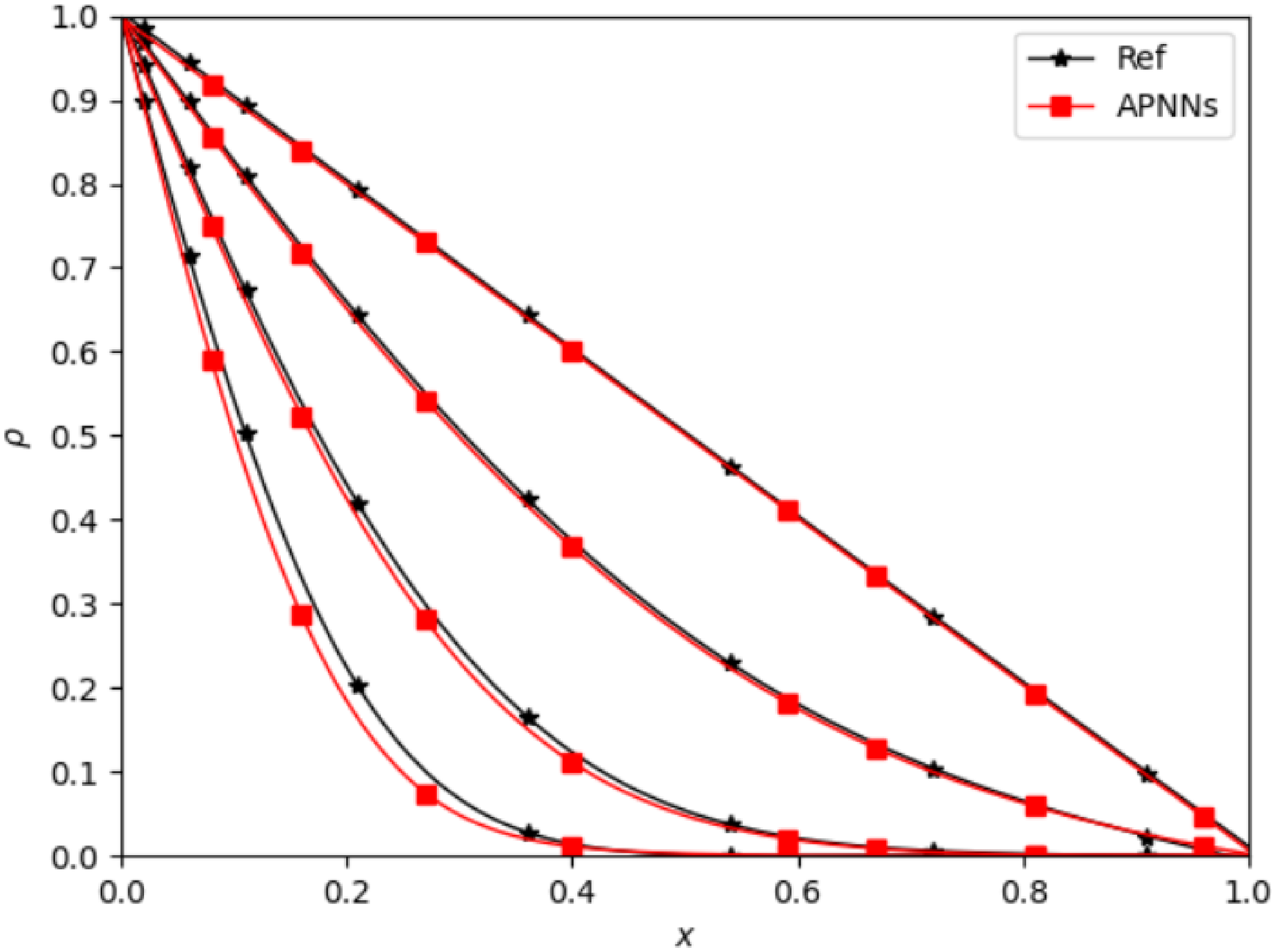}
		\end{minipage}
	}
	\caption{Diffusion regime with $\epsilon=10^{-8}$. The density $\rho$ at times $t=0.04, 0.1, 0.3, 2.0$. (Left) Ref v.s. PINNs. (Right) Ref v.s. APNNs.}
	\label{figure5.1}
\end{figure}

We apply PINNs and APNNs to solve the 1D linear radiative transfer equation \eqref{eq:GRTEs-1d-lin} under the diffusion regime with $\epsilon=10^{-8}$, i.e.,
\[\left\{
\begin{aligned}
&\epsilon \partial_t I + \mu \partial_x I =\frac{1}{\epsilon}\left(\left \langle I \right \rangle - I\right) && x\in[0,1], \ t\in[0,2], \ \mu \in[-1,1], \\
&I(t,0,\mu >0)=1 && t\in[0,2], \ \mu \in(0,1], \\
&I(t,1,\mu <0)=0 && t\in[0,2], \ \mu \in[-1,0), \\
&I(0,x,\mu )=0  && x\in[0,1], \ \mu \in[-1,1],
\end{aligned}
\right.
\]
where $c=1$ and $\sigma=1$.
We set $I^{nn}_{\theta}=[3,40,40,40,40,1]$ for PINNs, and $g^{nn}_{\theta_1}=[3,40,40,40,40,1]$, $\rho^{nn}_{\theta_2}=[2,40,40,40,40,1]$ for APNNs.
Set $N_{int}=16384, N_{sb}=N_{tb}=12288$ (the number of the interior, boundary and initial training points), and $\lambda_{1}=\lambda_{2}=\lambda_{3}=1$.
The activation function of the output layer are chosen to be $\sigma^o_I(X)=X,\sigma^o_g(X)=X, \sigma^o_{\rho}(X)=\ln(1+e^X)$ to ensure the nonnegativity.
No labeled data is utilized in this example.
The reference solution is obtained by the finite difference discretization of the diffusion limit equation. The reference solution and the prediction obtained by PINNs and APNNs are plotted in Figure~\ref{figure5.1}. The $L^2$ relative errors of PINNs and APNNs are shown in Table~\ref{Table5.1}, which show that the results of the micro-macro decomposition based APNNs are in a better agreement with the reference solutions than those of PINNs.
\begin{table}[!htbp]
	\centering
	\caption{Diffusion regime with $\epsilon=10^{-8}$. The errors of PINNs and APNNs.}
	\label{Table5.1}
	\begin{tabular}{lllll}
		\hline\noalign{\smallskip}
		$L^2_{\text{error}}(\rho)$ &$t=0.04$ & $t=0.1$ & $t=0.3$ & $t=2.0$ \\
		\noalign{\smallskip}\hline\noalign{\smallskip}
		PINNs &5.95e-01    &6.76e-01      &7.17e-01  &5.97e-01  \\
		APNNs &6.85e-02    &2.80e-02    &1.29e-02 & 6.71e-03  \\
		\noalign{\smallskip}\hline
	\end{tabular}
\end{table}
\subsubsection{Intermediate regime with a variable scattering frequency}
\begin{figure}[!htbp]
	\centering
	{
		\begin{minipage}{2.2in}
			\centering
			\includegraphics[width=2.2in]{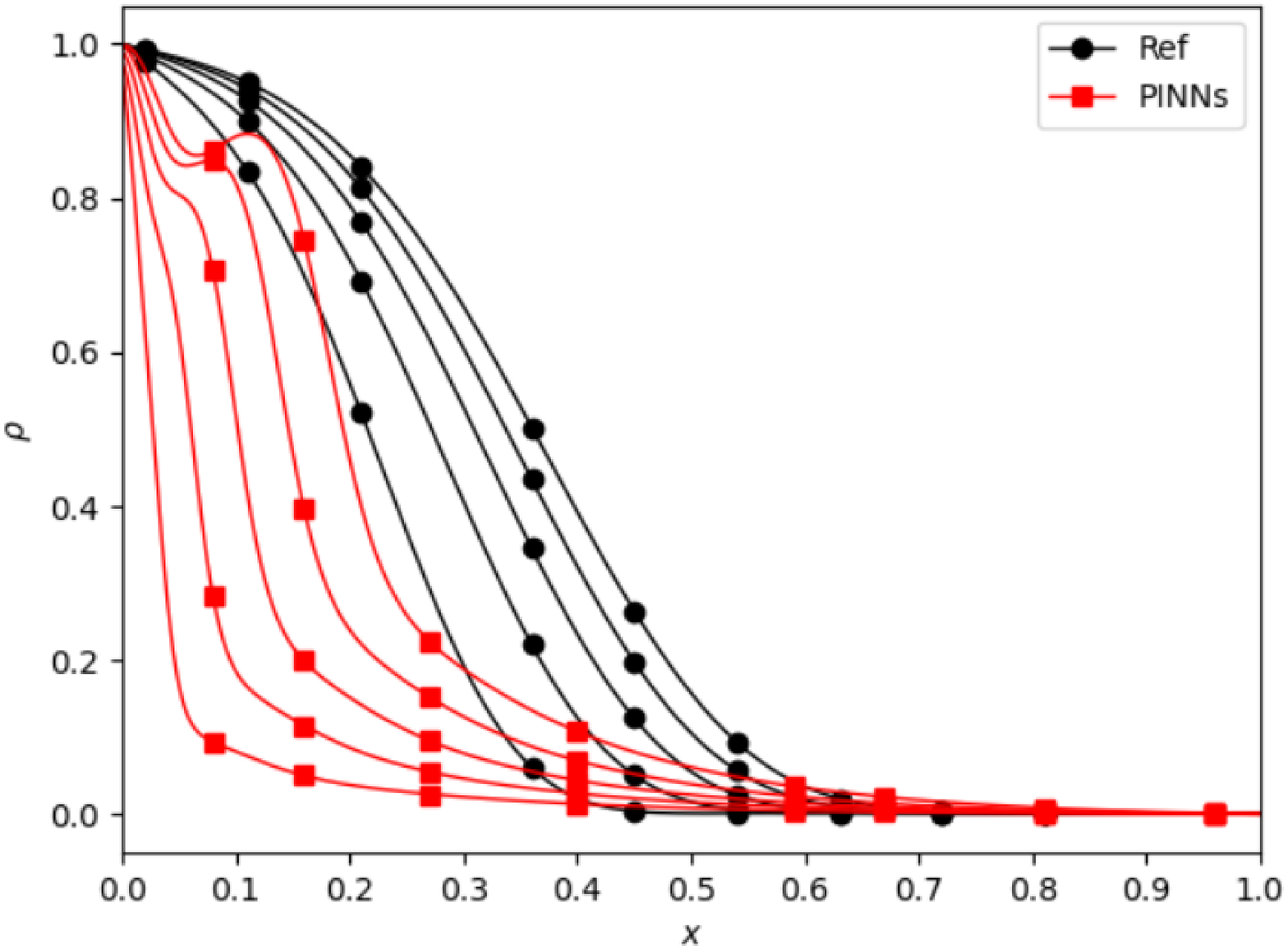}
		\end{minipage}
	}
	{
		\begin{minipage}{2.2in}
			\centering
			\includegraphics[width=2.2in]{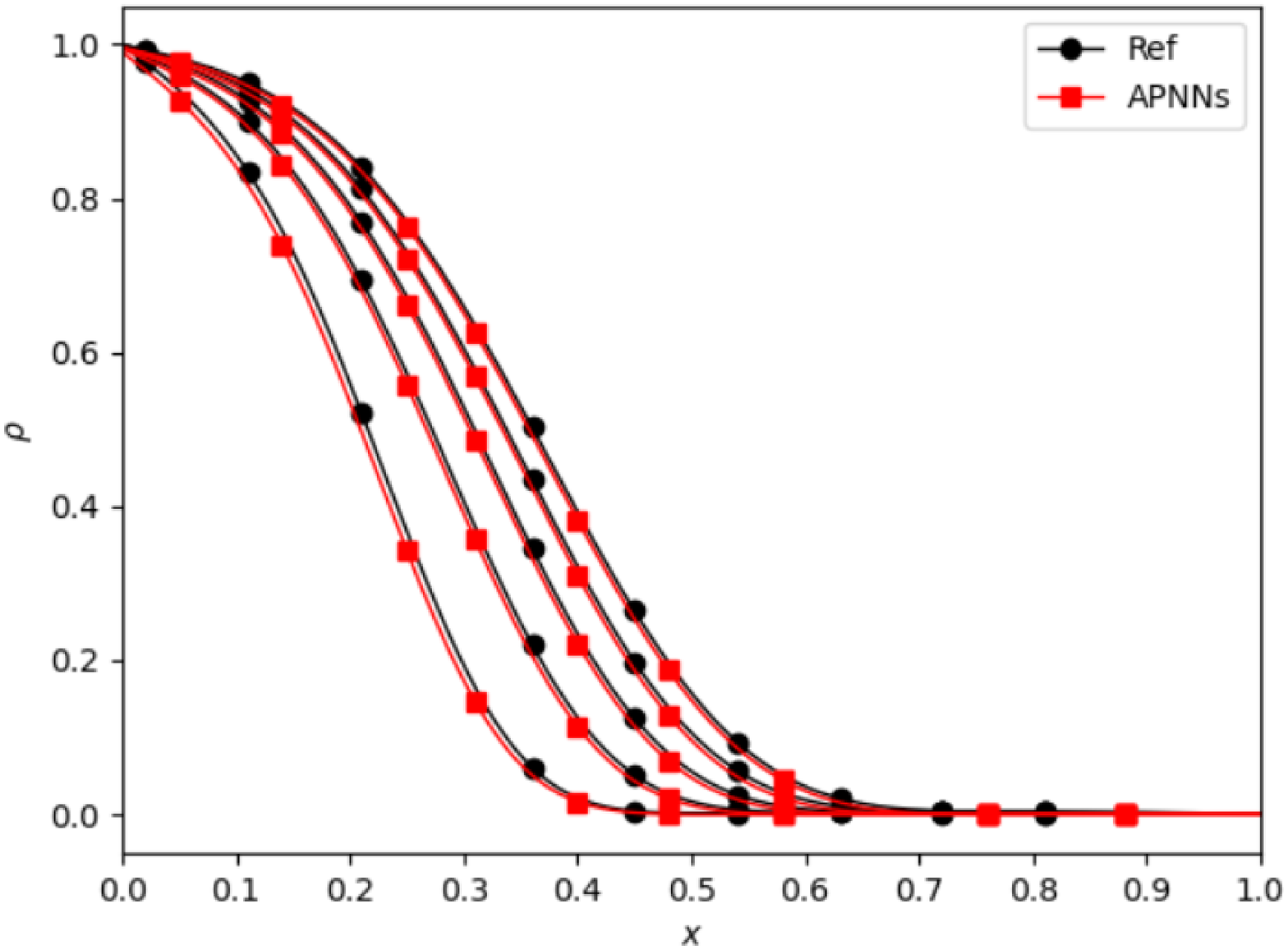}
		\end{minipage}
	}
	\caption{Intermediate regime with $\epsilon=10^{-2}$. The density $\rho$ at times $t=0.2, 0.4, 0.6, 0.8, 1.0$. (Left) Ref v.s. PINNs. (Right) Ref v.s. APNNs.}
	\label{figure5.2}
\end{figure}

We solve the 1D linear radiative transfer equation in the intermediate regime with $\epsilon=10^{-2}$.
\[\left\{
\begin{aligned}
& \epsilon \partial_t I + \mu \partial_x I =\frac{1+(10x)^2}{\epsilon}\left(\left \langle I \right \rangle - I\right) && x\in[0,1], \ t\in[0,1], \ \mu \in[-1,1], \\
& I(t,0,\mu >0)=1 && t\in[0,1], \ \mu \in(0,1], \\
& I(t,1,\mu <0)=0 && t\in[0,1], \ \mu \in[-1,0), \\
& I(0,x,\mu )=0  &&  x\in[0,1], \ \mu \in[-1,1],
\end{aligned}
\right.
\]
where $c=1$ and $\sigma=1+(10x)^2$.
We adopt the same settings as the previous example except that $N_{int}=8192$, $N_{sb}=N_{tb}=6144$, the activation function of the output layer for $\rho_{\theta_{2}}^{nn}$ is $\sigma_{\rho}^{o}(X)=X$,
and plot the reference solution and the prediction obtained by PINNs and APNNs in Figure~\ref{figure5.2}.
The reference solutions are obtained by the UGKS in \cite{sun2015asymptotic}.
The $L^2$-norm errors are presented in Table~\ref{Table5.2}, which indicates that the results of APNNs agree well with the reference solutions whereas PINNs have poor performance.

%
\begin{table}[!htbp]
	\centering
	\caption{Intermediate regime with $\epsilon=10^{-2}$. The errors of PINNs and APNNs.}
	\label{Table5.2}
	\scalebox{0.96}{
	\begin{tabular}{llllll}
		\hline\noalign{\smallskip}
		$L_{\text{error}}^2(\rho)$ &$t=0.2$ & $t=0.4$ & $t=0.6$ & $t=0.8$ & $t=1.0$ \\
		\noalign{\smallskip}\hline\noalign{\smallskip}
		PINNs &8.06e-01    &7.06e-01      &6.07e-01   &5.10e-01 &4.14e-01  \\
		APNNs &2.95e-02    &1.77e-02     &1.48e-02   &1.35e-02 &1.22e-02  \\
		\noalign{\smallskip}\hline
	\end{tabular}}
\end{table}
\subsection{The stationary nonlinear GRTEs}
\begin{figure}[!htbp]
	\centering
	{
		\begin{minipage}{2.2in}
			\centering
			\includegraphics[width=2.2in]{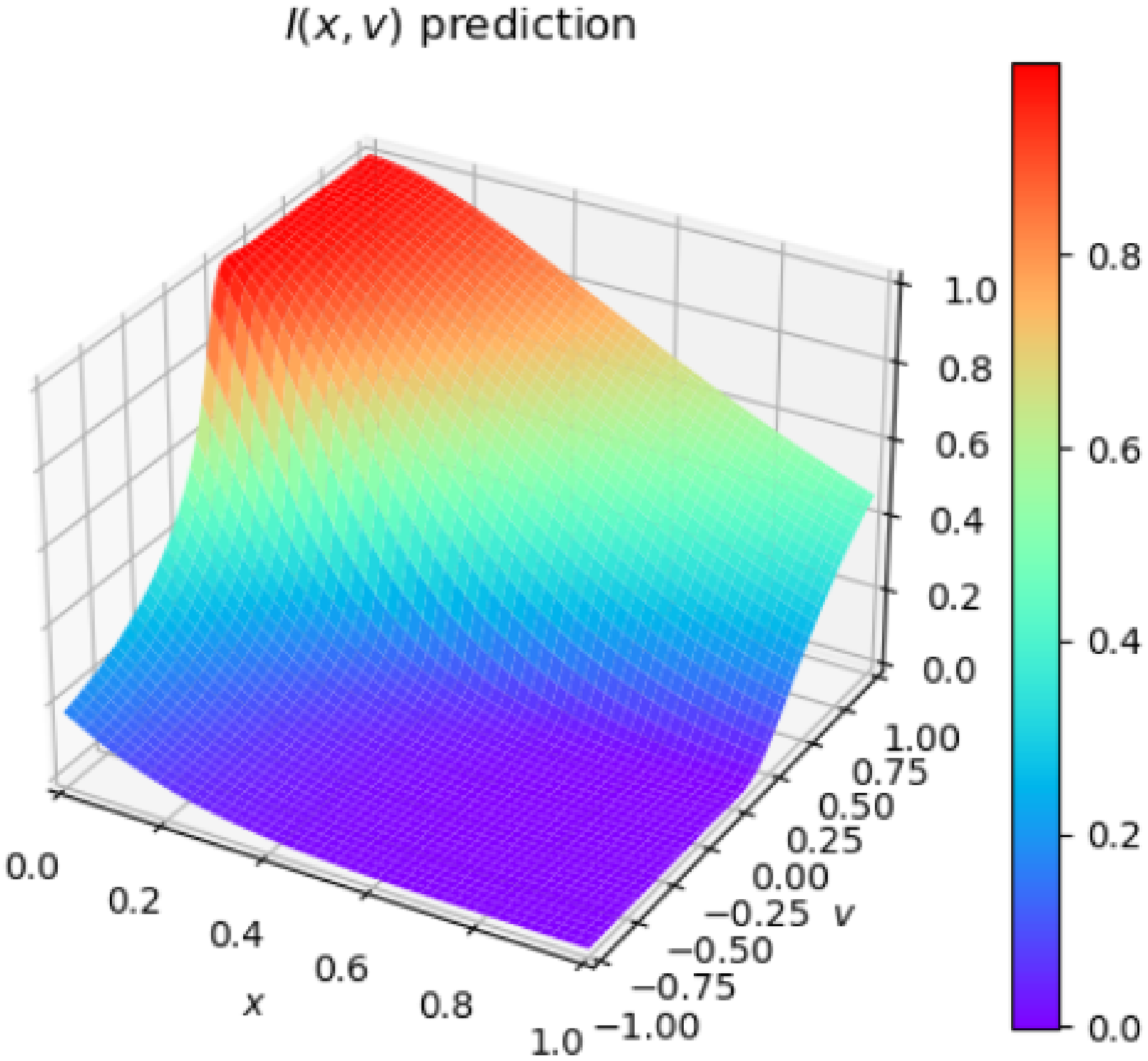}
		\end{minipage}
	}
	{
		\begin{minipage}{2.2in}
			\centering
			\includegraphics[width=2.2in]{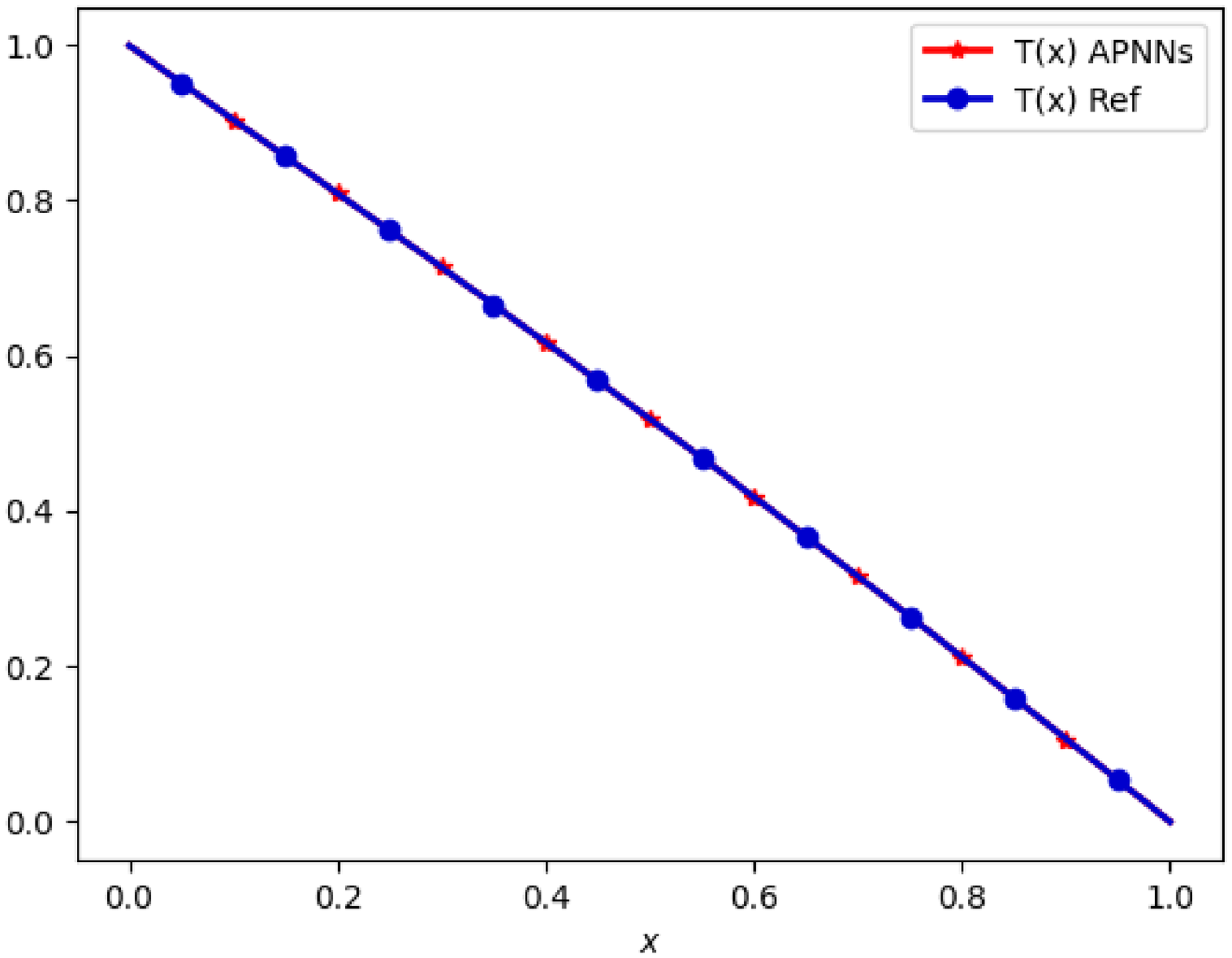}
		\end{minipage}
	}
	\caption{Kinetic regime with $\epsilon=1$. (Left) $I(x,v)$ of APNNs. (Right) $T(x)$:  Ref v.s. APNNs.}
	\label{figure5.3}
\end{figure}

\begin{figure}[!htbp]
	\centering
	{
		\begin{minipage}{2.2in}
			\centering
			\includegraphics[width=2.2in]{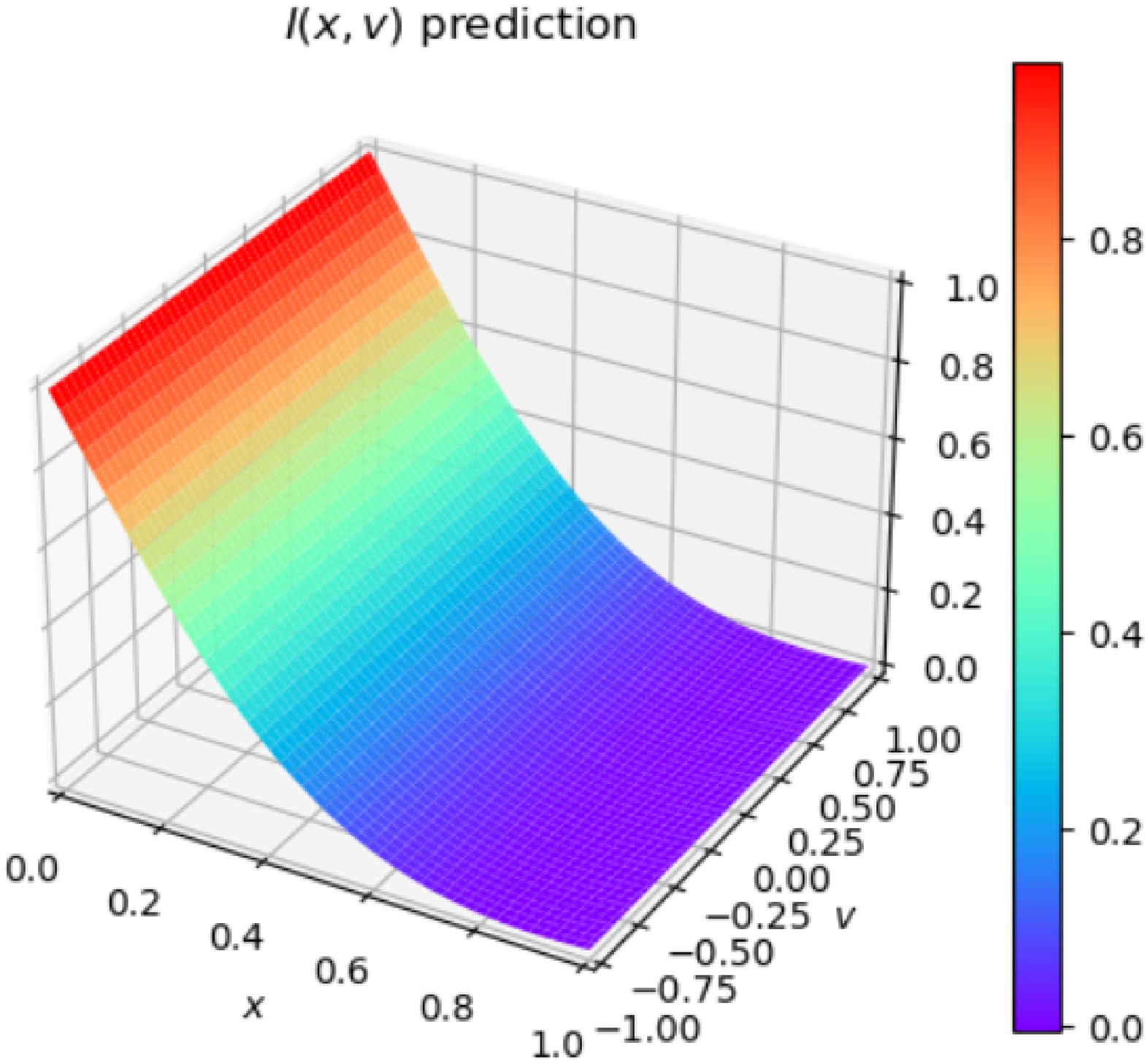}
		\end{minipage}
	}
	{
		\begin{minipage}{2.2in}
			\centering
			\includegraphics[width=2.2in]{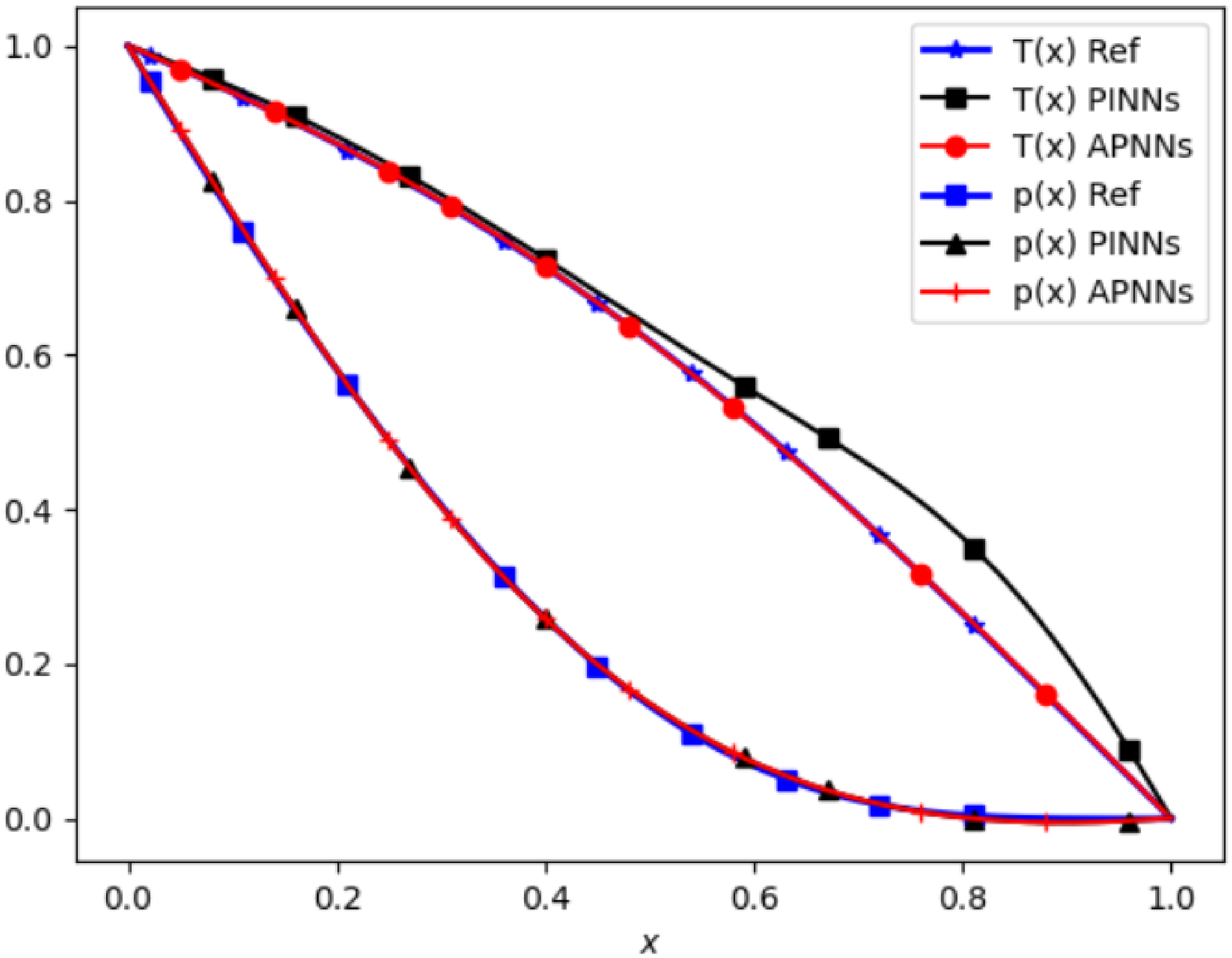}
		\end{minipage}
	}
	\caption{Diffusion regime with $\epsilon=10^{-3}$. (Left) $I$ of APNNs. (Right)  $T$ and $\rho$: Ref v.s. PINNs v.s. APNNs.}
	\label{figure5.4}
\end{figure}

We compare the performance between PINNs and APNNs for the 1D steady nonlinear radiative transfer equation \cite{lu2022solving}, i.e.,
\[
\left\{
\begin{aligned}
& \epsilon v \partial_x I(x,v)=\sigma\left(\frac{1}{2}acT^4(x)-I(x,v)\right) &&  x\in[0,1], \ v\in[-1,1], \\
& \epsilon^2 \partial_{xx}T(x)=-\sigma\left(acT ^4(x)-2\langle I(x,v) \rangle \right) && x\in[0,1], \ v\in[-1,1],\\
& I(0,v>0)=1, \qquad I(1,v<0)=0, && \\
& T(0)=1, \qquad \qquad  T(1)=0, &&
\end{aligned}
\right.
\]
where $a=c=\sigma=1$.
We set $I^{nn}_{\theta_1}=[2,50,50,50,50,1]$, $T^{nn}_{\theta_{2}}=[1,50,50,50,50,1]$ for PINNs, and $g^{nn}_{\theta_1}=[2,50,50,50,50,1]$, $(\rho,  T)^{nn}_{\theta_{2}}=[1,50,50,50,50,2]$ for APNNs.
No labeled data is used, and $N_{int}=4800$, $N_{sb}=3072$, $\lambda_{1}=\lambda_{3}=1$.
The activation functions of the output layers are all unit functions.
Figure~\ref{figure5.3} shows that the prediction obtained by APNNs agrees well with the reference solution in the kinetic regime with $\epsilon = 1$.
For the diffusion regime with $\epsilon = 10^{-3}$, Figure~\ref{figure5.4} indicates that the material temperature $T$ predicted by APNNs is more accurate than PINNs.

%
\subsection{The time-dependent nonlinear GRTEs}
We apply APNNs, MD-APNNs with a few labels, and Data-driven networks (using labeled data only) to the simulation of non-stationary nonlinear GRTEs.
The reference solutions are computed by the UGKS in \cite{sun2015asymptotic}. We set $g^{nn}_{\theta_1}=[3,50,50,50,50,1]$, $(\rho,  T)^{nn}_{\theta_{2}}=[2,50,50,50,50,2]$, and $N_{int}=16384$, $ N_{sb}=N_{sb1}+N_{sb2}=N_{tb}=12288$.
The activation functions of the output layers are set as $\sigma_g^o(X)=X$, $\sigma_{(\rho, T)}^o(X)=e^{-X}$.
For comparison, we plot the material temperature $T_e = T$, and the radiation temperature $T_r = \left(\frac{1}{ac}\int_{S^2}I\text{d}\Omega \right)^{\frac{1}{4}}$.

\subsubsection{Problem \uppercase\expandafter{\romannumeral1}}
We solve the 1D time-dependent GRTEs with temperature-independent opacity $\sigma=10\text{cm}^{-1}$ and heat capacity $C_v=0.01\text{GJ}/\text{cm}^3\text{KeV}$  on a slab of length $0.25$cm which is initially at equilibrium at $1$keV,
with the reflection condition and incident Planckian source condition on the left and right boundaries, respectively.
\begin{equation*}
\left\{
\begin{aligned}
& \frac{\epsilon^2}{c} \frac{\partial I}{\partial t}+\epsilon \mu \frac{\partial I}{\partial x}  =\sigma \left(\frac{1}{2}acT^4-I\right) && x\in[0,0.25], \ t\in[0,1], \ \mu\in[-1,1], \\
& \epsilon^2 C_v\frac{\partial T}{\partial t}=\sigma \left(2 \left \langle I \right \rangle -acT^4\right) && x\in[0,0.25], \ t\in[0,1], \ \mu\in[-1,1], \\
& I(t,0,\mu>0)=I(t,0,-\mu) && t\in[0,1], \ \mu\in(0,1], \\
& I(t, 0.25,\mu<0)=\frac{1}{2}ac(0.1)^4 && t\in[0,1], \ \mu\in[-1,0), \\
& I(0,x,\mu)=\frac{1}{2}acT(0,x)^4, \quad T(0,x)=1  &&  x\in[0,0.25],\mu\in[-1,1].
\end{aligned}
\right.
\end{equation*}

The empiricial loss of the governing equation in this problem is
\begin{align*}
&L_{\text{MD-APNNs},1g}^{\epsilon,nn} = \frac{1}{N_{\text{int}}}\sum_{j=1}^{N_{\text{int}}} (|\frac{\epsilon^2}{c}\partial_t g_{\theta_1}^{nn}(t_j^{\text{int}},x_j^{\text{int}},\mu_j^{\text{int}}) + \epsilon \mu \partial_x g_{\theta_1}^{nn}(t_j^{\text{int}},x_j^{\text{int}},\mu_j^{\text{int}}) \notag \\
&-\epsilon \frac{1}{2}\sum_{m=1}^{N_s}\mu_m^{\text{int}} \partial_x g_{\theta_1}^{nn}(t_j^{\text{int}},x_j^{\text{int}},\mu_m^{\text{int}})\omega_m + \sqrt{\sigma_0} \mu \partial_x \rho_{\theta_{21}}^{nn}(t_j^{\text{int}},x_j^{\text{int}}) + \sigma g_{\theta_1}^{nn}(t_j^{\text{int}},x_j^{\text{int}},\mu_j^{\text{int}}) |^2 \notag \\
&+|\frac{1}{c}\partial_t \rho_{\theta_{21}}^{nn}(t_j^{\text{int}},x_j^{\text{int}}) +\frac{1}{\sqrt{\sigma_0}} \frac{1}{2}\sum_{m=1}^{N_s}\mu_m^{\text{int}} \partial_x g_{\theta_1}^{nn}(t_j^{\text{int}},x_j^{\text{int}},\mu_m^{\text{int}})\omega_m + \frac{1}{2}C_v \partial_tT_{\theta_{22}}^{nn}(t_j^{\text{int}},x_j^{\text{int}}) |^2 \notag \\
&+|\epsilon^2 C_v \partial_tT_{\theta_{22}}^{nn}(t_j^{\text{int}},x_j^{\text{int}}) - \sigma\left(2\rho_{\theta_{21}}^{nn}(t_j^{\text{int}},x_j^{\text{int}})-ac{T^{nn}_{\theta_{22}}}^4(t_j^{\text{int}},x_j^{\text{int}})\right)|^2 ).
\end{align*}

The empiricial loss function of the boundary conditions is
\begin{align*}
L_{\text{MD-APNNs},1b}^{\epsilon,nn} =&\lambda_{1}(1)\frac{1}{N_{\text{sb1}}}\sum_{j=1}^{N_{\text{sb1}}} |\partial_x\rho_{\theta_{21}}^{nn}(t_j^{\text{sb1}},0) |^2 \\
& + \lambda_{1}(2)\frac{1}{N_{\text{sb1}}}\sum_{j=1}^{N_{\text{sb1}}} |g_{\theta_1}^{nn}(t_j^{\text{sb1}},0,\mu_j^{\text{sb1}}>0)-g_{\theta_1}^{nn}(t_j^{\text{sb1}},0,-\mu_j^{\text{sb1}})|^2   \\
&+\lambda_{1}(3)\frac{1}{N_{\text{sb2}}}\sum_{j=1}^{N_{\text{sb2}}} |\rho_{\theta_{21}}^{nn}(t_j^{\text{sb2}},0.25)-\frac{1}{2}ac(0.1)^4 |^2 \\ & + \lambda_{1}(4)\frac{1}{N_{\text{sb2}}}\sum_{j=1}^{N_{\text{sb2}}} |g_{\theta_1}^{nn}(t_j^{\text{sb2}},0.25,\mu_j^{\text{sb2}}<0)|^2.
\end{align*}

The empirical loss of the initial conditions is
\begin{align*}
L_{\text{MD-APNNs},1i}^{\epsilon,nn} &=\lambda_{2}(1)\frac{1}{N_{\text{tb}}}\sum_{j=1}^{N_{\text{tb}}} |\rho_{\theta_{21}}^{nn}(0,x_j^{\text{tb}})-\frac{1}{2}acT_{\theta_{22}}^{nn}(0,x_j^{\text{tb}})^4|^2  \notag \\
&+\lambda_{2}(2)\frac{1}{N_{\text{tb}}}\sum_{j=1}^{N_{\text{tb}}} |g_{\theta_1}^{nn}(0,x_j^{\text{tb}},\mu_j^{\text{tb}})|^2
+\lambda_{2}(3)\frac{1}{N_{\text{tb}}}\sum_{j=1}^{N_{\text{tb}}}|T_{\theta_{22}}^{nn}(0,x_j^{\text{tb}})-1|^2.
\end{align*}

The empiricial loss of the conservation law is
\begin{align*}
L_{\text{MD-APNNs},1c}^{\epsilon,nn} =\lambda_3 \frac{1}{N_{\text{int}}}\sum_{j=1}^{N_{\text{int}}}\left|\frac{1}{2}\sum_{m=1}^{N_s} g_{\theta_1}^{nn}(t_j^{\text{int}},x_j^{\text{int}},\mu_m^{\text{int}})\omega_m \right|^2.
\end{align*}

The empiricial loss of the data regularization term is
\begin{align*}
L_{\text{MD-APNNs},1l}^{\epsilon,nn} =\lambda_0\frac{1}{N_{0}}\sum_{i=1}^{N_{0}}|T_{\theta_{22}}^{nn}(t_i,x_i)-T^*(t_i,x_i)|^2.
\end{align*}

Finally, the total empirical loss function for training in $\text{MD-APNNs,1}$ is as follows
\begin{align*}
L_{\text{MD-APNNs,1}}^{\epsilon,nn}=\sum_{\mathrm{i}=\{g,b,i,c,l\}} L_{\text{MD-APNNs,1i}}^{\epsilon,nn}.
\end{align*}

Here and hereafter, without specifying, we set $\lambda_i=1$ (also $\lambda_i(j)=1$) for all $i,j$.

\begin{figure}[!htbp]
	\centering
	{
		\begin{minipage}{2.05in}
			\centering
			\includegraphics[width=2.05in]{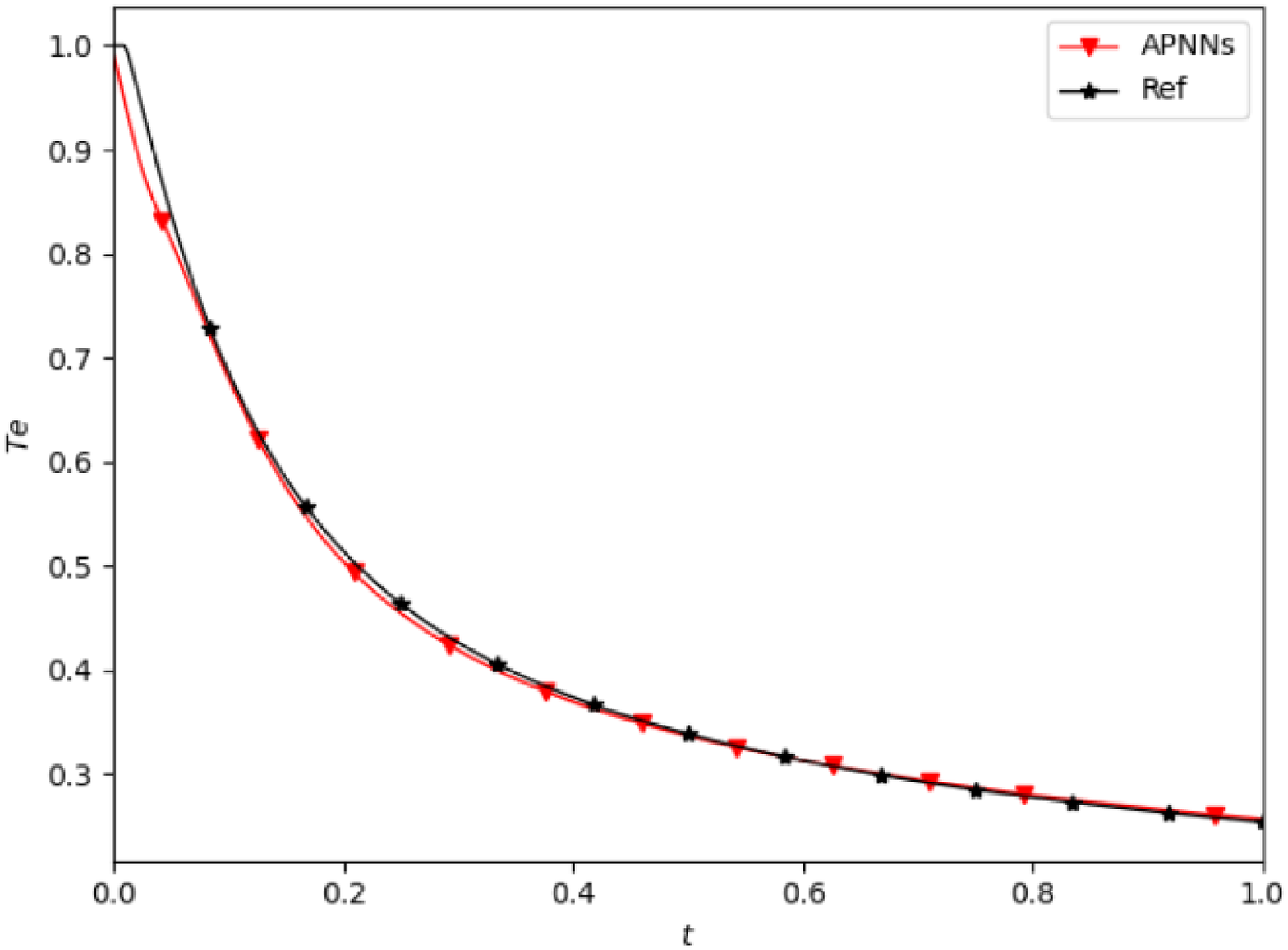}
		\end{minipage}
	}
	{
		\begin{minipage}{2.05in}
			\centering
			\includegraphics[width=2.05in]{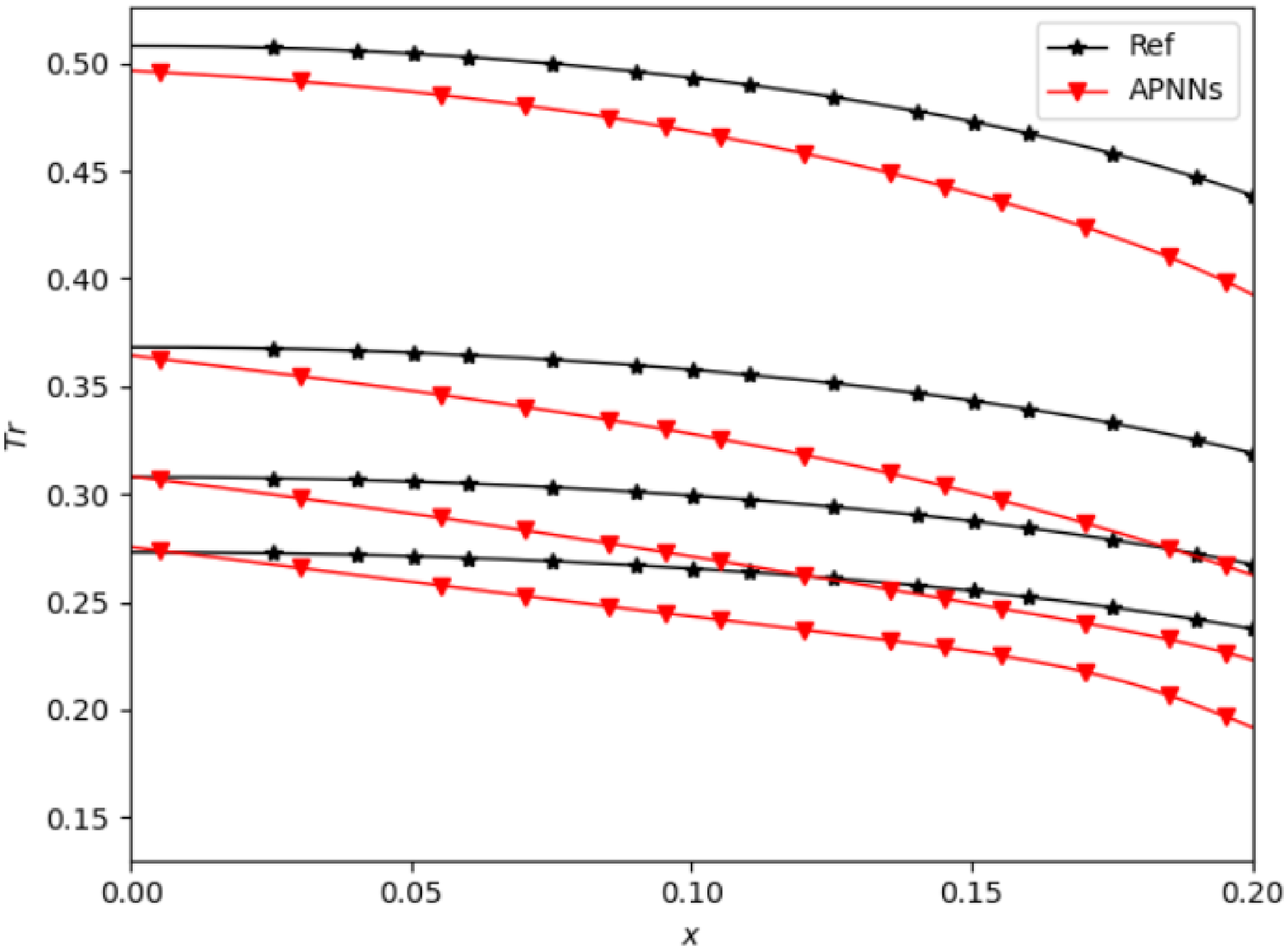}
		\end{minipage}
	}
	\caption{Kinetic regime with $\epsilon=1$. Ref v.s. APNNs. (Left) The material temperature $T_e$ at $x=0.0025$. (Right) The radiation temperature $T_r$ at times $t=0.2, 0.4, 0.6, 0.8$. $\lambda_{1}(1)=\lambda_0=0$, $\lambda_{1}(3)=10$. }
	\label{figure5.5}
\end{figure}

\begin{figure}[!htbp]
	\centering
	{
		\begin{minipage}{2.05in}
			\centering
			\includegraphics[width=2.05in]{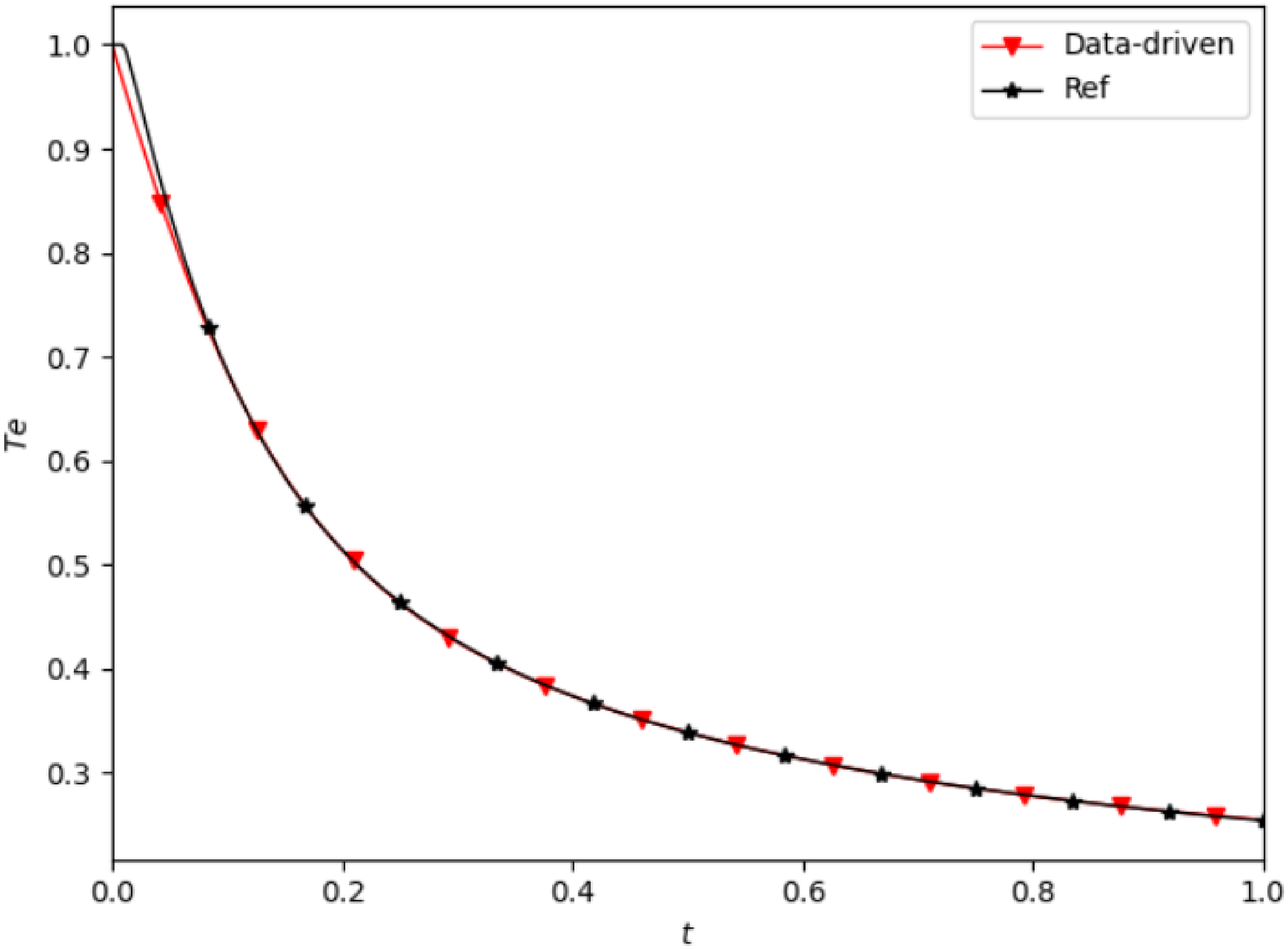}
		\end{minipage}
	}
	{
		\begin{minipage}{2.05in}
			\centering
			\includegraphics[width=2.05in]{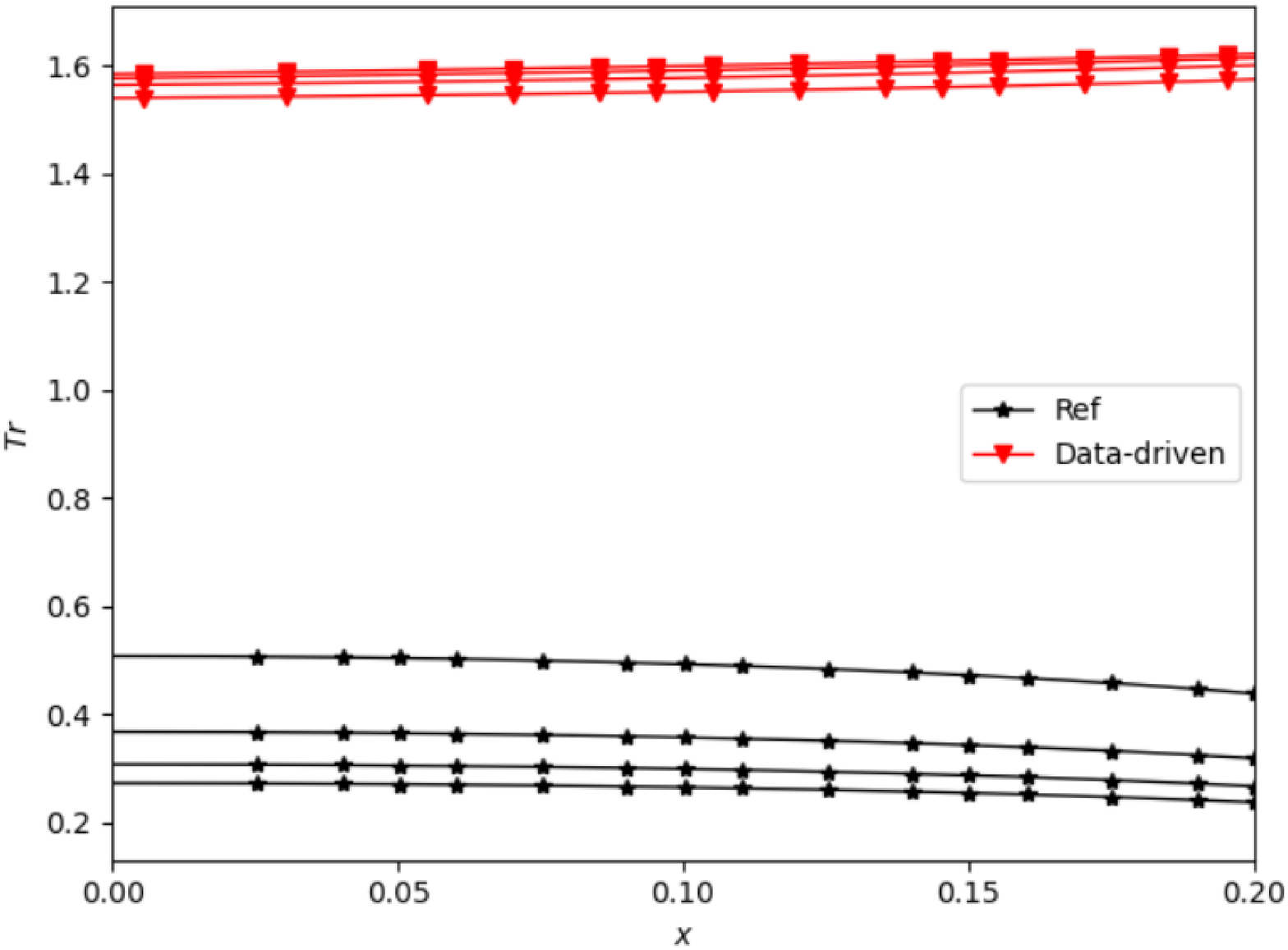}
		\end{minipage}
	}
	\caption{Kinetic regime with $\epsilon=1$. Ref v.s. Data-driven.  (Left) $T_e$ at $x=0.0025$. (Right) $T_r$ at times $t=0.2, 0.4, 0.6, 0.8$. $\lambda_0=1$ and other weight hyperparameters are 0.}
	\label{figure5.6}
\end{figure}

\begin{figure}[!htbp]
	\centering
	{
		\begin{minipage}{2.05in}
			\centering
			\includegraphics[width=2.05in]{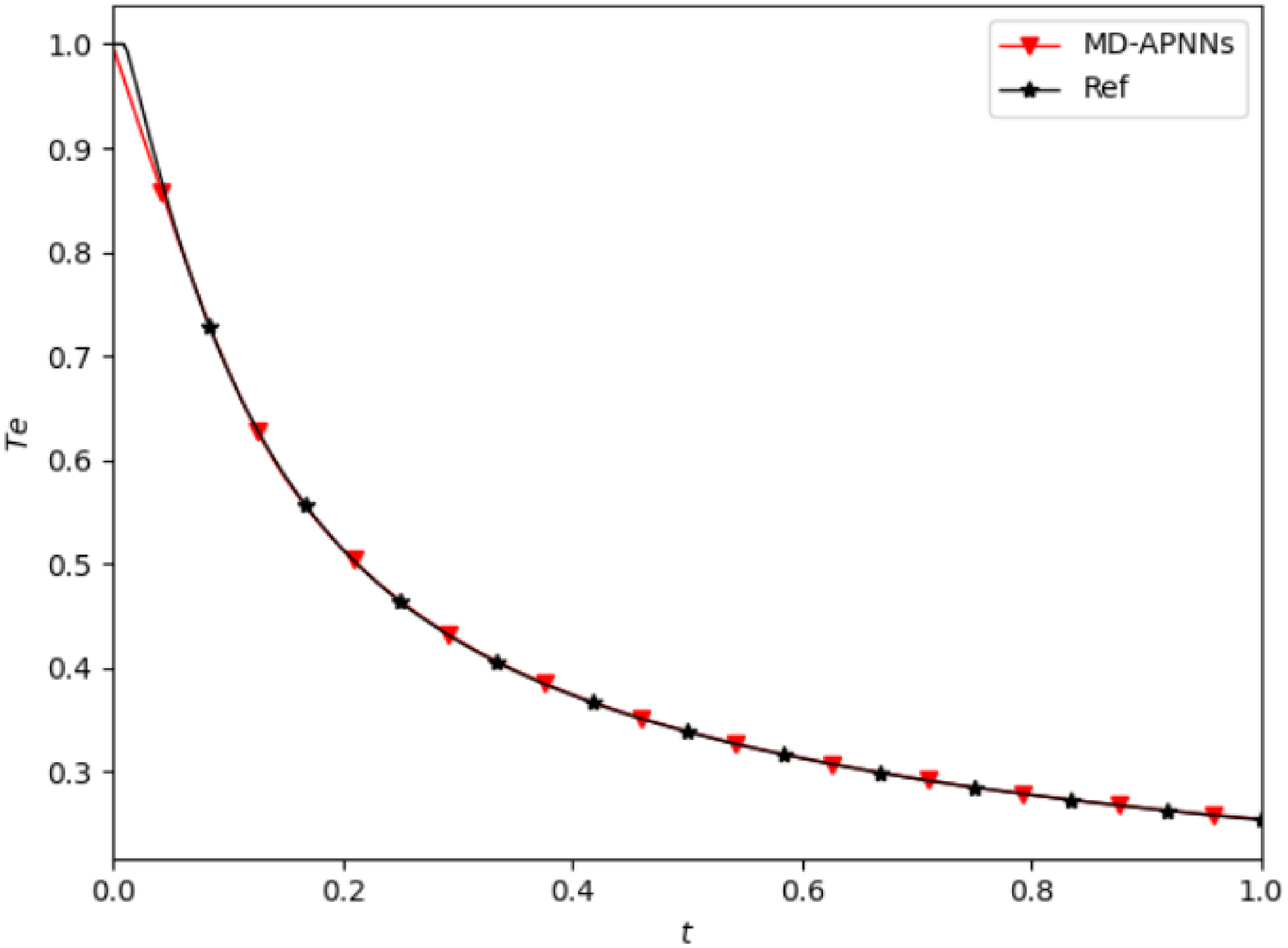}
		\end{minipage}
	}
	{
		\begin{minipage}{2.05in}
			\centering
			\includegraphics[width=2.05in]{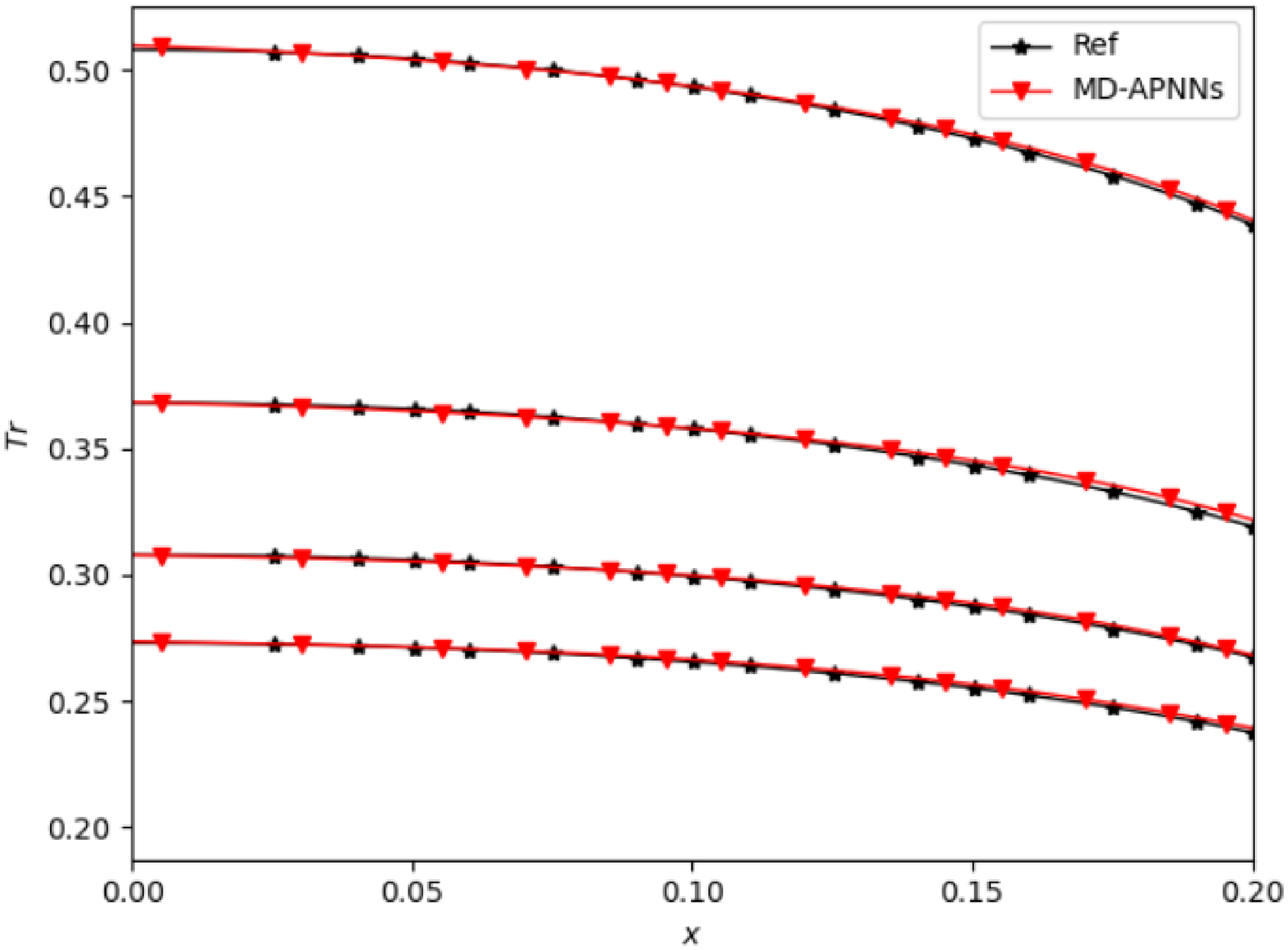}
		\end{minipage}
	}
	\caption{Kinetic regime with $\epsilon=1$. Ref v.s. MD-APNNs. (Left) $T_e$ at $x=0.0025$. (Right) $T_r$ at times $t=0.2, 0.4, 0.6, 0.8$. $\lambda_{1}(1)=0$, $\lambda_0=10$.}
	\label{figure5.7}
\end{figure}

{\bf (1) Kinetic regime with $\epsilon=1$.} We plot the prediction obtained by APNNs, Data-driven networks and MD-APNNs, together with the reference solution in Figures~\ref{figure5.5}, \ref{figure5.6} and \ref{figure5.7}, respectively.
For MD-APNNs and Data-driven networks, we use $N_0=60$ labeled data of the material temperature $T$.
We observe that the performance of APNNs is not satisfying (see Figure~\ref{figure5.5}).
Data-driven networks using the labeled data of $T$ can generate a good approximation of $T_e$ (Figure~\ref{figure5.6} (left)), but the prediction of radiation temperature $T_r$ are far from the reference solution (Figure~\ref{figure5.6} (right)).
The MD-APNN method gives accurate prediction compared to the reference solution.
The $L^2$-norm errors of APNNs, Data-driven networks, and MD-APNNs are presented in Table~\ref{Table5.3}.
\begin{table}[!htbp]
	\centering
	\caption{Kinetic regime with $\epsilon=1$: the errors of $T_e$ and $T_r$ (at $t=0.2, 0.4, 0.6, 0.8$) for APNNs, Data-driven networks and MD-APNNs.}
	\label{Table5.3}
	\scalebox{0.9}{
	\begin{tabular}{llllll}
		\hline\noalign{\smallskip}
		$L^2$ error &$T_e$ & $T_r(t=0.2)$ & $T_r(t=0.4)$ & $T_r(t=0.6)$ & $T_r(t=0.8)$  \\
		\noalign{\smallskip}\hline\noalign{\smallskip}
		Data-driven &1.49e-02 &2.30e+00 &3.62e+00  &4.57e+00  &5.30e+00 \\
		\noalign{\smallskip}\hline\noalign{\smallskip}
		APNNs &2.70e-02    &8.51e-02  &1.22e-01  &1.25e-01  &1.30e-01  \\
		\noalign{\smallskip}\hline\noalign{\smallskip}
		MD-APNNs & 1.16e-02    & 3.21e-03  &4.54e-03  &4.51e-03  &1.20e-02  \\
		\noalign{\smallskip}\hline
	\end{tabular}}
\end{table}

\begin{figure}[!htbp]
	\centering
	{
		\begin{minipage}{2.2in}
			\centering
			\includegraphics[width=2.2in]{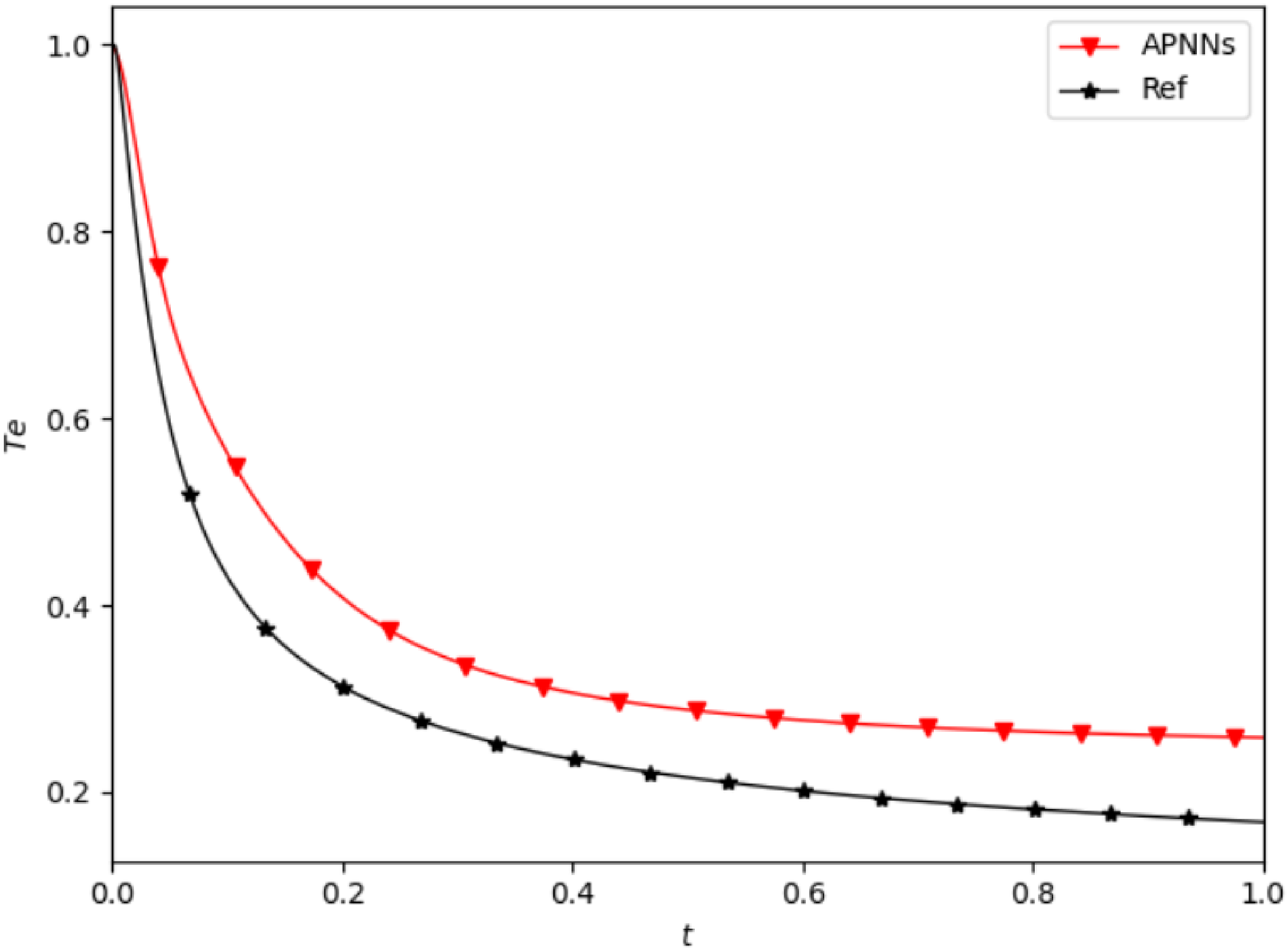}
		\end{minipage}
	}
	{
		\begin{minipage}{2.2in}
			\centering
			\includegraphics[width=2.2in]{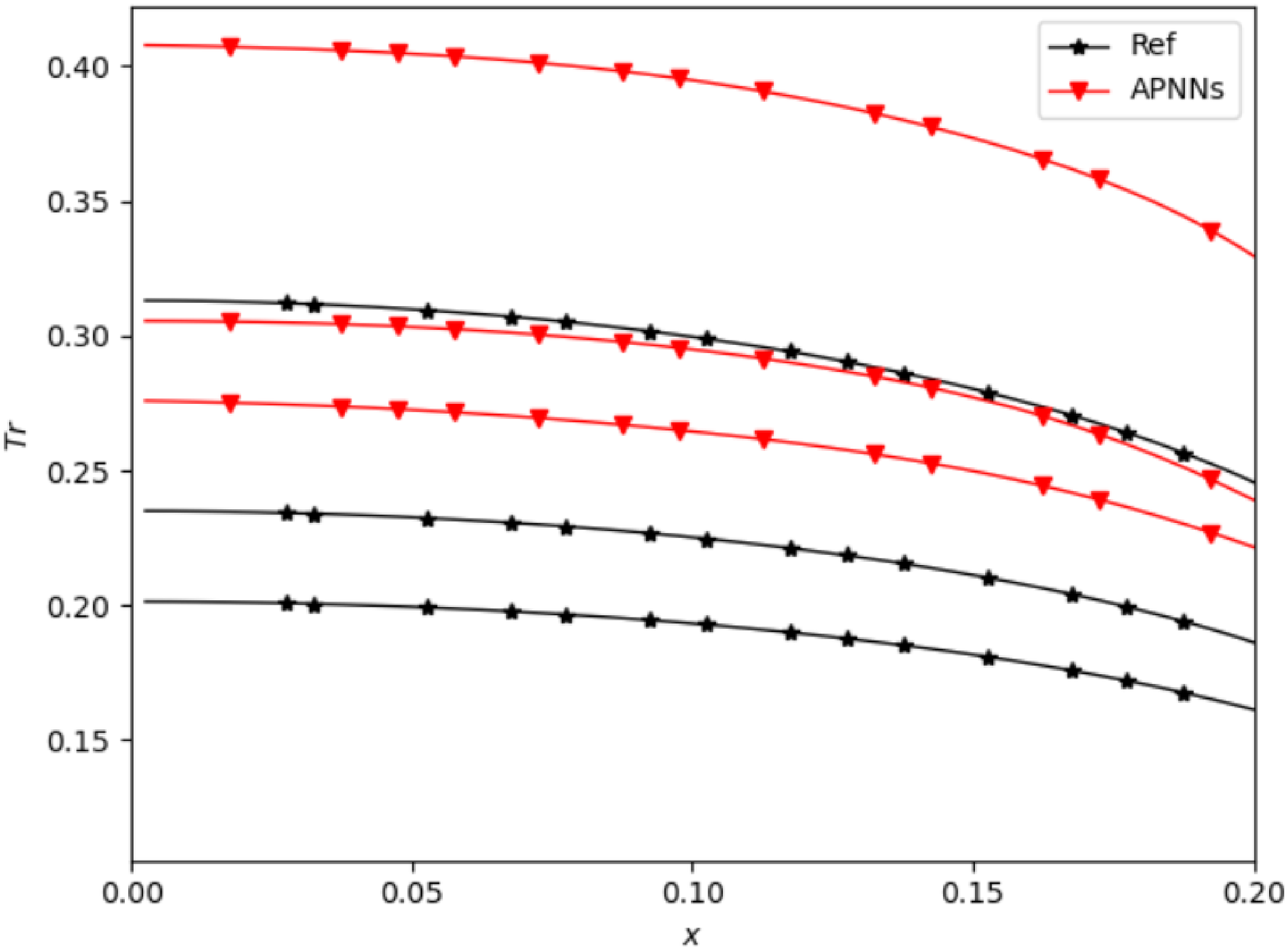}
		\end{minipage}
	}
	
	\caption{Diffusion regime with $\epsilon=10^{-6}$. Ref v.s. APNNs. (Left) $T_e$ at $x=0.0025$. (Right) $T_r$ at times $t=0.2, 0.4, 0.6$.  $\lambda_{1}(3)=10$, $\lambda_0=0$.}
	\label{figure5.8}
\end{figure}

\begin{figure}[!htbp]
	\centering
	{
		\begin{minipage}{2.2in}
			\centering
			\includegraphics[width=2.2in]{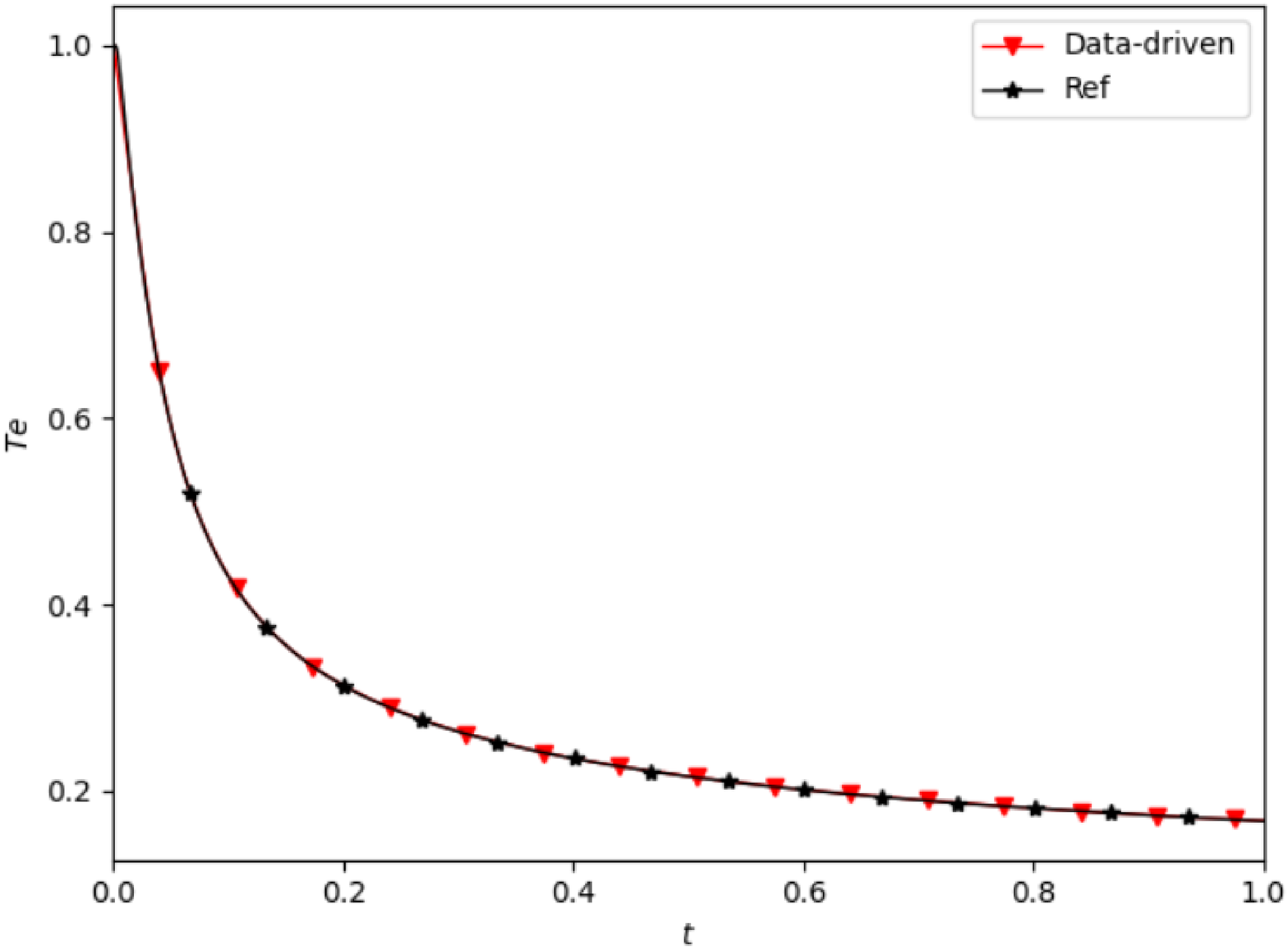}
		\end{minipage}
	}
	{
		\begin{minipage}{2.2in}
			\centering
			\includegraphics[width=2.2in]{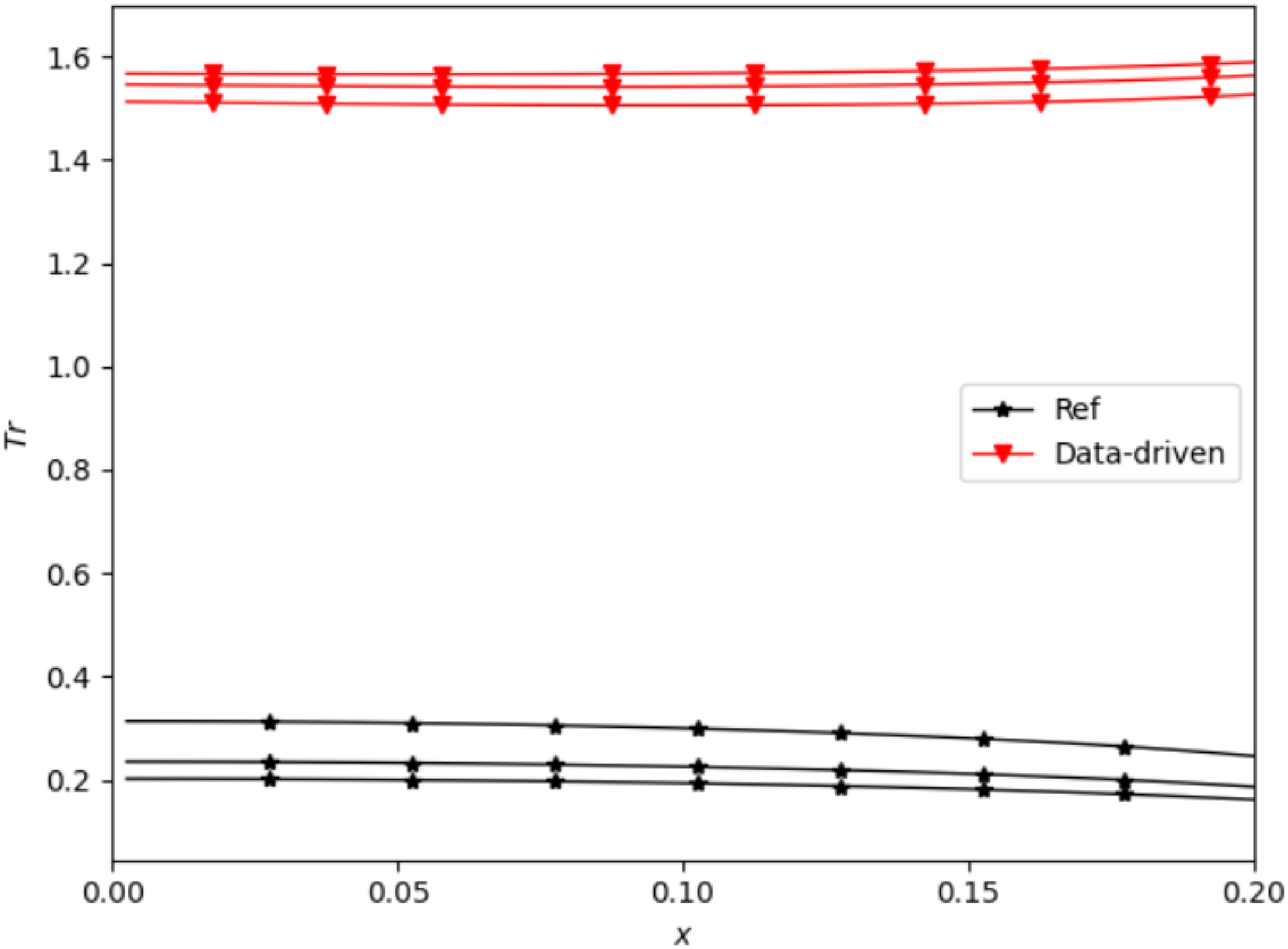}
		\end{minipage}
	}
	\caption{Diffusion regime with $\epsilon=10^{-6}$. Ref v.s. Data-driven. (Left) $T_e$ at $x=0.0025$. (Right) $T_r$ at times $t=0.2, 0.4, 0.6$. $\lambda_0=1$ and other weight hyperparameters are 0.}
	\label{figure5.9}
\end{figure}

\begin{figure}[!htbp]
	\centering
	{
		\begin{minipage}{2.2in}
			\centering
			\includegraphics[width=2.2in]{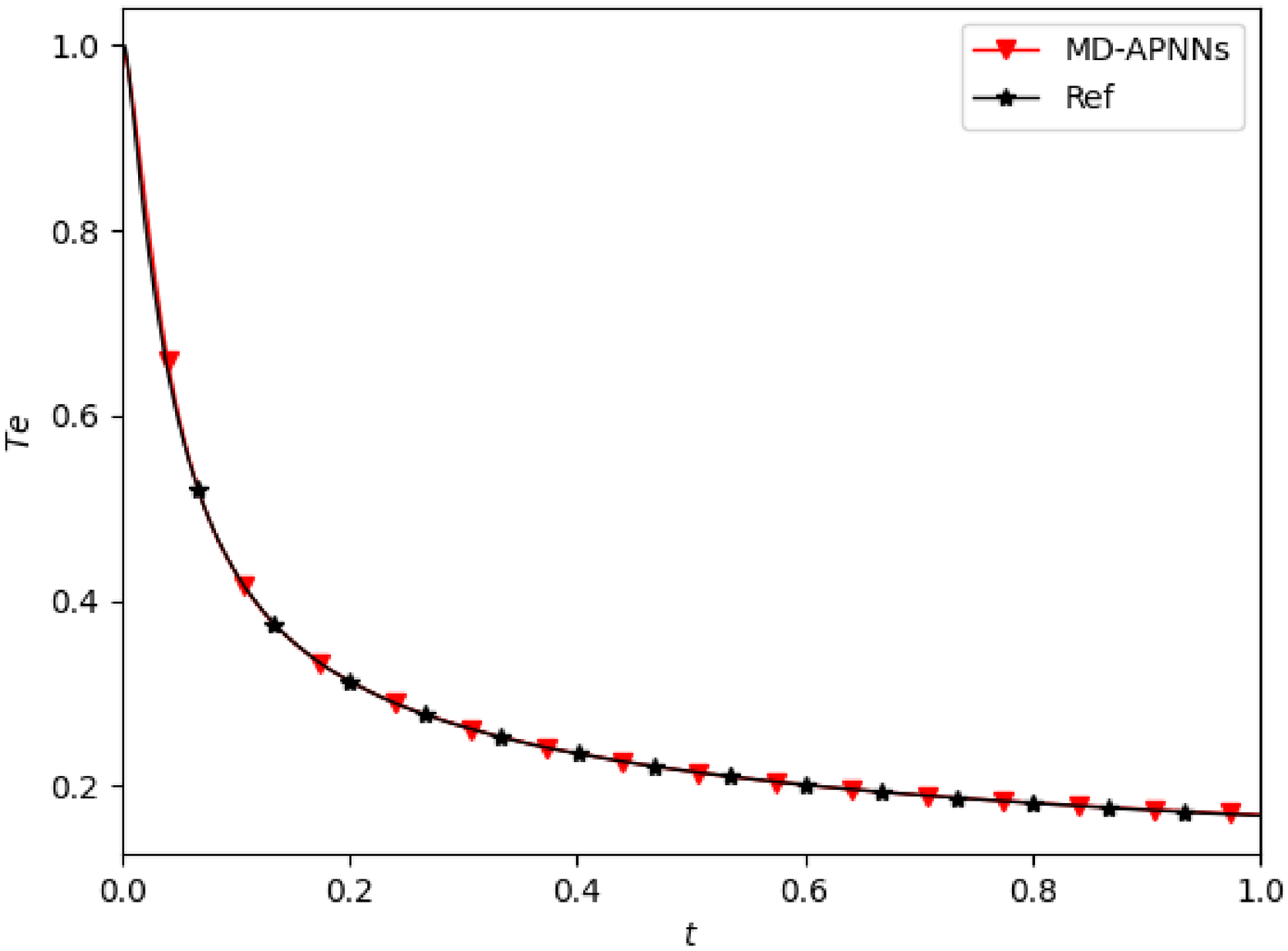}
		\end{minipage}
	}
	{
		\begin{minipage}{2.2in}
			\centering
			\includegraphics[width=2.2in]{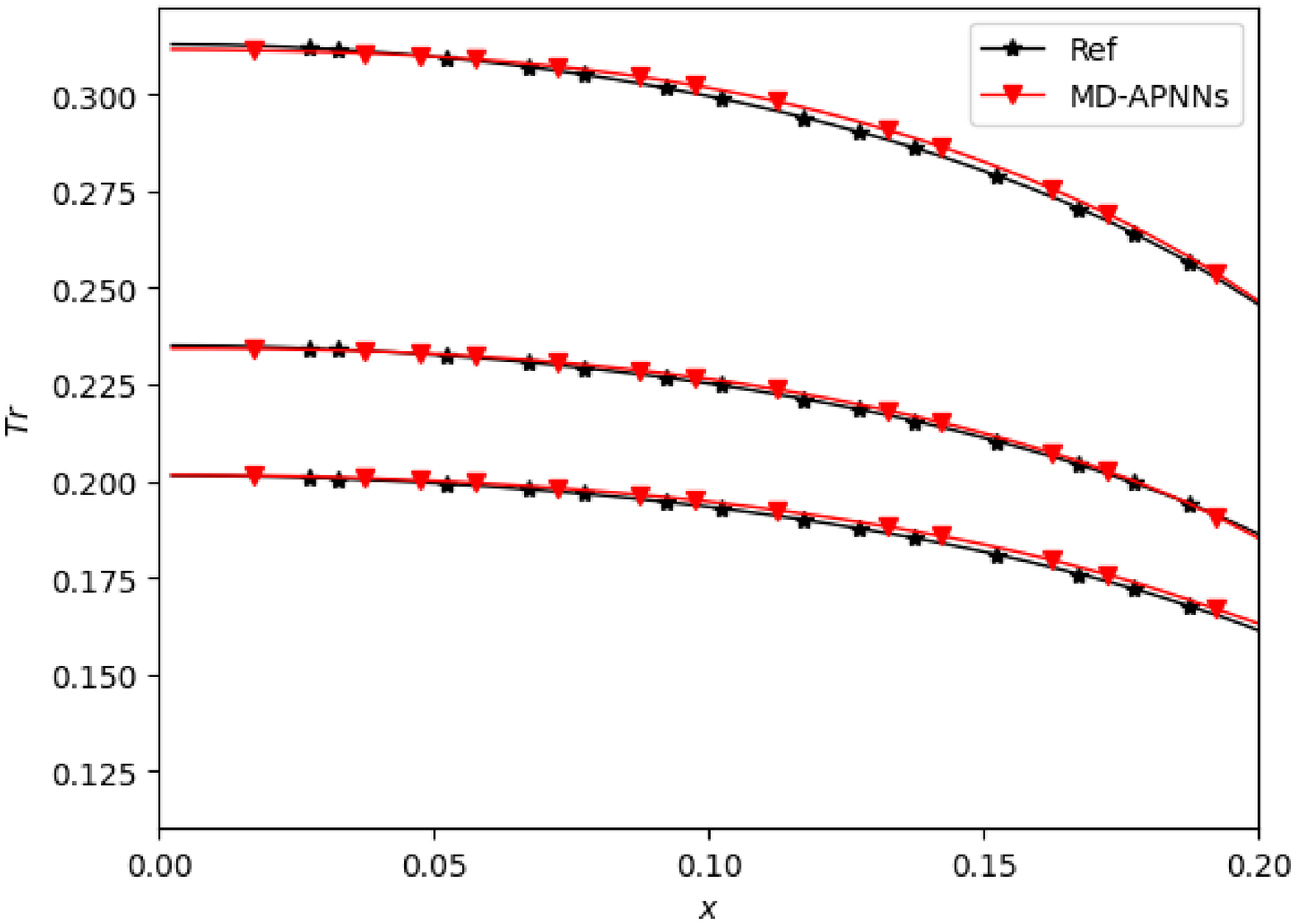}
		\end{minipage}
	}
	\caption{Diffusion regime with $\epsilon=10^{-6}$. Ref v.s. MD-APNNs. (Left) $T_e$ at $x=0.0025$. (Right) $T_r$ at times $t=0.2, 0.4, 0.6$. $\lambda_1(3)=10$.}
	\label{figure5.10}
\end{figure}

{\bf (2) Diffusion regime with $\epsilon=10^{-6}$.}
Now we use $N_0=100$ labeled data of the material temperature $T$ for Data-driven networks and MD-APNNs.
The prediction obtained by APNNs, Data-driven networks and MD-APNNs are plotted in Figures~\ref{figure5.8}, \ref{figure5.9} and \ref{figure5.10}, respectively. The $L^2$-norm errors of APNNs, Data-driven networks, and MD-APNNs are presented in Table~\ref{Table5.4}, which indicate that the MD-APNN method produces a better result compared to the APNN and Data-driven method.

\begin{table}[!htbp]
	\centering
	\caption{Diffusion regime with $\epsilon=10^{-6}$: the errors of $T_e$ and $T_r$ (at $t=0.2, 0.4, 0.6$) for APNNs, Data-driven networks and MD-APNNs.}
	\label{Table5.4}
	\scalebox{0.96}{
	\begin{tabular}{lllll}
		\hline\noalign{\smallskip}
		$L^2$ error &$T_e$ & $T_r(t=0.2)$ & $T_r(t=0.4)$ & $T_r(t=0.6)$  \\
		\noalign{\smallskip}\hline\noalign{\smallskip}
		Data-driven &7.33e-03 &4.44e+00 &6.38e+00  &7.70e+00  \\
		\noalign{\smallskip}\hline\noalign{\smallskip}
		APNNs &2.84e-01    &3.21e-01  &3.05e-01  &3.76e-01 \\
		\noalign{\smallskip}\hline\noalign{\smallskip}
		MD-APNNs &1.18e-02    &7.93e-03   &1.70e-02  &3.12e-02 \\
		\noalign{\smallskip}\hline
	\end{tabular}}
\end{table}

\subsubsection{Problem \uppercase\expandafter{\romannumeral2}}
We solve the 1D example \cite{fu2022asymptotic} with smooth initial conditions at the equilibrium and periodic boundary condition. Spatial domian is $[0, 2]$, time interval is $[0,0.5]$ and the parameters are set as $a=c=1, C_v=0.1, \sigma=10$.

\begin{equation*}
\left\{
\begin{aligned}
&\frac{\epsilon^2}{c}\partial_tI+\epsilon\mu\partial_xI=\sigma\left(\frac{1}{2}acT^4-I\right) && x\in[0,2], t\in[0,0.5], \mu\in [-1,1],\\
&\epsilon^2C_v\partial_tT=\sigma\left(\int_{-1}^{1}I \text{d}v-acT^4\right) && x\in[0,2], t\in[0,0.5],\\
&I(t,0,\mu)=I(t, 2 , \mu) && t\in[0,0.5],\mu\in [-1,1],\\
&I(0,x,\mu)=\frac{1}{2}acT(0,x)^4,\quad T(0,x)=\frac{3+\sin(\pi x)}{4} && x\in[0,2],\mu\in [-1,1].
\end{aligned}
\right.
\end{equation*}

The constructed empiricial loss functions of the governing equation, conservation law and data regularization term are the same to the previous example. The empiricial loss functions of the boundary and initial conditions are different and stated as follows
\begin{align*}
L_{\text{MD-APNNs},2b}^{\epsilon,nn} &=\lambda_{1}(1)\frac{1}{N_{\text{sb1}}}\sum_{j=1}^{N_{\text{sb1}}} |\rho_{\theta_{21}}^{nn}(t_j^{\text{sb1}},0)-\rho_{\theta_{21}}^{nn}(t_j^{\text{sb1}},2) |^2  \\ &+\lambda_{1}(2)\frac{1}{N_{\text{sb1}}}\sum_{j=1}^{N_{\text{sb1}}} |g_{\theta_1}^{nn}(t_j^{\text{sb1}},0,\mu_j^{\text{sb1}})-g_{\theta_1}^{nn}(t_j^{\text{sb1}},2,\mu_j^{\text{sb1}})|^2, \\
L_{\text{MD-APNNs},2i}^{\epsilon,nn}& =
\lambda_{2}(1)\frac{1}{N_{\text{tb}}}\sum_{j=1}^{N_{\text{tb}}} |\rho_{\theta_{21}}^{nn}(0,x_j^{\text{tb}})-\frac{1}{2}acT_{\theta_{22}}^{nn}(0,x_j^{\text{tb}})^4|^2 \\
&+\lambda_{2}(2)\frac{1}{N_{\text{tb}}}\sum_{j=1}^{N_{\text{tb}}} |g_{\theta_1}^{nn}(0,x_j^{\text{tb}},\mu_j^{\text{tb}})|^2 \\ &+\lambda_{2}(3)\frac{1}{N_{\text{tb}}}\sum_{j=1}^{N_{\text{tb}}}|T_{\theta_{22}}^{nn}(0,x_j^{\text{tb}})-\frac{3+\sin(\pi x_j^{\text{tb}} )}{4}|^2.
\end{align*}

Then the total emperical loss for $\text{MD-APNNs,2}$ is
\begin{align*}
L_{\text{MD-APNNs,2}}^{\epsilon,nn}=\sum_{\mathrm{i}=\{g,b,i,c,l\}} L_{\text{MD-APNNs,2i}}^{\epsilon,nn}.
\end{align*}

\begin{figure}[!htbp]
	\setlength{\abovecaptionskip}{0.cm}
	\setlength{\belowcaptionskip}{-0.cm}
	\centering

	{
		\begin{minipage}{2.2in}
			\centering
			\includegraphics[width=2.2in]{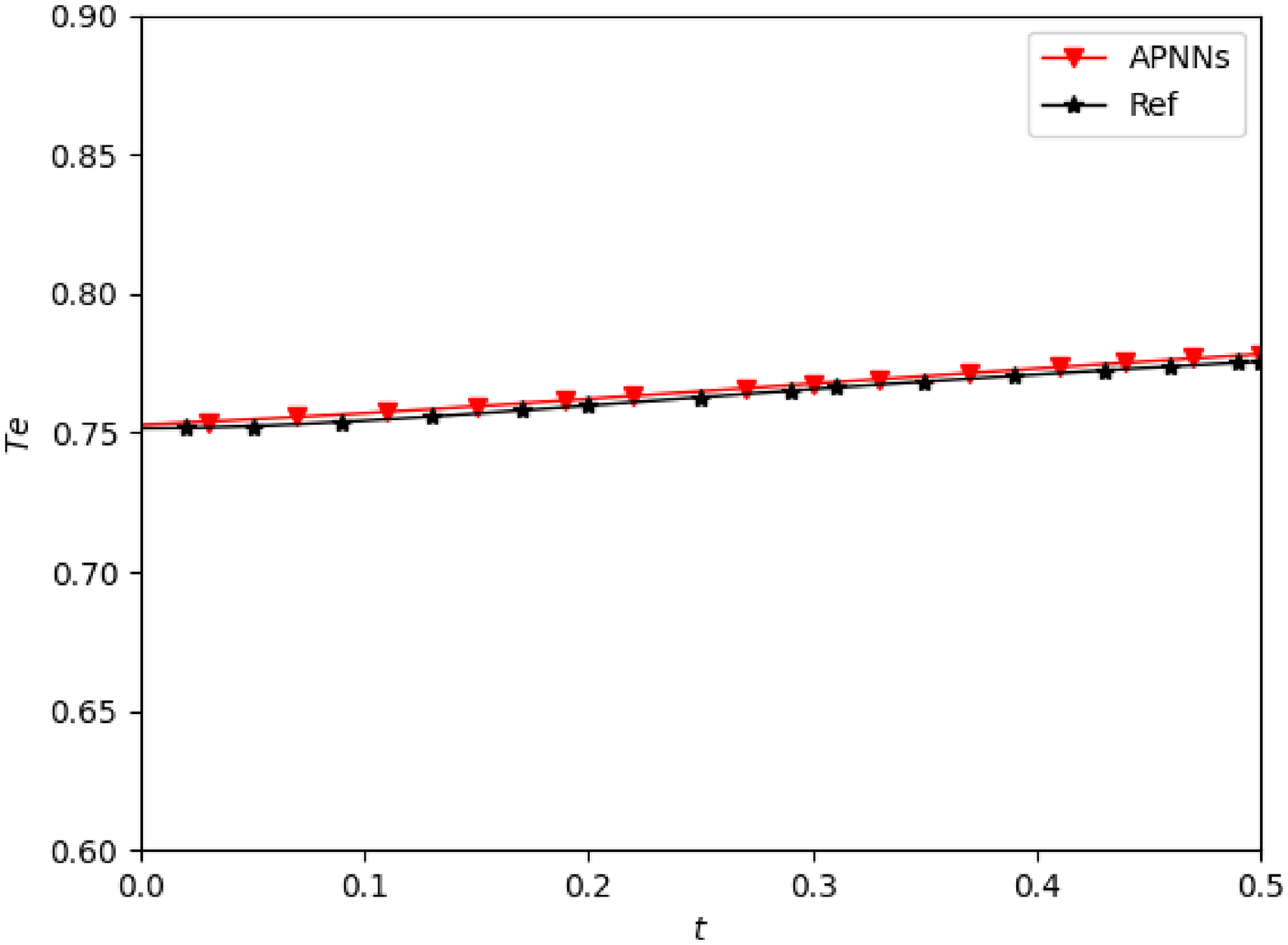}
		\end{minipage}
	}
	{
		\begin{minipage}{2.2in}
			\centering
			\includegraphics[width=2.2in]{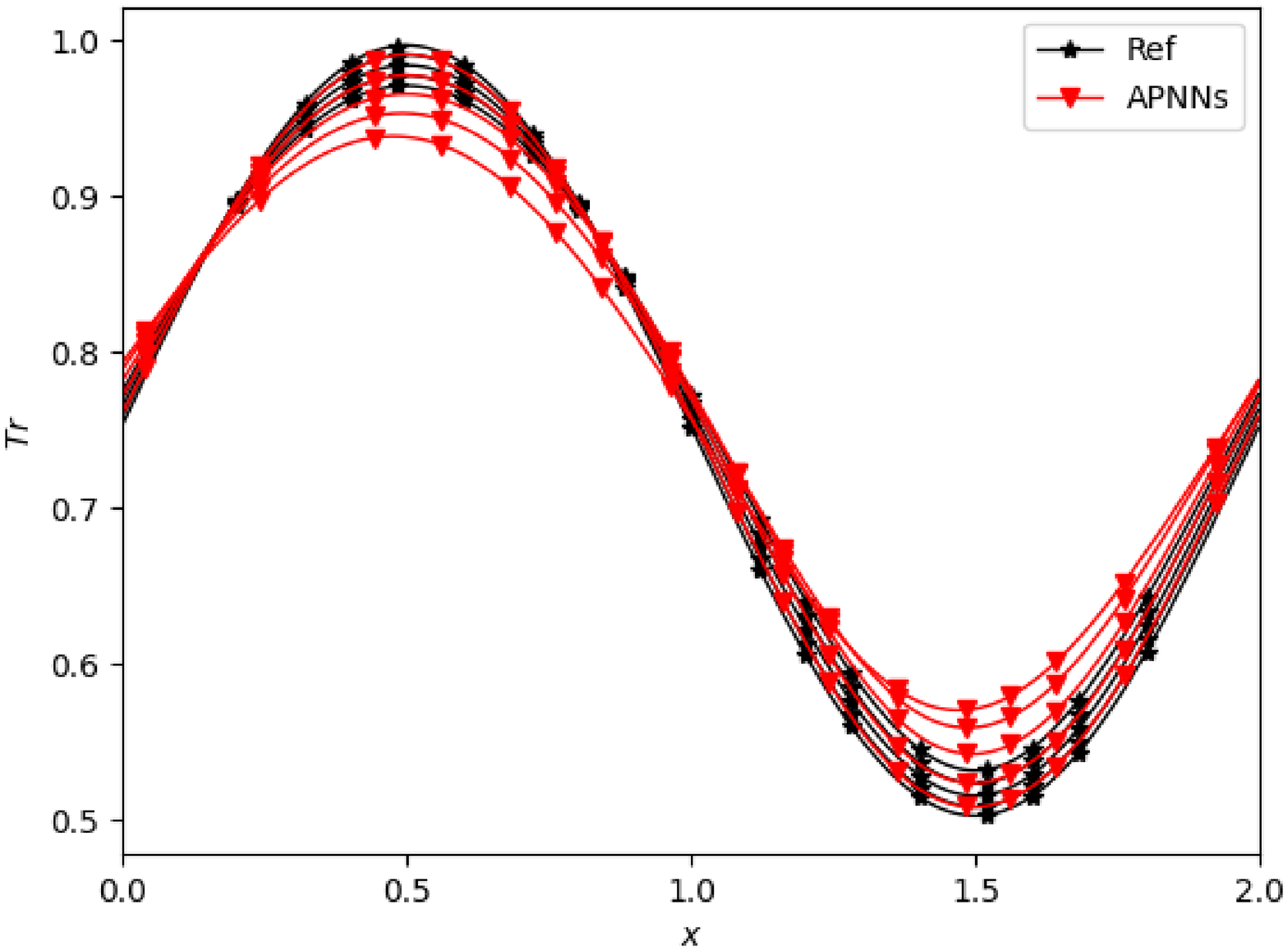}
		\end{minipage}
	}
	
	\caption{Kinetic regime with $\epsilon=1$. Ref v.s. APNNs. (Left) The material temperature $T_e$ at $x=0.0025$. (Right) The radiation temperature $T_r$ at times $t=0.2,\ 0.4,\ 0.6,\ 0.8,\ 1.0$. $\lambda_{2}(1)=20$, $\lambda_{2}(3)=10$, $\lambda_0=0$.}
	\label{figure5.11}
\end{figure}

\begin{figure}[!htbp]
	\setlength{\abovecaptionskip}{0.cm}
	\setlength{\belowcaptionskip}{-0.cm}
	\centering
	{
		\begin{minipage}{2.2in}
			\centering
			\includegraphics[width=2.2in]{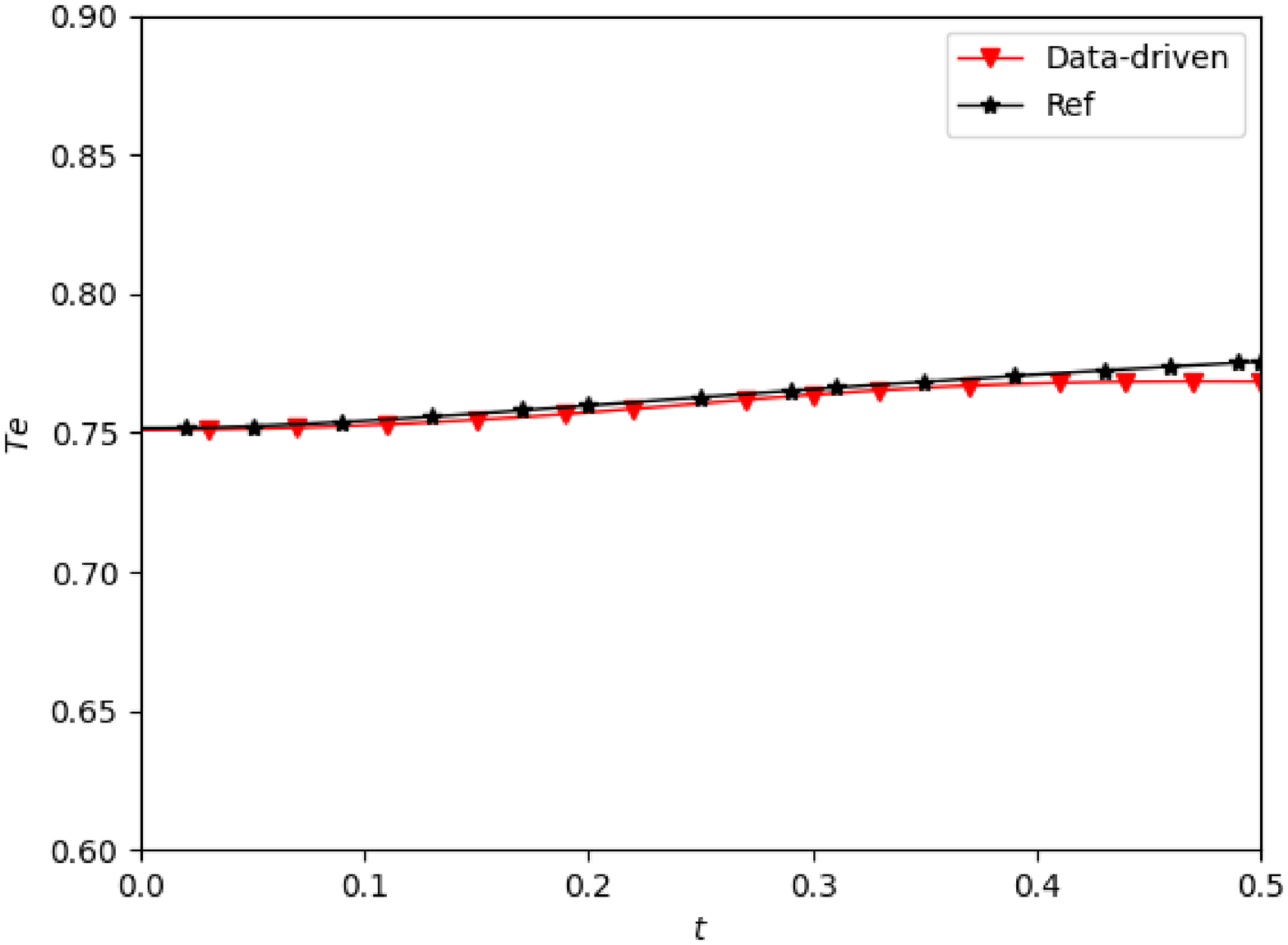}
		\end{minipage}
	}
	{
		\begin{minipage}{2.2in}
			\centering
			\includegraphics[width=2.2in]{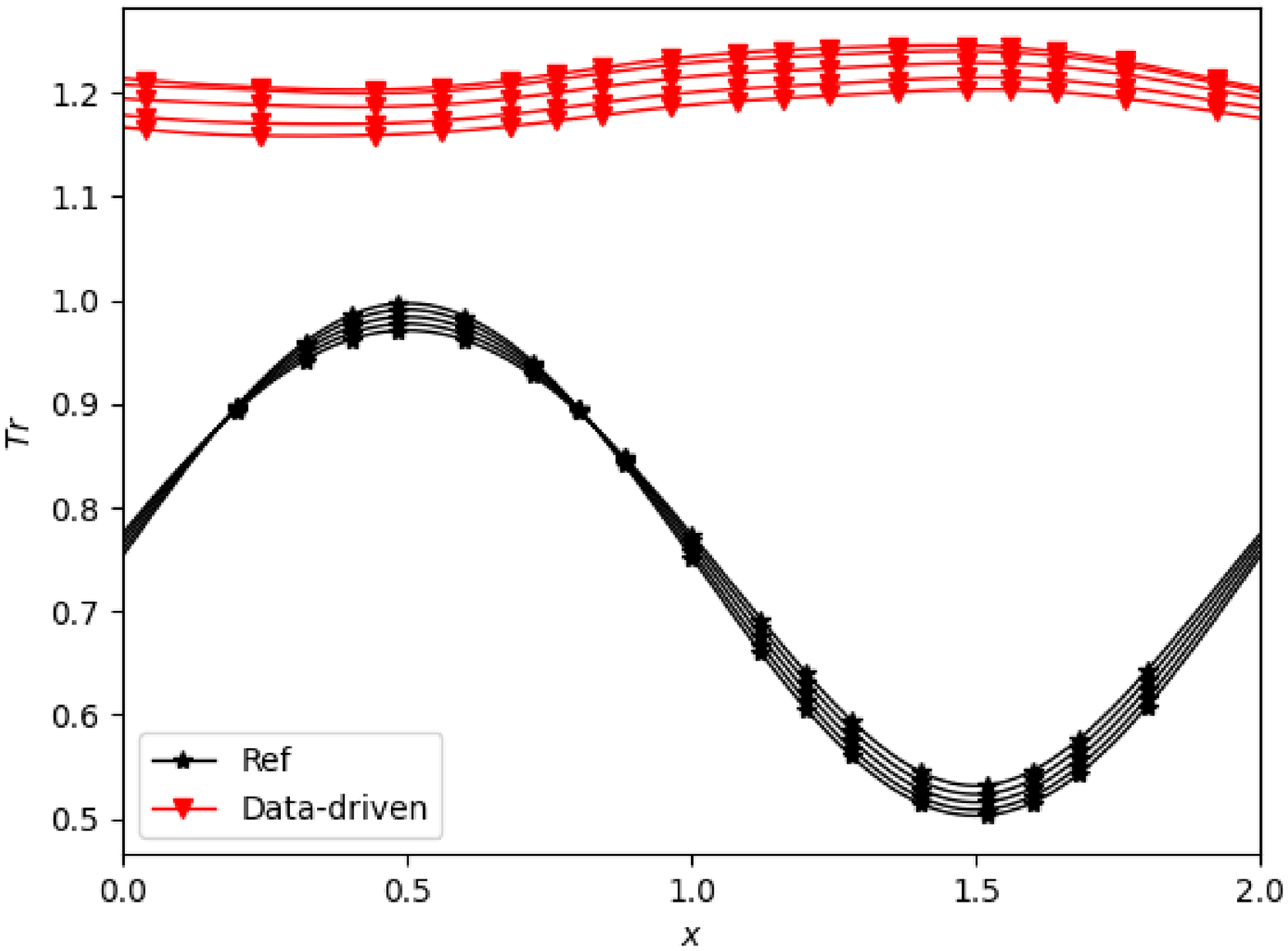}
		\end{minipage}
	}
	
	\caption{Kinetic regime with $\epsilon=1$. Ref v.s. Data-driven. (Left)$T_e$ at $x=0.0025$. (Right)$T_r$ at times $t=0.2,\ 0.4,\ 0.6,\ 0.8,\ 1.0$. $\lambda_{0}=1$ and other weight hyperparameters are 0. }
	\label{figure5.12}
\end{figure}

\begin{figure}[!htbp]
	\setlength{\abovecaptionskip}{0.cm}
	\setlength{\belowcaptionskip}{-0.cm}
	\centering
	{
		\begin{minipage}{2.2in}
			\centering
			\includegraphics[width=2.2in]{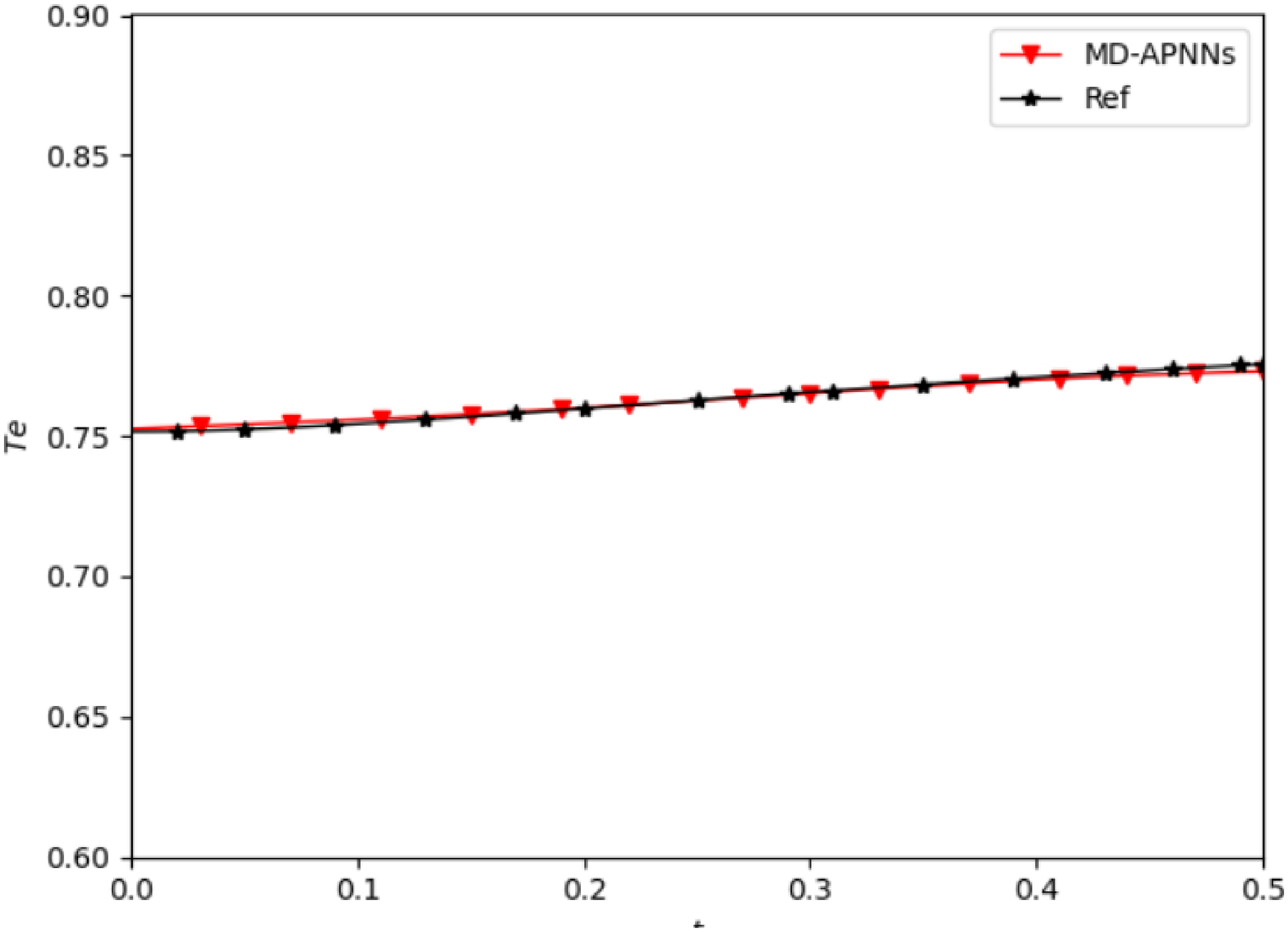}
		\end{minipage}
	}
	{
		\begin{minipage}{2.2in}
			\centering
			\includegraphics[width=2.2in]{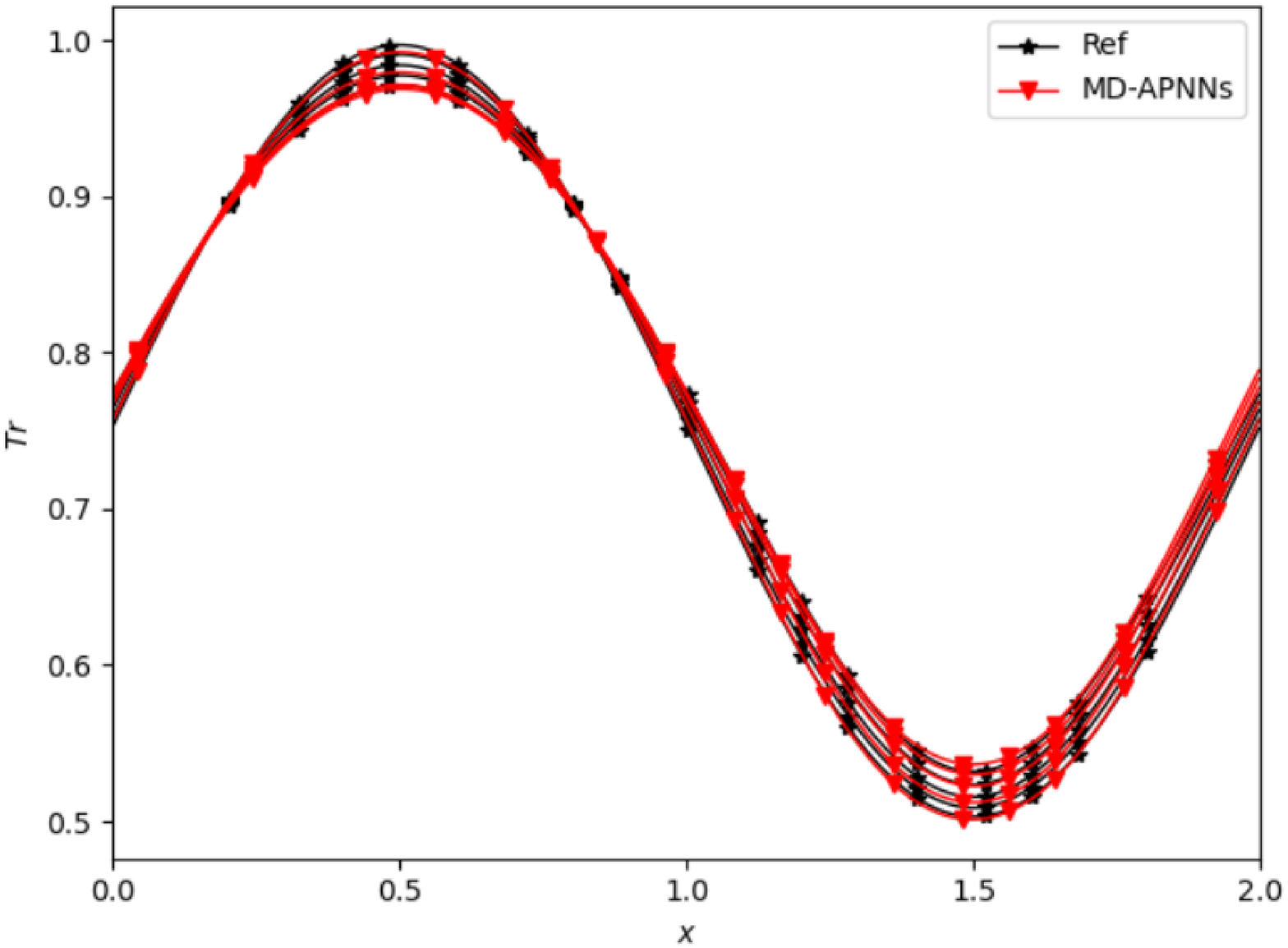}
		\end{minipage}
	}
	
	\caption{Kinetic regime with $\epsilon=1$. Ref v.s. MD-APNNs. (Left)$T_e$ at $x=0.0025$. (Right)$T_r$ at times $t=0.2,\ 0.4,\ 0.6,\ 0.8,\ 1.0$. $\lambda_{2}(1)=20$, $\lambda_{2}(3)=10$.}
	\label{figure5.13}
	
\end{figure}

{\bf (1) Kinetic regime with $\epsilon=1$.} We plot the prediction obtained by APNNs, Data-driven networks and MD-APNNs, together with the reference solution in Figures~\ref{figure5.11}, \ref{figure5.12} and \ref{figure5.13}, respectively.
For MD-APNNs and Data-driven networks, we use $N_0=63$ labeled data about the material temperature $T$.
We observe that APNNs perform well in the approximation of $T_e$ (Figure~\ref{figure5.11} (left)), but do a little poorly in the prediction of $T_r$ (Figure~\ref{figure5.11} (right)).
The Data-driven network method using only the labeled data of $T$ can generate a good approximation of $T_e$ (Figure~\ref{figure5.12} (left)), but the prediction of $T_r$ significantly deviate from the reference solution. (Figure~\ref{figure5.12} (right)).
The MD-APNN method gives consistent results of both $T_e$ and $T_r$ compared with the reference solutions.
The $L^2$-norm errors of APNNs, Data-driven networks, and MD-APNNs are presented in Table~\ref{Table5.5}.

\begin{table}[!htbp]
	\setlength{\abovecaptionskip}{0.cm}
	\setlength{\belowcaptionskip}{-0.cm}
	\centering
	\caption{Kinetic regime with $\epsilon=1$: the errors of $T_e$ and $T_r$ (at $t=0.2$, $0.4$, $0.6$, $0.8$, $1.0$) for APNNs, Data-driven networks and MD-APNNs.}
	\label{Table5.5}
	\setlength{\tabcolsep}{1.2mm}{
		\scalebox{0.97}{
	\begin{tabular}{lllllll}
		\hline\noalign{\smallskip}
		$L^2$ error &$T_e$ & $T_r(t=0.2)$ & $T_r(t=0.4)$ & $T_r(t=0.6)$ & $T_r(t=0.8)$ & $T_r(t=1.0)$ \\
		\noalign{\smallskip}\hline\noalign{\smallskip}
		Data-driven &4.08e-03 &6.59e-01 &6.46e-01  &6.24e-01  &5.99e-01  &5.78e-01    \\
		\noalign{\smallskip}\hline\noalign{\smallskip}
		APNNs &3.01e-03    &7.68e-03  & 1.65e-02  &2.41e-02  &2.91e-02  &3.44e-02   \\
		\noalign{\smallskip}\hline\noalign{\smallskip}
		MD-APNNs & 1.63e-03   &3.74e-03  &8.95e-03  &1.10e-02  &7.90e-03  &4.51e-03   \\
		\noalign{\smallskip}\hline
	\end{tabular}}}
	
\end{table}

\begin{figure}[!htbp]
	\setlength{\abovecaptionskip}{0.cm}
	\setlength{\belowcaptionskip}{-0.cm}
	\centering

	{
		\begin{minipage}{2.2in}
			\centering
			\includegraphics[width=2.2in]{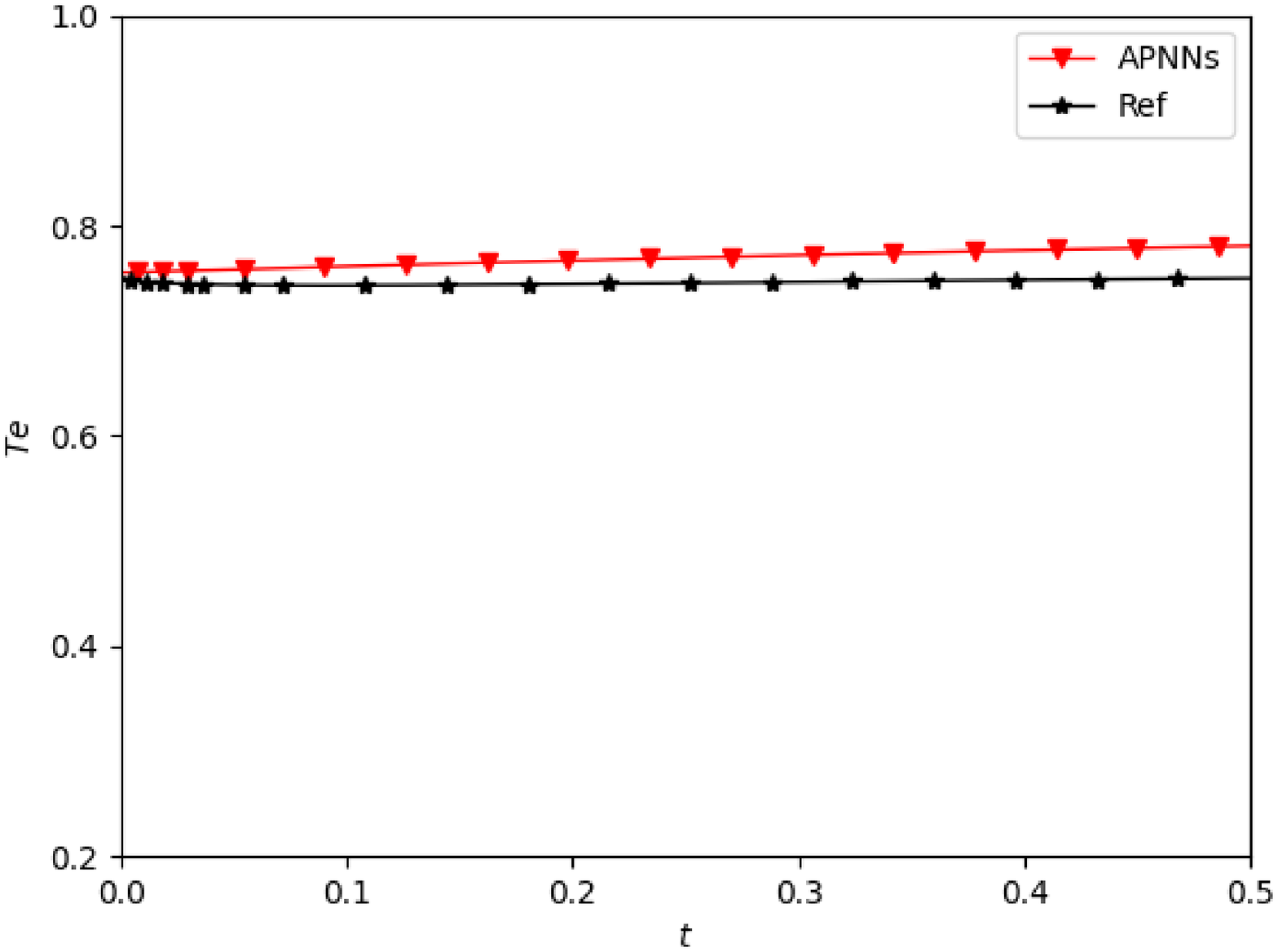}
		\end{minipage}
	}
	{
		\begin{minipage}{2.2in}
			\centering
			\includegraphics[width=2.2in]{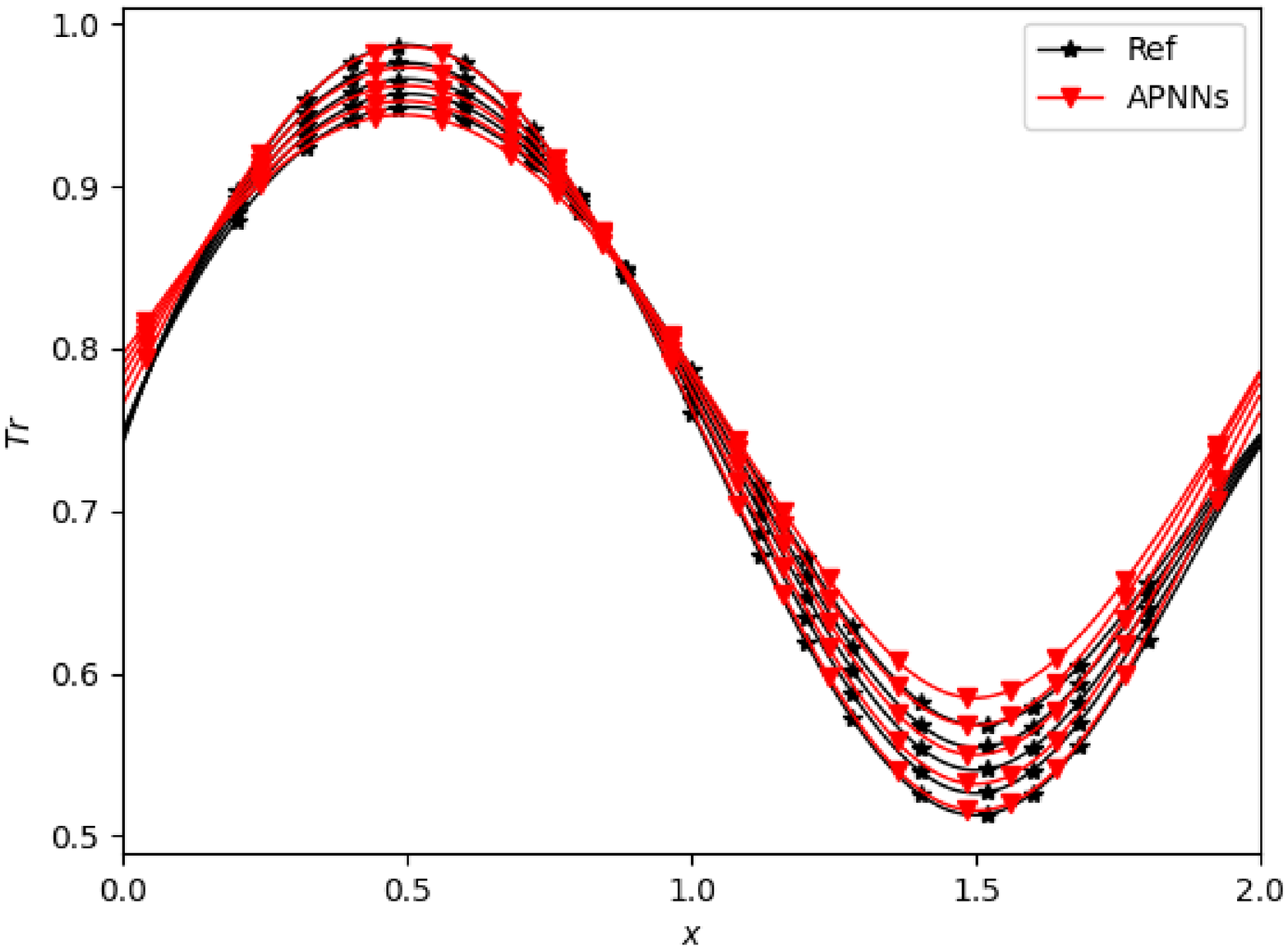}
		\end{minipage}
	}
	
	\caption{Diffusion regime with $\epsilon=10^{-3}$. Ref v.s. APNNs. (Left) $T_e$ at $x=0.0025$. (Right) $T_r$ at times $t=0.2,\ 0.4,\ 0.6,\ 0.8,\ 1.0$. $\lambda_{2}(1)=20$, $\lambda_{2}(3)=10$, $\lambda_0=0$.}
	\label{figure5.14}
\end{figure}

\begin{figure}[!htbp]
	\setlength{\abovecaptionskip}{0.cm}
	\setlength{\belowcaptionskip}{-0.cm}
	\centering
	{
		\begin{minipage}{2.2in}
			\centering
			\includegraphics[width=2.2in]{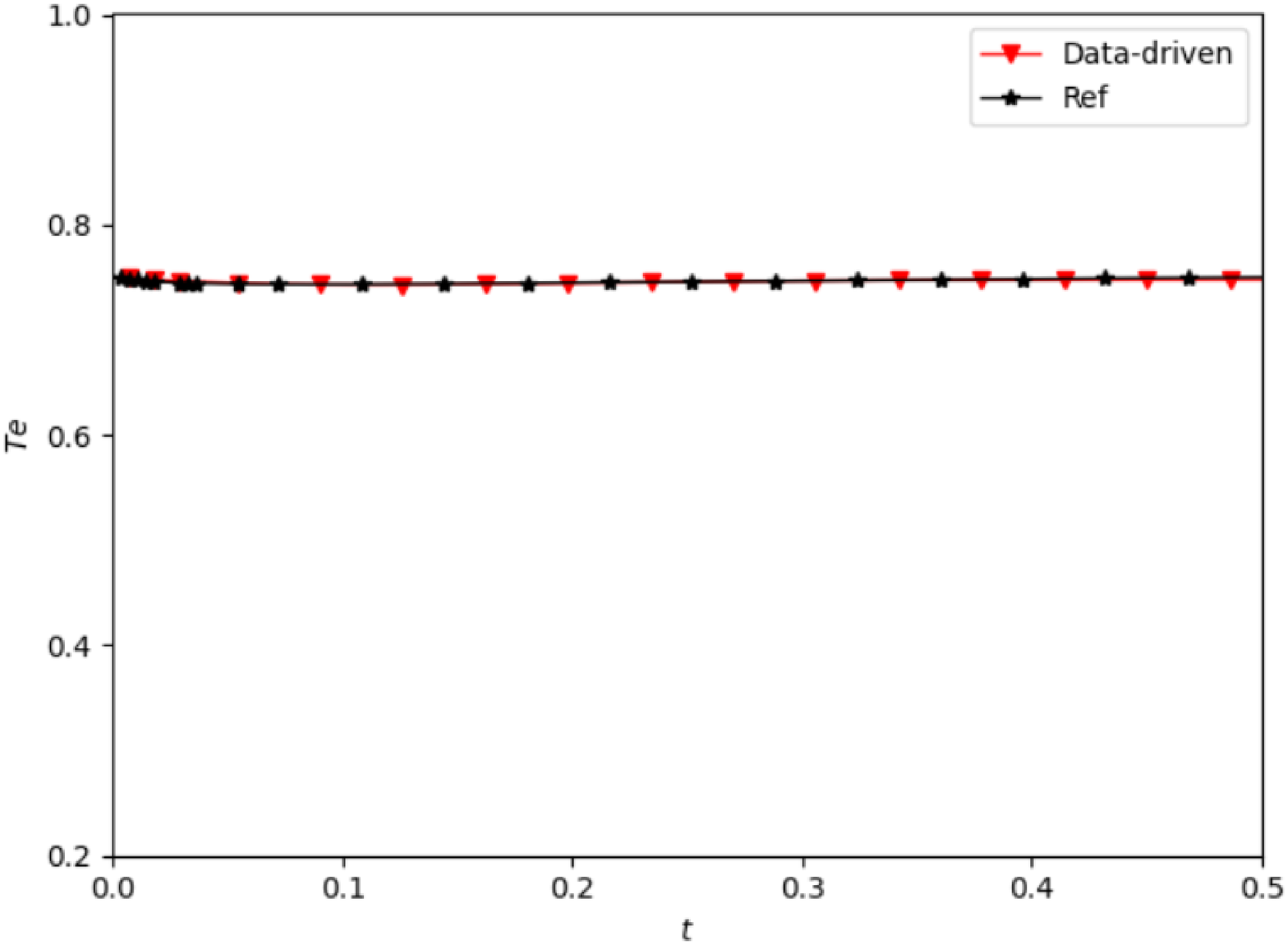}
		\end{minipage}
	}
	{
		\begin{minipage}{2.2in}
			\centering
			\includegraphics[width=2.2in]{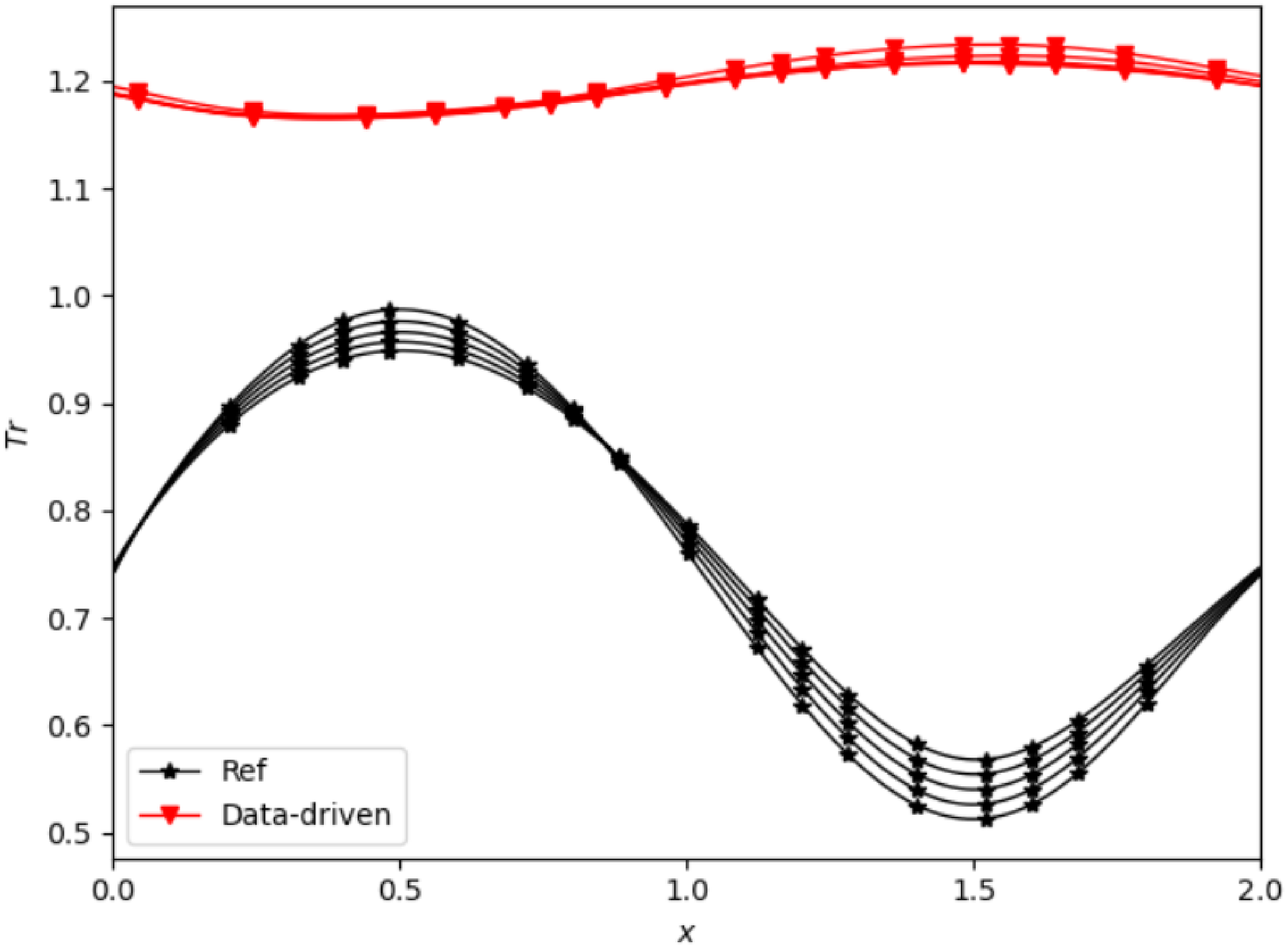}
		\end{minipage}
	}
	
	\caption{Diffusion regime with $\epsilon=10^{-3}$. Ref v.s. Data-driven. (Left) $T_e$ at $x=0.0025$. (Right) $T_r$ at times $t=0.2,\ 0.4,\ 0.6,\ 0.8,\ 1.0$. $\lambda_{0}=1$ and other weight hyperparameters are 0.}
	\label{figure5.15}
\end{figure}

\begin{figure}[!htbp]
	\setlength{\abovecaptionskip}{0.cm}
	\setlength{\belowcaptionskip}{-0.cm}
	\centering
	{
		\begin{minipage}{2.2in}
			\centering
			\includegraphics[width=2.2in]{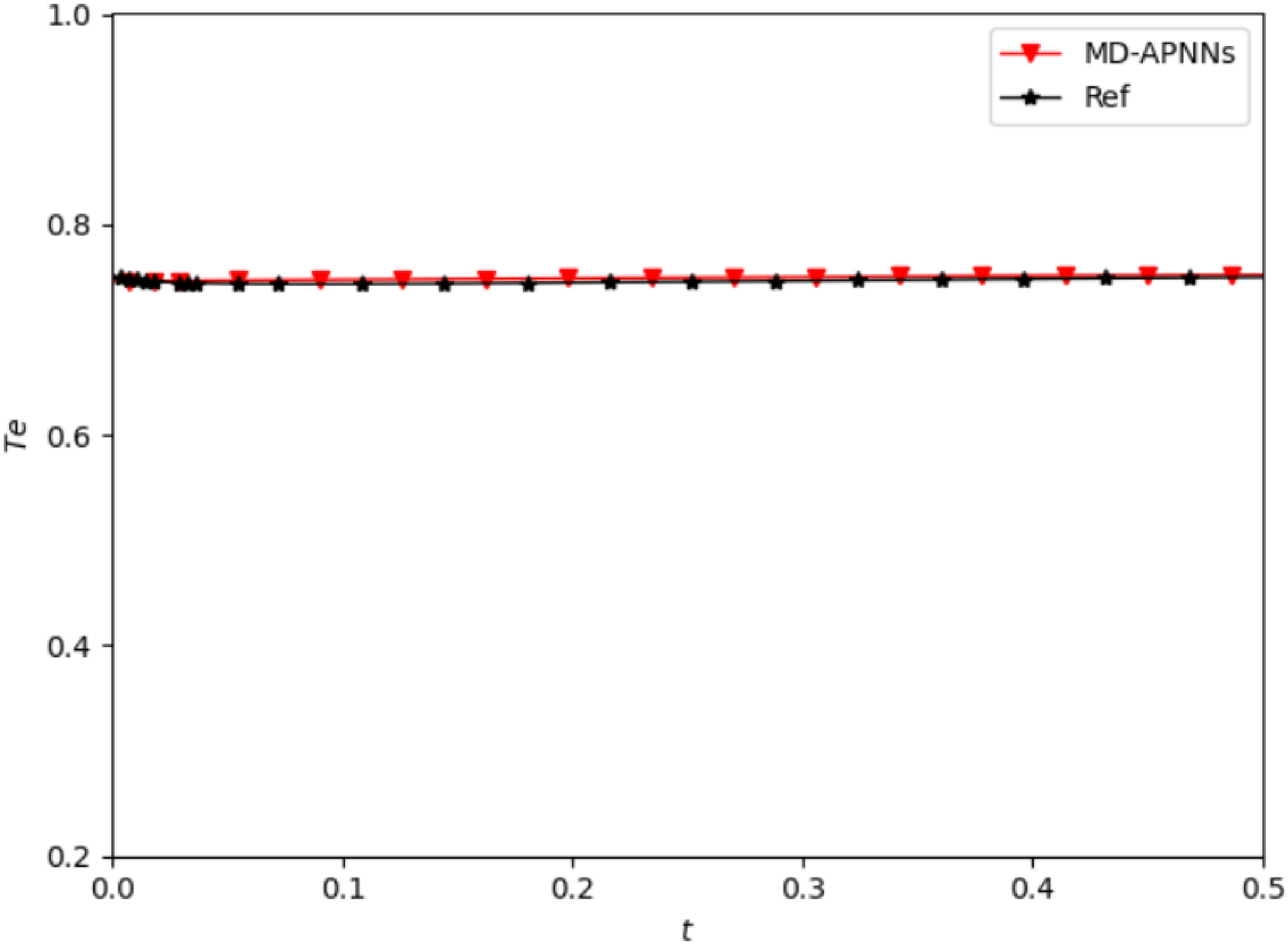}
		\end{minipage}
	}
	{
		\begin{minipage}{2.2in}
			\centering
			\includegraphics[width=2.2in]{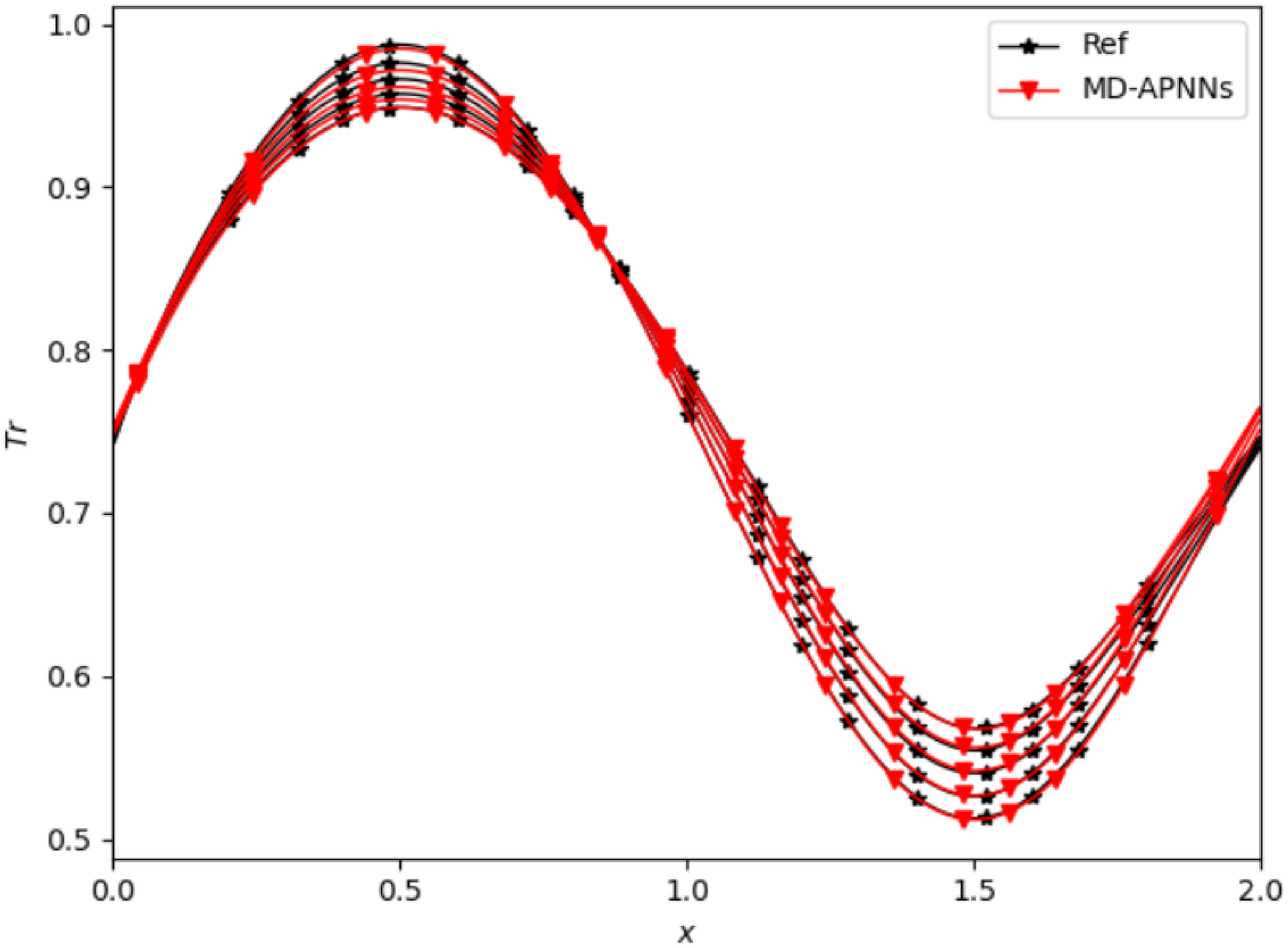}
		\end{minipage}
	}
	
	\caption{Diffusion regime with $\epsilon=10^{-3}$. Ref v.s. MD-APNNs. (Left) $T_e$ at $x=0.0025$. (Right) $T_r$ at times $t=0.2,\ 0.4,\ 0.6,\ 0.8,\ 1.0$. $\lambda_{2}(1)=20$, $\lambda_{2}(3)=10$.}
	\label{figure5.16}
\end{figure}

%
{\bf (2) Diffusion regime with $\epsilon=10^{-3}$.}
At this time, we use $N_0=70$ labeled data of the material temperature $T$ for Data-driven networks and MD-APNNs.
The prediction obtained by APNNs, Data-driven networks and MD-APNNs are plotted in Figures~\ref{figure5.14}, \ref{figure5.15} and \ref{figure5.16}, respectively.
The $L^2$-norm errors of APNNs, Data-driven networks and MD-APNNs are presented in Table~\ref{Table5.6}, which show that the MD-APNN method gives a better prediction compared to the APNN and Data-driven network method.

\begin{table}[!htbp]
	\setlength{\abovecaptionskip}{0.cm}
	\setlength{\belowcaptionskip}{-0.cm}
	\centering
	\caption{Diffusion regime with $\epsilon=10^{-3}$: the errors of $T_e$ and $T_r$ (at $t=0.2$, $0.4$, $0.6$, $0.8$, $1.0$) for APNNs, Data-driven and MD-APNNs.}
	\label{Table5.6}
	\setlength{\tabcolsep}{1.2mm}{
		\scalebox{0.95}{
	\begin{tabular}{lllllll}
		\hline\noalign{\smallskip}
		$L^2$ error &$T_e$ & $T_r(t=0.2)$ & $T_r(t=0.4)$ & $T_r(t=0.6)$ & $T_r(t=0.8)$ & $T_r(t=1.0)$ \\
		\noalign{\smallskip}\hline\noalign{\smallskip}
		Data-driven &1.25e-03 &6.28e-01 &6.08e-01  &5.95e-01  &5.86e-01  &5.79e-01    \\
		\noalign{\smallskip}\hline\noalign{\smallskip}
		APNNs & 3.16e-02    &6.60e-03  &1.11e-02  &1.48e-02  &1.76e-02  &1.94e-02   \\
		\noalign{\smallskip}\hline\noalign{\smallskip}
		MD-APNNs &4.74e-03   &2.94e-03  &4.49e-03  &5.45e-03  &5.15e-03  &4.06e-03   \\
		\noalign{\smallskip}\hline
	\end{tabular}}}
\end{table}

\section{Conclusion}\label{sec6}
In this work, we propose model-data asymptotic-preserving neural netwroks(MD-APNNs) based on macro-micro decomposition for solving the multi-scale nonlinear gray radiative transfer equations(GRTEs). Theoretically, we give the error estimation between neural network predicted solution and the reference solution, which shows that as the loss function vanishes, the approximated solution converges to the analytic solution. Numerically, in order to validate the advantages of the MD-APNN method, we first test linear transfer equations under diffusive scaling and then simulate some steady and unstationary GRTEs under both kinetic regime and diffusive regime. The results show the scheme can capture the quantities of interests accurately. So far we only consider the one-dimensional case with inflow isotropic boundary conditions, periodic boundary conditions and reflecting boundary conditions. In future work, we will focus on solving high-dimensional GRTEs and cases with boundary layer and temperature-dependent opacity.

\begin{appendices}
	\section{Proof of Theorem~\ref{th:error}}
	Set $\widetilde{\rho}=\rho -\rho^*, \widetilde{g}=g - g^*, \widetilde{T}=T - T^*$ and $\widetilde{I}=I - I^*=\widetilde{\rho}+\frac{\epsilon}{\sqrt{\sigma_0}}\widetilde{g}$. For the Dirichlet boundary case, $(\widetilde{\rho},\widetilde{g},\widetilde{T})$ satisfies
	\begin{subnumcases}{\label{eq:A.1}}
	\frac{1}{c}\partial_t\widetilde{\rho}+\frac{1}{\sqrt{\sigma_0}}\left \langle \Omega \cdot \nabla_r \widetilde{g} \right \rangle + \frac{1}{4\pi}C_v \partial_t \widetilde{T}=r_1(t,r) & $\text{in} \ \tau \times D$, \label{eq:A.1-a}\\
	\frac{\epsilon^2}{c}\partial_t\widetilde{g} + \epsilon \Omega \cdot \nabla_r \widetilde{g} - \epsilon \left \langle \Omega \cdot \nabla_r \widetilde{g} \right \rangle + \sqrt{\sigma_0} \Omega \cdot \nabla_r \widetilde{\rho} \notag \\
	+ \sigma \widetilde{g}=r_2(t,r,\Omega) & $\text{in} \ \tau \times D \times S^2$, \label{eq:A.1-b}\\
	\epsilon^2 C_v \partial_t \widetilde{T}- 4\pi \sigma \widetilde{\rho} + \sigma a c \left(T^4 - {T^*}^4\right)=r_3(t,r) & $\text{in} \ \tau \times D$, \label{eq:A.1-c}\\
	\left \langle \widetilde{g} \right \rangle = r_4(t,r) & $\text{in} \ \tau \times D$, \label{eq:A.1-d}\\
	\widetilde{I}=\widetilde{\rho}+\frac{\epsilon}{\sqrt{\sigma_0}}\widetilde{g}=r_5(t,r,\Omega) & $\text{on} \ \Gamma_{-}$,  \label{eq:A.1-e}\\
	\widetilde{I}=\widetilde{\rho}+\frac{\epsilon}{\sqrt{\sigma_0}}\widetilde{g}=r_6(r,\Omega)  & $\text{in} \ \left\{0\right\} \times D\times S^2$, \label{eq:A.1-f} \\
	\widetilde{T}=r_7(r) & $\text{in} \ \left\{0\right\} \times D$, \label{eq:A.1-g}
	\end{subnumcases}
	where $L_{\text{APNNs}}^{\epsilon}=\sum\limits_{i=1}^{7} \Vert r_i \Vert_{L^2}^2$.
	
	For the periodic B.C., Eq.\eqref{eq:A.1-e} is replaced by 
	\begin{align*}
	\widetilde{I}(t,r,\Omega)-\widetilde{I}(t,S(r),\Omega) &:= \left[\widetilde{\rho}(t,r)-\widetilde{\rho}(t,S(r))\right] +
	\left[\widetilde{g}(t,r,\Omega)-\widetilde{g}(t,S(r),\Omega)\right]\\
	&=r_{51}(t,r)+r_{52}(t,r,\Omega)=r_5(t,r,\Omega) \quad \text{on} \ \Gamma.
	\end{align*}
	
	For the reflecting B.C., we substitute Eq.\eqref{eq:A.1-e} by
	\begin{align*}
	\widetilde{I}(t,r,\Omega)-\widetilde{I}(t,r,\Omega^{'}):=\widetilde{g}(t,r,\Omega)-\widetilde{g}(t,r,\Omega^{'})=r_5(t,r,\Omega) \quad \text{on} \ \Gamma_{-}.
	\end{align*}
	
	Substituting $C_v \partial_t \widetilde{T}$ in Eq.\eqref{eq:A.1-a} with Eq.\eqref{eq:A.1-c}, multiplying Eq.\eqref{eq:A.1-b} by $\epsilon/\sqrt{\sigma_0}$, and adding Eq.\eqref{eq:A.1-a} together, we obtain the $(\widetilde{I},\widetilde{T})$ system
	
	\begin{subnumcases}{\label{eq:A.2}}
	\frac{\epsilon^2}{c}\partial_t\widetilde{I}+\sigma\widetilde{I}+\epsilon\Omega\cdot \nabla_r\widetilde{I} + \frac{1}{4\pi}r_3(t,r)-\frac{\sigma ac}{4\pi}\left(T^4-{T^*}^4\right)= \notag\\
	\qquad \quad\qquad \quad\epsilon^2 r_1(t,r)+\frac{\epsilon}{\sqrt{\sigma_0}}r_2(t,r,\Omega), \qquad \quad  \label{eq:A.2-a}\\
	\epsilon^2C_v\partial_t \widetilde{T} - 4\pi\sigma\widetilde{\rho}+\sigma a c \left(T^4-{T^*}^4\right) = r_3(t,r), \label{eq:A.2-b}
	\end{subnumcases}
	
	Testing Eq.\eqref{eq:A.2-a} by $\widetilde{I}$ and integrating it over $D \times S^2$, we obtain
	\begin{align}
	&\frac{\epsilon^2}{2c}\frac{d}{dt}\int_{D}\int_{S^2}\widetilde{I}^2\text{d}\Omega\text{d}r + \int_{D}\int_{S^2}\sigma \widetilde{I}^2\text{d}\Omega\text{d}r +\frac{\epsilon}{2}\int_{\Gamma_{+}}\left|\Omega\cdot n_r\right|\widetilde{I}^2 \text{d}s  \notag \\
	&=\epsilon^2\int_{D}\int_{S^2}r_1\widetilde{I}\text{d}\Omega\text{d}r + \frac{\epsilon}{\sqrt{\sigma_0}}\int_{D}\int_{S^2}r_2\widetilde{I}\text{d}\Omega\text{d}r + \frac{\epsilon}{2}\int_{\Gamma_{-}}\left|\Omega\cdot n_r\right|r_5^2 \text{d}s \notag \\
	&\qquad + \frac{ac}{4\pi}\int_{D}\int_{S^2}\sigma \left(T^4-{T^*}^4\right) \widetilde{I}\text{d}\Omega\text{d}r - \frac{1}{4\pi}\int_{D}\int_{S^2}r_3\widetilde{I}\text{d}\Omega\text{d}r, \label{A.3}
	\end{align}
	where we have applied the integral by parts
	\begin{align*}
	\int_{D}\int_{S^2}\Omega \cdot \nabla_r \widetilde{I} \widetilde{I} \text{d}\Omega \text{d}r &= \frac{1}{2}\int_{D}\int_{S^2} \nabla\cdot(\Omega \widetilde{I}^2) \text{d}\Omega \text{d}r \\
	&=\frac{1}{2}\int_{\Gamma_{+}}\left|\Omega\cdot n_r\right|\widetilde{I}^2 \text{d}s - \frac{1}{2}\int_{\Gamma_{-}}\left|\Omega\cdot n_r\right|\widetilde{I}^2 \text{d}s
	\end{align*}
	
	Due to the fact that $\int_{\Gamma_{+}}\left|\Omega\cdot n_r\right|\widetilde{I}^2 \text{d}s \geq 0$ and $\int_{\Gamma_{-}}\left|\Omega\cdot n_r\right|\widetilde{r_5}^2 \text{d}s \le \int_{\Gamma_{-}}\widetilde{r_5}^2 \text{d}s $, we have
	\begin{align}
	&\frac{\epsilon^2}{2c}\frac{d}{dt}\Vert\widetilde{I}\Vert_{L^2(D\times S^2)}^2 + \int_{D}\int_{S^2}\sigma \widetilde{I}^2\text{d}\Omega\text{d}r \le \epsilon^2\left( r_1,\widetilde{I}\right)_{L^2(D\times S^2)} + \frac{\epsilon}{\sqrt{\sigma_0}} \left(r_2,\widetilde{I}\right)_{L^2(D\times S^2)} \notag \\
	& \qquad \ \ +\frac{1}{4\pi}\left(r_3,\widetilde{I}\right)_{L^2(D\times S^2)}+\frac{\epsilon}{2}\Vert r_5\Vert_{L^2(\Gamma_{-})}^2 + \frac{ac}{4\pi}\int_{D}\int_{S^2}\sigma \widetilde{I} \widetilde{T} G(T,T^*) \text{d}\Omega\text{d}r, \label{A.4}
	\end{align}
	where $(\cdot,\cdot)_{L^2(w)}$ represents the inner product of $L^2(w)$, and $G(T,T^*)= T^3+{T^*}^3+{T^*}^2T+T^{*}T^2$. According to our assumption (\ref{sigmaT-range-b}), $0\textless 4T_{\min}^3\le G(T,T^*)\le 4T_{\max}^3 \textless \infty$, together with (\ref{sigmaT-range-a}) we get
	\begin{align}
	&\frac{\epsilon^2}{2c}\frac{d}{dt}\Vert\widetilde{I}\Vert_{L^2(D\times S^2)}^2 + \sigma_{\min}\Vert\widetilde{I}\Vert_{L^2(D\times S^2)}^2 \le \epsilon^2 \left(r_1,\widetilde{I}\right)_{L^2(D\times S^2)} + \frac{\epsilon}{\sqrt{\sigma_0}} \left(r_2,\widetilde{I}\right)_{L^2(D\times S^2)} \notag \\
	&\qquad \ \ +\frac{1}{4\pi}\left(r_3,\widetilde{I}\right)_{L^2(D\times S^2)}+\frac{\epsilon}{2}\Vert r_5\Vert_{L^2(\Gamma_{-})}^2 + \frac{ac\sigma_{\max}T_{\max}^3}{\pi}\left(\widetilde{T}, \widetilde{I}\right)_{L^2(D\times S^2)}.
	\label{A.5}
	\end{align}	
	Using the Young's inequality with parameter $\alpha_1>0$, we further have
	\begin{align}
	& \frac{\epsilon^2}{2c}\frac{d}{dt}\Vert\widetilde{I}\Vert_{L^2(D\times S^2)}^2 + \sigma_{\min}\Vert\widetilde{I}\Vert_{L^2(D\times S^2)}^2 \le \frac{\pi\epsilon^4}{\alpha_1} \Vert r_1\Vert_{L^2(D)}^2   \notag \\
	& \qquad \ \ +\frac{\epsilon^2}{4\alpha_1\sigma_0}\Vert r_2\Vert_{L^2(D\times S^2)}^2 + \frac{\epsilon}{2}\Vert r_5\Vert_{L^2(\Gamma_{-})}^2 + \frac{\left(ac\sigma_{\max}T_{\max}^3\right)^2}{\alpha_1\pi}\Vert\widetilde{T}\Vert_{L^2(D)}^2   \notag \\
	& \qquad \ \ +\frac{1}{16\alpha_1\pi}\Vert r_3\Vert_{L^2(D)}^2+4\alpha_1\Vert\widetilde{I}\Vert_{L^2(D\times S^2)}^2. \label{A.6}
	\end{align}
	
	The cases periodic and reflecting B.C. are similar. In fact, for periodic B.C., since $n_r = -n_{S(r)}$
	\begin{align*}
	\int_{\Gamma}\Omega\cdot n_r \widetilde{I}^2 \text{d}s&=\frac{1}{2}\int_{\Gamma}\Omega\cdot n_r \left(\widetilde{I}^2(t,r,\Omega)- \widetilde{I}^2(t,S(r),\Omega)\right)\text{d}s \\
	&=\frac{1}{2}\int_{\Gamma}\Omega\cdot n_r r_5 \left(\widetilde{I}(t,r,\Omega)+ \widetilde{I}(t,S(r),\Omega)\right)\text{d}s \\
	&=\frac{1}{2}\int_{\Gamma}\Omega\cdot n_r r_5 \left(2\widetilde{I}(t,r,\Omega)-r_5\right)\text{d}s.
	\end{align*}
	
	And for reflecting B.C., since $\Omega\cdot n_r = -\Omega^{'}\cdot n_r$
	\begin{align*}
	\int_{\Gamma_{-}}\Omega\cdot n_r\widetilde{I}^2\text{d}s &= \int_{\Gamma_{-}}\Omega\cdot n_r\left(\widetilde{I}^2(t,r,\Omega)-\widetilde{I}^2(t,r,\Omega^{'})\right)\text{d}s \\
	&=\int_{\Gamma_{-}}\Omega\cdot n_r\left(\widetilde{I}(t,r,\Omega)-\widetilde{I}(t,r,\Omega^{'})\right)\left(\widetilde{I}(t,r,\Omega)+\widetilde{I}(t,r,\Omega^{'})\right)\text{d}s \\
	&= \int_{\Gamma_{-}}\Omega\cdot n_r r_5\left(2\widetilde{I}^2(t,r,\Omega)-r_5\right)\text{d}s.
	\end{align*}
	
	Furthermore testing Eq.\eqref{eq:A.2-b} by $\widetilde{T}$ and integrating it over $D$ lead to
	\begin{align}
	\frac{\epsilon^2}{2}C_v \frac{d}{dt}\int_{D} \widetilde{T}^2 \text{d}r - 4\pi \int_{D}\sigma \widetilde{\rho}\widetilde{T}\text{d}r + ac\int_{D}\sigma \widetilde{T}\left(T^4-{T^*}^4\right)\text{d}r = \int_{D} r_3 \widetilde{T} \text{d}r, \label{A.7}
	\end{align}
	
	Substituting $\widetilde{\rho}=\left \langle \widetilde{I}\right \rangle-\frac{\epsilon}{\sqrt{\sigma_0}} \left \langle \widetilde{g} \right \rangle=\left \langle \widetilde{I}\right \rangle-\frac{\epsilon}{\sqrt{\sigma_0}}r_4$ into the above equation, we get
	\begin{align}
	& \frac{\epsilon^2}{2}C_v\frac{d}{dt}\Vert\widetilde{T}\Vert_{L^2(D)}^2 + ac \int_{D}\sigma \widetilde{T}^2 G(T,T^*)\text{d}r = \notag \\
	&\qquad \ \ \left(r_3, \widetilde{T}\right)_{L^2(D)} + 4 \pi \int_{D}\sigma \left(\left \langle \widetilde{I}\right \rangle-\frac{\epsilon}{\sqrt{\sigma_0}}r_4\right) \widetilde{T} \text{d}r \le  \notag \\
	&\qquad \ \ \left(r_3, \widetilde{T}\right)_{L^2(D)} + \int_{D}\int_{S^2} \sigma \widetilde{I} \widetilde{T}\text{d}\Omega \text{d}r + \frac{4\pi\epsilon}{\sqrt{\sigma_0}}\int_{D}\sigma r_4 \widetilde{T}\text{d}r.
	\label{A.8}
	\end{align}	
	
	Under Assumption 4.1, applying Young's inequality with parameter $\alpha_2>0$, 
	\begin{align}
	\frac{\epsilon^2}{2}C_v\frac{d}{dt}\Vert\widetilde{T}\Vert_{L^2(D)}^2 +4ac\sigma_{\min}T_{\min}^3 \Vert\widetilde{T}\Vert_{L^2(D)}^2 \le \frac{1}{4\alpha_2}\Vert r_3 \Vert_{L^2(D)}^2 + \frac{\sigma_{\max}^2}{4\alpha_2}\Vert\widetilde{I}\Vert_{L^2(D\times S^2)}^2  \notag \\
	+\frac{1}{4\alpha_2}\left(\frac{4\pi \epsilon \sigma_{\max}}{\sqrt{\sigma_0}}\right)^2 \Vert r_4 \Vert_{L^2(D)}^2 + (4\pi+2)\alpha_2\Vert\widetilde{T}\Vert_{L^2(D)}^2. \label{A.9}
	\end{align}

	Adding Eq.\eqref{A.6} and Eq.\eqref{A.9} together gives us
	\begin{align}
	&\frac{\epsilon^2}{2c}\frac{d}{dt}\Vert\widetilde{I}\Vert_{L^2(D\times S^2)}^2 + \frac{\epsilon^2}{2}C_v\frac{d}{dt}\Vert\widetilde{T}\Vert_{L^2(D)}^2 + \left(\sigma_{\min}-4\alpha_1-\frac{\sigma_{\max}^2}{4\alpha_2}\right) \Vert\widetilde{I}\Vert_{L^2(D\times S^2)}^2  \notag \\ &+\left(4ac\sigma_{\min}T_{\min}^3-\frac{\left(ac\sigma_{\max}T_{\max}^3\right)^2}{\alpha_1\pi}-(4\pi+2)\alpha_2\right) \Vert\widetilde{T}\Vert_{L^2(D)}^2
	\le \frac{\pi\epsilon^4}{\alpha_1} \Vert r_1\Vert_{L^2(D)}^2   \notag \\
	&\qquad \ \ + \frac{\epsilon^2}{4\alpha_1\sigma_0}\Vert r_2\Vert_{L^2(D\times S^2)}^2+\left(\frac{1}{4\alpha_2}+\frac{1}{16\alpha_1\pi}\right)\Vert r_3 \Vert_{L^2(D)}^2 \notag \\
	&\qquad \ \ + \frac{1}{4\alpha_2}\left(\frac{4\pi \epsilon \sigma_{\max}}{\sqrt{\sigma_0}}\right)^2 \Vert r_4 \Vert_{L^2(D)}^2
	+ \frac{\epsilon}{2}\Vert r_5\Vert_{L^2(\Gamma_{-})}^2.  \label{A.10}
	\end{align}
	
	Set the notations
	\begin{align*}
	&A(t):=\frac{1}{2c}\Vert\widetilde{I}\Vert_{L^2(D\times S^2)}^2 + \frac{C_v}{2}\Vert\widetilde{T}\Vert_{L^2(D)}^2, \quad C_1:=\sigma_{\min}-4\alpha_1-\frac{\sigma_{\max}^2}{4\alpha_2}, \\
	&C_2:=4ac\sigma_{\min}T_{\min}^3-\frac{\left(ac\sigma_{\max}T_{\max}^3\right)^2}{\alpha_1\pi}-(4\pi+2)\alpha_2, \quad C_3:=\max\left\{\left|C_1\right|2c, \left|C_2\right|\frac{2}{C_v} \right\} \textgreater 0.
	\end{align*}
	
	Then Eq.\eqref{A.10} becomes
	\begin{align}
	&\epsilon^2 \frac{d}{dt}A(t)-C_3A(t) \le \frac{\pi\epsilon^4}{\alpha_1} \Vert r_1\Vert_{L^2(D)}^2 + \frac{\epsilon^2}{4\alpha_1\sigma_0}\Vert r_2\Vert_{L^2(D\times S^2)}^2 \notag \\
	&\qquad \ \ +\left(\frac{1}{4\alpha_2}+\frac{1}{16\alpha_1\pi}\right)\Vert r_3 \Vert_{L^2(D)}^2 \notag \\
	&\qquad \ \ + \frac{1}{4\alpha_2}\left(\frac{4\pi \epsilon \sigma_{\max}}{\sqrt{\sigma_0}}\right)^2 \Vert r_4 \Vert_{L^2(D)}^2
	+ \frac{\epsilon}{2}\Vert r_5\Vert_{L^2(\Gamma_{-})}^2. \label{A.11}
	\end{align}
	
	Using the Gronwall's inequality, (or multiplying Eq.\eqref{A.11} by $e^{-\frac{C_3t}{\epsilon^2}}$ and integrating it with respect to $t$) we get
	\begin{align}
	e^{-\frac{C_3t}{\epsilon^2}} A(t) - A(0) \le C_4(\epsilon) \int_{0}^{t} e^{-\frac{C_3s}{\epsilon^2}} \left( \Vert r_1\Vert_{L^2(D)}^2 + \Vert r_2\Vert_{L^2(D\times S^2)}^2 + \Vert r_3 \Vert_{L^2(D)}^2  \right. \notag \\
	\left. +\lambda_3\Vert r_4 \Vert_{L^2(D)}^2 + \lambda_1\Vert r_5\Vert_{L^2(\Gamma_{-})}^2 \right) \text{d}s \le C_4(\epsilon) \left(L_{\text{APNNs},g}^{\epsilon} +L_{\text{APNNs},c}^{\epsilon} + L_{\text{APNNs},b}^{\epsilon}\right),
	\label{A.12}
	\end{align}
	where $\lambda_1, \lambda_3 \geq 1$ and 
	\begin{align*}
	C_4(\epsilon)=\max\left\{\frac{\pi\epsilon^2}{\alpha_1}, \frac{1}{4\alpha_1\sigma_0}, \left(\frac{1}{4\alpha_2}+\frac{1}{16\pi\alpha_1}\right)/\epsilon^2, \frac{16\pi^2\sigma_{\max}^2}{4\alpha_2\sigma_0}, \frac{1}{2\epsilon} \right\}>0.
	\end{align*}
	
	Then we have
	\begin{align}
	A(t) &\le e^{\frac{C_3t}{\epsilon^2}} \left(\frac{1}{2c}\Vert\widetilde{I}(0,\cdot,\cdot)\Vert_{L^2(D\times S^2)}^2 + \frac{C_v}{2}\Vert\widetilde{T}(0,\cdot)\Vert_{L^2(D)}^2\right) \notag \\
	& \quad + C_4(\epsilon)e^{\frac{C_3t}{\epsilon^2}} \left(L_{\text{APNNs},g}^{\epsilon} +L_{\text{APNNs},c}^{\epsilon} + L_{\text{APNNs},b}^{\epsilon}\right),
	\label{A.13}
	\end{align}
	
	Set
	\begin{align*}
	C_5(\epsilon):=\max\left\{\frac{1}{2c}e^{\frac{C_3t}{\epsilon^2}},\frac{C_v}{2}e^{\frac{C_3t}{\epsilon^2}}, C_4(\epsilon) e^{\frac{C_3t}{\epsilon^2}} \right\}.
	\end{align*}
	
	Finally we get for all $t\in \tau$,
	\begin{align}
	& \frac{1}{2c}\Vert\widetilde{I}(t,\cdot,\cdot)\Vert_{L^2(D\times S^2)}^2 + \frac{C_v}{2}\Vert\widetilde{T}(t,\cdot)\Vert_{L^2(D)}^2\le C_5(\epsilon) \left(L_{\text{APNNs},g}^{\epsilon} +L_{\text{APNNs},c}^{\epsilon} + L_{\text{APNNs},b}^{\epsilon}\right) \notag \\ &\qquad\quad +C_5(\epsilon) \lambda_2\left(\Vert\widetilde{I}(0,\cdot,\cdot)\Vert_{L^2(D\times S^2)}^2 +\Vert\widetilde{T}(0,\cdot)\Vert_{L^2(D)}^2\right) \notag \\
	&\qquad \le C_5(\epsilon)\left(L_{\text{APNNs},i}^{\epsilon}+L_{\text{APNNs},g}^{\epsilon} + L_{\text{APNNs},c}^{\epsilon} + L_{\text{APNNs},b}^{\epsilon}\right)=C_5(\epsilon)L_{\text{APNNs}}^{\epsilon},
	\label{A.14}
	\end{align}	
	where $\lambda_2 \ge 1$. Taking $L^{\infty}$ norm on Eq.\eqref{A.14} for $\forall t\in \tau$, we conclude that
	\begin{align}
	&\Vert\widetilde{I}\Vert_{L^{\infty}\left(\tau ;L^2(D\times S^2)\right)}^2 \le 2c C_5(\epsilon) L_{\text{APNNs}}^{\epsilon}, \quad 
	\Vert\widetilde{T}\Vert_{L^{\infty}\left(\tau ;L^2(D)\right)}^2 \le \frac{2C_5(\epsilon)}{C_v} L_{\text{APNNs}}^{\epsilon}. \label{A.15}
	\end{align}	
	
\end{appendices}

\section*{Acknowledgements}
This work was supported by the National Key R\&D Program (2020YFA0712200),
	National Key Project (GJXM92579), the Sino-German Science Center (Grant No. GZ 1465)
	and the ISF-NSFC joint research program (Grant No. 11761141008) for S.Jiang, by CAEP foundation (No. CX20200026) and National Key Project (GJXM92579) for W.Sun, by NSFC (No. 12071060) for L.Xu and by NSFC General Projects (No. 12171071) for G.Zhou.

\section*{References}
\bibliography{CMEMA_MDAPNN}
\end{document}